\documentclass[12pt,a4paper,fleqn,english,intoc,bibliography=totoc,index=totoc,BCOR10mm,captions=tableheading,titlepage]{scrbook}
\usepackage{lmodern}

\usepackage[T1]{fontenc}
\usepackage[latin9]{inputenc}
\usepackage{fancyhdr}
\usepackage{babel}
\usepackage{amsmath}
\usepackage{amssymb}
\usepackage[unicode=true,
 bookmarks=true,bookmarksnumbered=true,bookmarksopen=true,bookmarksopenlevel=1,
 breaklinks=false,pdfborder={0 0 0},pdfborderstyle={},backref=false,colorlinks=false]
 {hyperref}
\hypersetup{pdftitle={Your title},
 pdfauthor={Your name},
 pdfpagelayout=OneColumn, pdfnewwindow=true, pdfstartview=XYZ, plainpages=false}

\makeatletter

\AtBeginDocument{\renewcommand{\ref}[1]{\mbox{\autoref{#1}}}}
\newlength{\abc}
\AtBeginDocument{%
\addto\extrasenglish{

}
}

\usepackage[figure]{hypcap}

\let\myTOC\tableofcontents
\renewcommand\tableofcontents{%
  \frontmatter
  \pdfbookmark[1]{\contentsname}{}
  \myTOC
  \mainmatter }

\setkomafont{captionlabel}{\bfseries}
\setcapindent{1em}

\usepackage{calc}

\let\mySection\section\renewcommand{\section}{\suppressfloats[t]\mySection}

\makeatother

\usepackage[notcite]{showkeys}
\usepackage{hyperref}
\begin{document}

\lhead{}

\rhead[\leftmark]{}

\lfoot[\thepage]{}

\cfoot{}

\rfoot{\thepage}
\begin{center}
\textbf{\Large{}Abstract Corrected Iterations}{\Large\par}
\par\end{center}

\begin{center}
{\large{}Haim Horowitz and Saharon Shelah}{\large\par}
\par\end{center}

\begin{center}
\textbf{\small{}Abstract}{\small\par}
\par\end{center}

\begin{center}
{\small{}We consider $(<\lambda)$-support iterations of a version of $(<\lambda)$-strategically
complete $\lambda^+$-c.c. definable forcing notions along partial
orders. We show that such iterations can be corrected to yield an
analog of a result by Judah and Shelah for finite support iterations
of Suslin ccc forcing, namely that if $(\mathbb{P}_{\alpha}, \underset{\sim}{\mathbb Q_{\beta}} : \alpha \leq \delta, \beta <\delta)$
is a FS iteration of Suslin ccc forcing and $U\subseteq \delta$ is
sufficiently closed, then letting $\mathbb{P}_U$ be the iteration
along $U$, we have $\mathbb{P}_U \lessdot \mathbb{P}_{\delta}$.}\footnote{Date: November 10, 2024



Publication 1204 of the second author

Research partially supported by the Israel Science Foundation (ISF) grant no: 1838/19 (2019-2023)

The Israel Science Foundation (ISF) grant no: 2320/23 (2023-2027)

The United States-Israel Binantional Science Foundation (BSF)}{\small\par}
\par\end{center}

\textbf{\large{}0. Introduction}{\large\par}

Our motivation is the following result by Judah and Shelah:

\textbf{Theorem A ({[}JuSh292{]}):} Let $(\mathbb{P}_{\alpha}, \underset{\sim}{\mathbb Q_{\beta}} : \alpha \leq \delta, \beta <\delta)$
be a finite support iteration of Suslin ccc forcing notions (assume
for simplicity that the definitions are without parameters). For a
given $U\subseteq \delta$, let $\mathbb{P}_U$ be the induced iteration
along $U$, then $\mathbb{P}_U \lessdot \mathbb{P}_{\delta}$.

Recent years have witnessed a proliferation of results in generalized
descriptive set theory and set theory of the $\lambda$-reals, and
so an adequate analog of the above-mentioned result for the higher
setting is naturally desirable. Such an analog was crucial for proving
the consistency of $cov(meagre_{\lambda})<\mathfrak{d}_{\lambda}$
in {[}Sh:945{]}. It is not clear that the straightforward analogous
statement holds in the $\lambda$-context, however, it turns out that
the desirable result can be obtained by passing to an appropriate
``correction'' of the original iteration. This was obtained in {[}Sh:1126{]}
for the specific forcing that was relevant for the result in {[}Sh:945{]}.
Our main goal in this paper is to extend the result for a large class
of definable $(<\lambda)$-support iterations of $\lambda^+$-c.c.
forcing. Namely, our mail result will be a more concrete form of the
following:

\textbf{Theorem }(Informal): There is an operation (a ``correction'')
$\mathbb{P} \mapsto \mathbb{P}^{cr}$ on $(<\lambda)$-support iterations
of $(<\lambda)$-strategically complete reasonably definable $\lambda^+$-c.c.
forcing notions along well-founded partial orders, such that $\mathbb{P}^{cr}$
adds the same generics as $\mathbb P$, and if $U$ is an adequate
subset of the set of indices for the iteration, then $\mathbb{P}_U^{cr} \lessdot \mathbb{P}^{cr}$.

Note that even for $\lambda=\aleph_0$ we shall obtain consequences not covered by [JuSh292], as our result includes also iterations with partial memory. Our definability requirements are also much more general than [JuSh292], as instead of analytic definitions we only require that the definitions are reasonably absolute (e.g., in the case of $\lambda=\aleph_0$ and under sufficiently strong large cardinal assumptions, our result covers iterations of forcings defined in $L(\mathbb R)$). The complete formulation of our main result can be found in Conclusions 2.26, 3.12 and 3.13. In order to get a further taste of the main result, we shall illustrate here a less general (but somewhat more formal than before) consequence:
\\
\\
\textbf{Theorem B}: (A) implies (B) where:
\\
\\
A. Let $\lambda$ be a cardinal satisfying $\lambda=\lambda^{<\lambda}$ and let $\bold q$ consist of the following:
\\
\\
a. An ordinal $\alpha(*)$.
\\
\\
b. $\bar u = (u_{\alpha} : \alpha<\alpha(*))$ where $u_{\alpha} \subseteq \alpha$.
\\
\\
c. $\bar{\varphi}=(\varphi_{\alpha} : \alpha<\alpha(*))$ where each $\varphi_{\alpha}$ is a definition of a forcing notion $\mathbb{Q}=\mathbb{Q}_{\varphi_{\alpha}}$ with a generic $\underset{\sim}{\eta_{\alpha}}$, whose members are of the form $p=(tr(p), \bold{B}(...,\underset{\sim}{\eta_{\beta(\epsilon, p)}},...)_{\epsilon<\zeta(p)})$,  where $tr(p)$ is a function from some $v\in [\lambda]^{<\lambda}$ to $H(\lambda)$, $\zeta(p) \leq \lambda$, $\bold B$ is a $\lambda$-Borel function from $(2^{\lambda})^{ \zeta(p)}$ to $H(\lambda)^{\lambda}$, $\beta(\epsilon, p) \in u_{\alpha}$ and $\Vdash_{\mathbb Q} "\underset{\sim}{\eta_{\alpha}}=\cup \{ tr(p) : p\in \underset{\sim}{G} \}"$,
\\
\\
d. If $p\leq_{\mathbb Q_{\varphi_{\alpha}}} q$ then $tr(p) \subseteq tr(q)$.
\\
\\
e. If $\{ p_i : i<j \} \subseteq \mathbb{Q}_{\varphi_{\alpha}}$, $tr(p_i)=\eta$ for all $i<j$, and $j\leq lg(\eta)$, then $\{ p_i : i<j \}$ has a common upper bound that is $\lambda$-Borel computable from $\{p_i : i<j \}$.
\\
\\
f. The forcing notions $\mathbb{Q}_{\varphi_{\alpha}}$ are $(<\lambda)$-strategically complete and satisfy a strengthening of $\lambda^+$-cc called "$(\lambda, D)$-cc" (to be defined later).
\\
\\
g. For each $\mathbb{Q}_{\varphi_{\alpha}}$, the trunks and the generic satisfy a few additional reasonable requirements (to be specified in Definition 1.4).
\\
\\
h. The definitions $\varphi_{\alpha}$ and their relevant properties (e.g. compatibility of conditions, the trunk of a condition being a specific $\eta$, etc) are absolute between models of the form $V^{\mathbb{P}_1}$ and $V^{\mathbb{P}_2}$ where $\mathbb{P}_1 \lessdot \mathbb{P}_2$ are $(<\lambda)$-strategically complete and $\lambda^+$-cc.
\\
\\
B. There is $(\mathbb{P}^{cr}_{\bold q}, \bar{\underset{\sim}{\eta}}^*)=(\mathbb P, \underset{\sim}{\bar \eta}^*)$ where:
\\
\\
a. $\mathbb P$ is $(<\lambda)$-strategically complete and $\lambda^+$-cc.
\\
\\
b. $\underset{\sim}{\bar \eta}^*=(\underset{\sim}{\eta_{\alpha}^*} : \alpha<\alpha(*))$ is a sequence of $\mathbb P$-names of $\lambda$-reals.
\\
\\
c. For each $\alpha<\alpha(*)$, let $V^{\alpha}:=V[...,\underset{\sim}{\eta_{\beta}^*},...]_{\beta \in u_{\alpha}}$, then $\underset{\sim}{\eta_{\alpha}^*}$ is "somewhat generic" for $\mathbb{Q}_{\varphi_{\alpha}}^{V^{\alpha}}$ in the sense that if $I$ is an antichain in $\mathbb{Q}^{V^{\alpha}}_{\varphi_{\alpha}}$ that is absolutely maximal, then $\underset{\sim}{\eta_{\alpha}^*}$ satisfies some $p\in I$. 
\\
\\
d. If $U\subseteq \alpha(*)$ and $\alpha \in U \rightarrow u_{\alpha} \subseteq U$, then $\bold q \restriction U$ is naturally defined and $(\mathbb{P}_{\bold q \restriction U}^{cr}, {\underset{\sim}{\bar{\eta}^*}} \restriction U)$ are as above for $\bold q \restriction U$.
\\
\\
e. If $U_1, U_2 \subseteq \alpha(*)$ are as in (d) and $\pi : U_1 \rightarrow U_2$ is an isomorphism such that $\alpha \in u_{\beta} \leftrightarrow \pi(\alpha) \in u_{\pi (\beta)}$ and such that $\varphi_{\alpha}=\varphi_{\pi (\alpha)}$ for all $\alpha\in U_1$, then there are $\mathbb{P}_l \lessdot \mathbb P$ $(l=1,2)$ such that $\underset{\sim}{\eta_{\alpha}^*}$ is a $\mathbb{P}_l$-name for every $\alpha \in U_l$, $\mathbb{P}_l=\mathbb{P}_{\bold q \restriction U_l}^{cr}$ and $\pi(\underset{\sim}{\eta_{\alpha}^*})=\underset{\sim}{\eta_{\pi(\alpha)}^*}$ for every $\alpha \in U_1$.
\\
\\
We expect our general result to be applicable in numerous contexts. As mentioned above, a specific case was applied in [Sh:945] to obtain the consistency of a new inequality of cardinal invariants for the $\lambda$-reals. We expect also applications to cardinal invariants of the continuum, as indicated by the following immediate corollary:
\\
\\
\textbf{Theorem C}: Let $\frak{x}_1,...,\frak{x}_n$ be cardinal invariants of the continuum such that the consistency of $\aleph_1<\frak{x}_1<...<\frak{x}_n<\frak c$ can be forced over a model of $CH$ using a FS iteration over a well-founded partial order of definable forcing notions satisfying the assumptions of our main theorem, then it is also consistent that $\frak s=\aleph_1 <\frak{x}_1 <...<\frak{x}_n <\frak c$.
\\
\\
The above theorem follows from the proof from [JuSh292] of the fact that FS iterations of Suslin ccc forcing notions over a model of $CH$ preserve $\frak s=\omega_1$. The proof relies on the aforementioned result about subiterations of Suslin ccc forcing, and so it follows for FS iterations over a well-founded partial order of suitable forcing notions by using the corresponding corrected iteration and the main result of this paper.
\\
\\
We shall start by defining our building blocks, namely forcing templates
and iteration templates. These will allow for a much larger variety
of examples than what appears in {[}Sh:1126{]} (in particular, an
iteration may involve forcing notions with different definitions).
One of the differences between the current work and {[}Sh:1126{]}
is that our forcing notions might be definable using parameters that
don't belong to $V$, and so this will require the introduction of
a new type of memory (``weak memory'') that will allow the computation
of the relevant parameters.

We then continue by introducing the class $\mathbf M$ of iteration
parameters, from which we shall practically construct our iterations.
We shall then consider the notion of an existentially closed iteration
parameter, and we shall isolate a property of iteration parameters
that guarantee the existence of an existentially closed erxtension.
We shall then obtain our desired corrected iteration from those existentially
closed extensions by taking an appropriate closure under $\mathbb{L}_{\lambda^+}$.
\\
\\
\textbf{Notation and conventions D}: Throughout the paper, ordinals will be denoted by lowercase Greek letters, with the exceptions of the letters $\kappa, \lambda$, $\mu$ (and sometimes $\theta$ and $\chi$) that will be used for cardinals, and $\varphi$, $\psi$ (and sometimes $\theta$ and $\chi$) which will be used to denote formulas. For regular $\kappa<\lambda$ we denote the set $\{ \delta <\lambda : cf(\delta)=\kappa\}$ by $S_{\kappa}^{\lambda}$. Forcing templates will be denoted by $\bold p$ and iteration templates will be denoted by $\bold q$. Forcing notions will be denoted by $\mathbb P$ and $\mathbb Q$, where typically $\mathbb P$ will be used for iterations and $\mathbb Q$ will be used for iterands. We adhere to the Jerusalem tradition according to which "$p\leq q$" means that the forcing condition $q$ is stronger than $p$. We shall work with the following modification of $H(\kappa)$:
\\
\\
\textbf{Definition E: }A) Given two sets $X$ and $x$, $trcl_X(x)=trcl(x,X)$
will be defined as the minimal set $u$ such that:

1. $x\in u$.

2. $y\subseteq u$ for every $y\in u\setminus X$.

B) For a cardinal $\kappa$ and a set $X$ we define $H_{\leq \kappa}(X)$
as the collection of sets $x$ such that $|trcl(x,X)|\leq \kappa$
and $\emptyset \notin trcl(x,X)$.

C) $X$ is called $\kappa$-flat if $x\notin H_{\leq \kappa}(X\setminus \{x\})$
for every $x\in X$ (we may use $X$ as a set of atoms as in Definition 1.1(B)).
\\
\\
D) Given a cardinal $\lambda$, an ordinal $\zeta<\lambda^+$ and a set $X$, we define $H_{\leq \lambda, \zeta}(X)$ as follows: $H_{\leq \lambda, 0}:=X$, and for $\zeta>0$, letting $H_{\leq \lambda, <\zeta}(X):=\underset{\xi<\zeta}{\cup}H_{\leq \lambda, \xi}(X)$, we define $H_{\leq \lambda ,\zeta}:=[H_{\leq \lambda, <\zeta}(X)]^{<\lambda}$. So $H_{\leq \lambda}(X)=H_{\leq \lambda, <\lambda^+}(X)$.
\\
\\
Throughout the paper, we shall use the notion of $\lambda$-Borel functions. Our definitions will be somewhat nonstandard. Below we provide two possible versions for what is meant by a $\lambda$-Borel function:
\\
\\
\textbf{Nonstandard Definition F}: 
A. We say that $\bold B$ is a $\lambda$-Borel function if:
\\
\\
(Version 1) There are sets $X$ and $Y$ such that:
\\
\\
a. $\bold B$ is a definition of a partial function from $H_{\leq \lambda}(X)$ to $H_{\leq \lambda}(Y)$.
\\
\\
b. If $\mathbb{P}_1 \lessdot \mathbb{P}_2$ are  (relatives of) $(<\lambda)$-strategically complete forcing notions satisfying $\lambda^+$-cc (or $(\lambda, D)$-cc, which will be defined later in the paper), then $\bold{B}^{V^{\mathbb{P}_1}}=\bold{B}^{V^{\mathbb{P}_2}} \restriction V^{\mathbb{P}_1}$.
\\
\\
(Version 2) There are two sets $X$ and $Y$ such that:
\\
\\
a. $\bold B=(\bold{B}_{x, \zeta, y, \xi} : x\in [X]^{\leq \lambda}, y\in [Y]^{\leq \lambda}, \zeta, \xi<\lambda^+)$ where each $\bold{B}_{x, \zeta, y, \xi}$ is the $\lambda$-analog of the  ord-hc Borel operations from [Sh630] (to be defined in Clause (B) below).
\\
\\
b. ($x_1 \subseteq x_2) \wedge (y_1 \subseteq y_2) \wedge (\zeta_1 \leq \zeta_2) \wedge (\xi_1 \leq \xi_2) \rightarrow \bold{B}_{x_1, \zeta_1, y_1, \xi_1} \subseteq \bold{B}_{x_2, \zeta_2, y_2, \xi_2}$.
\\
\\
c. Given $z\in H_{\leq \lambda}(X)$, $\bold{B}(z)=\bold{B}_{x, \zeta, y, \xi}(z)$ whenever RHS is defined.
\\
\\
B (following [Sh630]). We define the $\lambda$-analog of the family of ord-hc Borel operations as the minimal family $\mathcal F$ of functions satisfying the following:
\\
\\
a. Each $\bold B \in \mathcal F$ is a function with $\leq \lambda$ coordinates, where the possible inputs for each coordinate are sets from $H_{\leq \lambda}(X)$ where $|X| \leq \lambda$, ordinals, truth values, sequences of ordinals of length $\leq \lambda$ and sequences of truth values of length $\leq \lambda$.
\\
\\
b. The range of each $\bold B \in \mathcal F$ consists of elements from $H_{\leq \lambda}(Y)$ (for some $Y$ satisfying $|Y| \leq \lambda$), ordinals and truth values.
\\
\\
c. $\mathcal F$ is closed under composition.
\\
\\
d. $\mathcal F$ contains the following atomic functions:
\\
\\
1. $\neg x$ for a truth value $x$.
\\
2. $x_1 \vee x_2$ for truth values $x_1$ and $x_2$.
\\
3. $\underset{i<\alpha}{\wedge} x_i$ for $\alpha \leq \lambda$ and truth values $x_i$.
\\
4. The constant values $True$ and $False$.
\\
5. For all $\alpha \leq \lambda$, $x_{\gamma}$ varying on truth values and for all $y_{\gamma}$ varying on sets from $H_{\leq \lambda}(X)$ (for $\gamma<\lambda$):
\\
- If $x_{\gamma}$ but not $x_{\delta}$ for $\delta < \gamma$ then $y_{\gamma}$.
\\
- If $\neg x_{\gamma}$ for every $\gamma<\alpha$ then $y_{\alpha}$.
\\
6. Similarly for ordinals.
\\
7. $\{ y_i : i<\alpha, x_i=T\}$ where $\alpha \leq \lambda$ and each $y_i$ varies on $H_{\leq \lambda}(X)$-sets or on ordinals, $x_n$ on truth values.
\\
8. The truth value of "$x$ is an ordinal" where $x$ varies on $H_{\leq \lambda}(X)$-sets.
\\
\\
$\textbf{Remark G}$: The reason for the second version of the definition is that for the $\lambda$-analog of the ord-hc Borel operations from [Sh630] we would like to have functions from $H_{\leq \lambda}(X)$ to $H_{\leq \lambda}(Y)$ where $|X|,|Y| \leq \lambda$. But as it might be the case that $|X|, |Y| > \lambda$, the formulation in the second version is required.

\textbf{\large{}1. Preliminary definitions, assumptions and facts}{\large\par}

\textbf{\large{}Forcing templates}{\large\par}

In this section we shall define the templates from which individual
forcing notions in the iteration shall be constructed. As we don't
have a general preservation theorem for $\lambda^+$-c.c. in $(<\lambda)$-support
iterations (see [Sh1036] and history there), we shall use the notion of $(\lambda, D)$-chain condition
for a filter $D$ (to be defined later) for which we have a preservation
result, and so the templates will include an appropriate filter to
witness this. Similarly to {[}Sh:630{]}, the forcing templates will
consist of a model $\mathfrak{B}_{\mathbf p}$ and formulas that will
define the forcing inside it. The forcing will be defined using a
parameter, which shall be a function whose domain is denoted $I_{\mathbf p}^0$.
The generic element will be a function whose domain is the set $I_{\mathbf p}^1$.
Additional formulas will provide winning strategies for strategic
completeness and will provide a compatibility relation on the forcing
that will satisfy the $(\lambda, D)$-chain condition.

\textbf{Hypothesis 0: }Throughout this paper, we assume that:

a. $\lambda$ is a cardinal satisying $\lambda=\lambda^{<\lambda}$ 

b. $D$ is a $\lambda$-complete filter on $\lambda^+ \times \lambda^+$
satisfying the following:

1. $\{(\alpha,\beta) : \alpha<\beta<\lambda^+\} \in D$.

2. If $u_{\alpha} \in [Ord]^{<\lambda}$ $(\alpha<\lambda^+)$, $g:\underset{\alpha<\lambda^+}{\cup}u_{\alpha} \rightarrow D$
and $f_{\alpha} : u_{\alpha} \rightarrow Ord$ has range $\subseteq \lambda$
$(\alpha<\lambda^+)$, then the following set belongs to $D$: $\{(\alpha,\beta) : \alpha<\beta<\lambda^+,$
$(f_{\alpha},f_{\beta})$ is a $\Delta-$system pair (see Definition
1.2 below), $\xi \in u_{\alpha} \cap u_{\beta} \rightarrow (\alpha,\beta) \in g(\xi)\}$.
\\
\\
3. $(\lambda^+ \setminus \gamma) \times (\lambda^+ \setminus \gamma) \in D$ for every $\gamma<\lambda^+$.
\\
\\
The following will serve to define the forcing notions that we intend to iterate:

\textbf{Definition 1.1: }Given a cardinal $\kappa>\lambda$. We call
$\mathbf{p}=(\lambda_{\mathbf p},\kappa_{\mathbf p}, \mathbf{U_p}, \mathbf{I_p}, \mathfrak{B}_{\mathbf p}^0, I_{\mathbf p}^{0}, I_{\mathbf p}^{1},\bar{\varphi},D_{\mathbf p}, \mathfrak{B}_{\mathbf p}, \bold{T_p}, R_{\bold p})$
a $(\lambda,D)$-forcing template if:

A) $\lambda=\lambda_{\mathbf p}<\kappa=\kappa_{\mathbf p}$.

B) $I_{\mathbf p}^{0} \cup I_{\mathbf p}^{1}\subseteq H_{\leq \lambda}(\mathbf{U_p} \cup \mathbf{I_p})$
where $\mathbf{U}=\mathbf{U_p}$ and $\mathbf{I}=\mathbf{I_p}$ are disjoint
sets of atoms. 
\\
$[$Motivation: $I^0_{\bold p}$ will serve as the domain of the "input" for the definition of the forcing, i.e. the parameters used in the definition of the forcing. $I^1_{\bold p}$ will serve as the "output", i.e. the domain of the generic.$]$

C) $\mathfrak{B}_{\mathbf p}$ is the expansion of $(H_{\leq \lambda}(\mathbf{U_p} \cup \mathbf{I_p}),\in)$
by adding the relations $|\mathfrak{B}_{\mathbf p}^{0}|$ and $P^{\mathfrak{B}_{\mathbf p}^{0}}$
for every $P\in \tau(\mathfrak{B}_{\mathbf p}^{0})$ for a model $\mathfrak{B}_{\mathbf p}^{0}$
with universe $\mathbf{I} \cup \mathbf{U}$.
$[$This will be the structure inside of which the definition of the forcing will be interpreted.$]$

D) $\bar{\varphi}=(\varphi_l(\bar{x_l},\bar y) : l<7)$ is a sequence
of first order formulas from $\mathbb{L}(\tau_{\mathfrak{B}_{\mathbf p}})$
and $lg(\bar{x_l})=k_l$ where $k_0=1$, $k_1=2$, $k_2=3$, $k_3=3$,
$k_4=2$, $k_5=2$, $k_6=2$. We allow the $\varphi_i$ to include a second
order symbol $F$ (over which we shall not quantify) that will be
interpreted as a function $h: I^0_{\mathbf p} \rightarrow \lambda$.
$[$These will be the formulas defining the forcing and its relevant features.$]$

E) $D_{\mathbf p}=D$ is a $\lambda$-complete filter as in Hypothesis
0 above.
\\
\\
F) $\bold{T_p}$ is a set that contains all possible trunks for conditions in the forcing, each is a function from some $u\in [I^1_{\bold p}]^{<\lambda}$ to $H(\lambda)$.
\\
\\
G) $R_{\bold p}$ is a reflexive binary relation on $\bold{T_p}$.
\\
\\
H) If $\{ t_{\alpha} : \alpha<\lambda^+ \} \subseteq \bold{T_p}$, then $\{(\alpha, \beta) : \alpha<\beta<\lambda^+, t_{\alpha}R_{\bold p} t_{\beta}\} \in D$.

Remark: We may omit the index $\mathbf{p}$ whenever the identity of
$\mathbf{p}$ is clear from the context.

\textbf{Definition 1.2: }Suppose that $u_l \in [Ord]^{<\lambda}$
$(l=1,2)$. A pair of functions $f_l: u_l \rightarrow Ord$ $(l=1,2)$
is called a $\Delta$-system pair if $otp(u_1)=otp(u_2)$, and for
every $\alpha \in u_1 \cap u_2$, $otp(u_1 \cap \alpha)=otp(u_2 \cap \alpha)$
and $f_1(\alpha)=f_2(\alpha)$.

\textbf{Claim/Example 1.3: }Let $D_{\lambda}^0$ be the collection
of subsets $X\subseteq \lambda^+ \times \lambda^+$ such that for
some club $E\subseteq \lambda^+$ and regressive function $g: S_{\lambda}^{\lambda^+} \rightarrow \lambda^+$,
$\{(\alpha, \beta) : \alpha<\beta<\lambda^+, \alpha \in S_{\lambda}^{\lambda^+} \cap E, \beta \in S_{\lambda}^{\lambda^+} \cap E, g(\alpha)=g(\beta) \}  \subseteq X$,
then $D_{\lambda}^0$ is as required in definition 1.1(E).

\textbf{Proof: }Clearly, $\emptyset \notin D_{\lambda}^0$. Let $(u_{\alpha} : \alpha<\lambda^+)$,
$(f_{\alpha} : \alpha<\lambda^+)$ and $g$ be as in definition 1.1(E),
then for every $\xi \in \underset{\alpha<\lambda^+}{\cup}u_{\alpha}$
there is a club $E_{\xi} \subseteq \lambda^+$ and a regressive function
$h_{\xi}: S_{\lambda}^{\lambda^+} \rightarrow \lambda^+$ such that
$X_{\xi} \subseteq g(\xi)$ where: $X_{\xi}:=\{(\alpha,\beta) : \alpha<\beta<\lambda^+, \alpha \in S_{\lambda}^{\lambda^+} \cap E_{\xi}, \beta \in S_{\lambda}^{\lambda^+} \cap E_{\xi}, h_{\xi}(\alpha)=h_{\xi}(\beta)\}$.
For every $\alpha<\lambda^+$ let $S_{\alpha}:=\underset{\beta<\alpha}{\cup}u_{\beta}$,
$E_{\alpha}^*:= \cap \{E_{\xi} : \xi \in S_{\alpha}\}$ and let $E_*:=\underset{\alpha<\lambda^+}{\Delta}E_{\alpha}^*$,
so $E_{\alpha}^*$ $(\alpha<\lambda^+)$ and $E_* \subseteq \lambda^+$
are clubs. For every $\delta \in E_* \cap S_{\lambda}^{\lambda^+}$
define:

1. $u_{\delta}^*:=u_{\delta} \cap S_{\delta}$.

2. $h_{\delta}^*: u_{\delta}^* \rightarrow \delta$ is defined by
$h_{\delta}^*(\xi):=h_{\xi}(\delta)$ (recaling that $h_{\xi}(\delta)$
is well-defined and is $<\delta$).

3. $y_{\delta}^*=\{(otp(u_{\delta} \cap \zeta),f_{\delta}(\zeta)) : \zeta \in u_{\delta}^*\}$.

4. $S_{\delta}^2:=\{(h_*,y_*) : h_*$ is a function with domain $\in [S_{\delta}]^{<\lambda}$
and range $\subseteq \delta$, $y_* \in [\lambda \times (\lambda+1)]^{<\lambda} \}$.

Note that $\alpha<\beta \rightarrow S_{\alpha}^2 \subseteq S_{\beta}^2$
and that $|S_{\alpha}^2|\leq \lambda$ for every $\alpha$. Note also
that $S_{\alpha}^2=\underset{\beta<\alpha}{\cup}S_{\beta}^2$ when
$cf(\alpha)=\lambda$.

Now define a regressive function $g_*$ on $S_{\lambda}^{\lambda^+} \cap E_*$
such that $g_*(\delta_1)=g_*(\delta_2)$ iff $h_{\delta_1}^*=h_{\delta_2}^*$
and $y_{\delta_1}^*=y_{\delta_2}^*$ (this can be done as in the proof
of the $\lambda$-completeness of $D_{\lambda}^0$, see below). Let
$X=\{(\delta_1,\delta_2) : \delta_1<\delta_2 \in S_{\lambda}^{\lambda^+} \cap E_* \wedge g_*(\delta_1)=g_*(\delta_2)\}$,
then $X\in D_{\lambda}^0$ as witnessed by $E_*$ and $g_*$. Therefore
it's enough to prove that every $(\delta_1,\delta_2)\in X$, $(f_{\delta_1},f_{\delta_2})$
is a $\Delta$-system pair and $\xi \in u_{\delta_1} \cap u_{\delta_2}$
implies $(\delta_1,\delta_2) \in g(\xi)$. Indeed, as $g_*(\delta_1)=g_*(\delta_2)$,
it follows that $h_{\delta_1}^*=h_{\delta_2}^*$ and $y_{\delta_1}^*=y_{\delta_2}^*$,
hence $u_{\delta_1}^*=Dom(h_{\delta_1}^*)=Dom(h_{\delta_2}^*)=u_{\delta_2}^*$.
Note also that if $\zeta \in Dom(f_{\delta_1}) \cap Dom(f_{\delta_2})=u_{\delta_1} \cap u_{\delta_2}$,
then as $\delta_1<\delta_2$, it follows that $\zeta \in u_{\delta_2}^*=Dom(h_{\delta_1}^*)$.
Therefore $Dom(f_{\delta_1}) \cap Dom(f_{\delta_2})=Dom(h_{\delta_1}^*)$,
and it follows that $(f_{\delta_1},f_{\delta_2})$ is a $\Delta$-system
pair. If $\xi \in u_{\delta_1}\cap u_{\delta_2}=Dom(f_{\delta_1}) \cap Dom(f_{\delta_2})=Dom(h_{\delta_1}^*)=Dom(h_{\delta_2}^*)$,
then as $h_{\delta_1}^*=h_{\delta_2}^*$, it follows that $h_{\xi}(\delta_1)=h_{\delta_1}^*(\xi)=h_{\delta_2}^*(\xi)=h_{\xi}(\delta_2)$.
Therefore, $(\delta_1,\delta_2) \in X_{\xi} \subseteq g(\xi)$ and
we're done.

In order to show that $D_{\lambda}^0$ is $\lambda$-complete, let
$\zeta<\lambda$ and let $\{X_{\xi} : \xi<\zeta\} \subseteq D_{\lambda}^0$,
we shall prove that $\underset{\xi<\zeta}{\cap}X_{\xi} \in D_{\lambda}^0$.
For each $\xi<\zeta$, there are $E_{\xi}$ and $g_{\xi}$ as in the
definition of $D_{\lambda}^0$ witnessing that $X_{\xi} \in D_{\lambda}^0$.
Fix a bijection $f: (\lambda^+)^{<\lambda} \rightarrow \lambda^+$
and let $E=\{ \delta<\lambda^+ : \delta$ is a limit ordinal, and
for every $\alpha<\delta$ and $\eta \in \alpha^{<\lambda}$, $f(\eta)<\delta \}$,
then $E\subseteq \lambda^+$ is a club. Let $\delta \in E\cap S_{\lambda}^{\lambda^+}$,
then $f(\eta)<\delta$ for every $\eta \in \delta^{<\lambda}$. Define
a function $g: S_{\lambda}^{\lambda^+} \rightarrow \lambda^+$ as
follows: if $\delta \in S_{\lambda}^{\lambda^+} \cap E$, we let $g(\delta)=f((g_{\xi}(\delta) : \xi<\zeta))$.
Otherwise, we let $g(\delta)=0$. $g$ is a well-defined regressive
function. Let $E'=E\cap (\underset{\xi<\zeta}{\cap}E_{\xi})$, then
$E'\subseteq \lambda^+$ is a club. Let $X=\{(\alpha, \beta) : \alpha<\beta<\lambda^+, \alpha,\beta \in E' \cap S_{\lambda}^{\lambda^+}, g(\alpha)=g(\beta)\}$,
then as $X\in D_{\lambda}^0$, it suffices to show that $X\subseteq X_{\xi}$
for every $\xi<\zeta$. As $E'\subseteq E_{\xi}$ for every $\xi<\zeta$,
if $\alpha,\beta \in E'\cap S_{\lambda}^{\lambda^+}$ and $g(\alpha)=g(\beta)$,
then $g_{\xi}(\alpha)=g_{\xi}(\beta)$. This implies that $X\subseteq X_{\xi}$,
as required. This completes the proof of the claim. $\square$

\textbf{Definition 1.4: }Given a $(\lambda,D)-$forcing template $\mathbf p$
and a funtion $h: I_{\mathbf p}^{0}\rightarrow H(\lambda)$, we say that
the pair $(\mathbf{p},h)$ is \textbf{active} if:

A) $(\mathbb{Q}_{\mathbf p,h},\leq_{\mathbf p,h})$ is a forcing notion
where $\mathbb{Q}_{\mathbf p,h}=\{a\in H_{\leq \lambda}(\mathbf{U} \cup \mathbf{I}) : \mathfrak{B}_{\mathbf p}\models \varphi_0(a,h)\}$,
$\leq_{\mathbb{Q}_{\mathbf p,h}}=\{(a,b) : \mathfrak{B}_{\mathbf p} \models \varphi_1(a,b,h)\}$.

B) For every $\gamma<\lambda$ and $p\in \mathbb{Q}_{\mathbf{p},h}$
the formula $\varphi_2(-,\gamma,p,h)$ defines a winning strategy
for the player COM in the game $G_{\gamma}(p,\mathbb{Q}_{\mathbf p,h})$
(see definition 1.14 below).

Remark: The strategy may not provide a unique move and we shall allow
the completeness player to extend the condition given by the strategy.

C) Each element of $\mathbb{Q}_{\mathbf p,h}$ is a function of size
$\lambda$ with domain $\subseteq I_{\mathbf p}^{1}$ and range $\subseteq H(\lambda)$ (so this includes conditions that are sequences, trees, etc).

D) $\varphi_4(-,-,h)$ defines a function $tr$ such that $Dom(tr)=\mathbb{Q}_{\mathbf p,h}$
and for every $p\in \mathbb{Q}_{\mathbf p,h}$, $tr(p) \in \bold{T_p}$ is a function
with domain $X$ for some $X\in[I_{\mathbf p}^{1}]^{<\lambda}$ and
range $\subseteq H(\lambda)$, such that the following conditions hold: 
\\
\\
1) $p\leq q \rightarrow tr(p)\subseteq tr(q)$.

2) The formula $\varphi_5(-,-,h)$ defines a binary compatibility
relation $com\subseteq \mathbb{Q}_{\mathbf p, h} \times \mathbf{T}_{\mathbf p}$
(note that, in contrast with (6) below, this is a relation between conditions and trunks).

3) If $com(p, \eta)$ then:

a. There is $q$ such that $p\leq q$ such that $tr(q)=\eta$.

b. If $q\leq p$ then $com(q,\eta)$.

4) $\leq_{\mathbf p, h}$ is a partial ordering of $\mathbf{T}_{\mathbf{p}}$ such that
$\eta_1 \leq \eta_2 \rightarrow \eta_1 \subseteq \eta_2$.

5) If $p_1,p_2 \in \mathbb{Q}_{\mathbf p}$ and $tr(p_1)\mathbf{R}_{\mathbf p}tr(p_2)$
then $p_1,p_2 \in \mathbb{Q}_{\mathbf p,h}$ have a common upper bound
$q$. This is defined by $\varphi_6(-,-,h)$.

6) If $\eta \in \bold{T}_{\bold p, h}$, $j<|Dom(\eta)|$, $\{ p_i : i<j\}$ are conditions and $\underset{i<j}{\wedge}tr(p_i)=\eta$ then:

a. There is $q$ such that $\underset{i<j}{\wedge}(p_i \leq q)$.

b. There is a $\lambda$-Borel function $\mathbf{C}_{\mathbf p,h,j}$
such that $q=\mathbf{C}_{\mathbf p,h,j}(...,p_i,...)_{i<j}$ (recalling Clause (C) above) and $q$ is a least upper bound for $\{  p_i : i<j\}$.
\\
$[$This could be simplified by replacing "$j<|Dom(\eta)|$" by "$j<\lambda$", but that would exclude, e.g., random real forcing and the forcing $\mathbb{Q}_{\bar \theta}$ from [Sh1126] $]$
\\
\\
c. $tr(q)=tr(p_i)$ for all $i<j$.

7) $[$Follows from Definition 1.1(H)$]$ $\mathbb{Q}_{\mathbf p,h}$ satisfies the $(\lambda,D)$-chain condition:
if $p_{\alpha} \in \mathbb{Q}_{\mathbf p,h}$ $(\alpha<\lambda^+)$
then $\{(\alpha,\beta) : tr(p_{\alpha})R_{\mathbf p}tr(p_{\beta})\} \in D$. In Requirement 1.18 below we shall actually strengthen this condition and require that it holds in an absolute way as described there.

8) (Relevant for $\lambda>\aleph_0$) For every $\delta<\lambda$
and a play $(p_i, q_i : i<\delta)$ of length $<\lambda$ chosen according
to the winning strategy for the game in clause (B), there is a bound
$p_{\delta}$ given by the strategy such that $tr(p_{\delta})=\underset{i<\delta}{\cup}tr(p_i)$. 
\\
\\
9) For every $a\in I^1_{\bold p}$ and $x\in H(\lambda)$, there is some $p_{a, x} \in \mathbb{Q}_{\bold p, h}$ such that $\Vdash_{\mathbb{Q}_{\bold p, h}} "p_{a,x} \in \underset{\sim}{G}$ iff $\underset{\sim}{\eta_{\bold p, h}}(a)=x$" (where $\underset{\sim}{\eta_{\bold p, h}}$ is defined in the next clause).

E) 1. $\Vdash_{\mathbb{Q}_{\mathbf p,h}} "Dom(\underset{\sim}{\eta_{\mathbf p}})=I_{\mathbf p}^{1}"$
where $\underset{\sim}{\eta_{\mathbf p}}=\underset{\sim}{\eta_{\mathbf p,h}}$
is the $\mathbb{Q}_{\mathbf p,h}$-name of $\cup \{tr(q) : q\in \underset{\sim}{G_{\mathbb{Q}_{\mathbf p}}}\}$.

2. For every $b\in I_{\mathbf p}^1$ and $p\in \mathbb{Q}_{\mathbf p,h}$
then there is $\eta \in \mathbf{T}_{\bold p}$ such that $b\in Dom(\eta) \wedge com(p,\eta)$. Moreover, in Clause (D)(6), if we are given in addition some $a\in I_{\bold p}^1 \setminus Dom(\eta)$, then there is $\nu$ extending $\eta$ such that $a\in Dom(\nu)$ and $com(p_i, \nu)$ for every $i<j$ (and so there exists $q_i$ above $p_i$ such that $tr(q_i)=\nu$ for every $i<j$).

F) $\underset{\sim}{\eta_{\mathbf p}}$ is generic for $\mathbb{Q}_{\mathbf p,h}$,
i.e. there is a $\lambda$-Borel function $\mathbf{B}$ defined in $V$
such that $\Vdash "p\in \underset{\sim}{G}$ iff $\mathbf{B}(p,\underset{\sim}{\eta_{\mathbf p}})=true"$
for every $p\in \mathbb{Q}_{\mathbf p,h}$.

G) If $p$ and $q$ are incompatible and $tr(p)\subseteq tr(q)$,
then $p\Vdash_{\mathbb{Q}_{\mathbf p,h}} "tr(q) \nsubseteq \underset{\sim}{\eta_{\mathbf p}}"$.
In this case we shall say that $p$ and $tr(q)$ are incompatible. 

H) If $j<\lambda$, $p_i \in \mathbb{Q}_{\mathbf p,h}$ $(i<j)$
and $q$ are as in 1.4(D)(6) and $p$ is a condition such that $tr(q)\subseteq tr(p)$
and such that $q$ and $tr(p)$ are incompatible, then there is $i<j$
such that $\{p_i, tr(p)\}$ are incompatible.

Remarks: 1. If $(\bold p, h)$ is not active, then we let $\mathbb{Q}_{\bold p, h}$ be trivial.
\\
2. The reader may wonder where the properties of forcing templates, their trunks, etc, are used in the construction of the iterations that will follow. This will play a major role in the proof of Claim 2.10. 
\\
3. Clauses (G)+(H) will be used later, for example, in Claim
4.1.
\\
\\
Below we shall give several examples of concrete forcing notions as the realizations of forcing templates.
\\
\\
\textbf{Example 1.4(A)}: Let $\lambda$ be either an inaccessible cardinal or $\aleph_0$ and assume that $P, g$ and $h$ are functions with domain $\lambda$ such that:
\\
\\
a. For every $\alpha<\lambda$, $P(\alpha)$ is a partial order of cardinality $<\lambda$.
\\
\\
b. For every $\alpha<\lambda$, $g(\alpha)$ is a regular cardinal from $(\alpha, \lambda)$ (relevant in the inaccessible case).
\\
\\
c. For every $\alpha<\lambda$, $h(\alpha): P(\alpha) \rightarrow g(\alpha)$ is a function such that $P(\alpha) \models a\leq b \rightarrow h(\alpha)(a) \leq h(\alpha)(b)$. 
\\
\\
d. If $\lambda>\aleph_0$ then for every $\alpha<\lambda$, $g(\alpha)=cf(g(\alpha))>\alpha$ and $P(\alpha)$ is $(<g(\alpha))$-directed. If $\lambda=\aleph_0$ then $P(\alpha)$ has a maximal element.
\\
\\
Let $\mathbb Q=\mathbb{Q}_{P,g,h}$ be the following forcing notion:
\\
\\
1. $p\in \mathbb Q$ iff: 
\\
\\
a. $p=(\eta, \rho, \nu)=(\eta_p, \rho_p, \nu_p)$.
\\
\\
b. $\rho \in \underset{\alpha \in [lg(\eta), \lambda)}{\prod}g(\alpha)$.
\\
\\
c. $\nu \in \underset{\alpha \in [lg(\eta), \lambda)}{\prod}P(\alpha)$.
\\
\\
d. If $\alpha \in [lg(\eta), \lambda)$ then $h(\alpha)(\nu(\alpha)) \leq \rho(\alpha)$.
\\
\\
e. $\eta \in \underset{\alpha <lg(\eta)}{\prod}P(\alpha)$
\\
\\
f. If $\lambda=\aleph_0$, then $\underset{i}{lim}(g(i)-\rho(i))=\infty$.
\\
\\
g. We let $tr(p):=\eta$.
\\
\\
2. Given $p,q \in \mathbb Q$, $p\leq q$ iff $\eta_p \subseteq \eta_q$, $\rho_p(\alpha) \leq \rho_q(\alpha)$ for every $\alpha \in [lg(\eta_q), \lambda)$, $P(\alpha) \models \nu_p(\alpha) \leq \nu_q(\alpha)$ for every $\alpha \in [lg(\eta_q),\lambda)$ and $P(\alpha) \models \nu_p(\alpha) \leq \eta_q(\alpha)$ for every $\alpha \in [lg(\eta_p),lg(\eta_q))$.
\\
\\
We shall now define a forcing template $\bold p$ that gives rise to the above forcing:
\\
\\
a. $\lambda_{\bold p}:=\lambda$, $\kappa_{\bold p}=\lambda^+$.
\\
\\
b. $I^0_{\bold p}=I^1_{\bold p}=\lambda$.
\\
\\
c. $\mathfrak{B}_{\bold p}$ and $\mathfrak{B}^0_{\bold p}$ will be trivial, i.e. $(H(\lambda^+), \in)$.
\\
\\
d. Denote by $h^*$ the function $h: I^0_{\bold p} \rightarrow H(\lambda)$ in the definition of active forcing templates. $h^*$ here will be given here by $h^*(\alpha)=(P(\alpha), g(\alpha), h(\alpha))$.
\\
\\
e. The formulas $\varphi_k$ will then define $\mathbb{Q}_{P,g,h}$ as described above using the parameter $h^*$. Denote the trunks in this case by $tr_{\bold p, h^*}(p)$.
\\
\\
f. $\bold{T}_{\bold p}=\{ tr_{\bold p, h^*}(p) : h^*, p$ as above$\}$.
\\
\\
g. $R_{\bold p}= \{ (\eta_1, \eta_2) \in \bold{T}_{\bold p} \times \bold{T}_{\bold p} : \eta_1 =\eta_2 \}$.
\\
$[$Note that while we allow the parameter $h^*$ to be a name, $\bold{T_p}$ and $R_{\bold p}$ are objects.$]$
\\
\\
h. $D_{\bold p}$ will be the filter $D^0_{\lambda}$ from Claim 1.3.
\\
\\
For a typical example of a triple $(P, g, h)$, consider a sequence $(\theta_{\alpha}, \sigma_{\alpha} : \alpha<\lambda)$ where $\alpha<\sigma_{\alpha}<\theta_{\alpha} <\lambda$. For each $\alpha$ let $P(\alpha)=([\theta_{\alpha}]^{<\sigma_{\alpha}}, \subseteq )$. For every $\alpha< \lambda$ let $g(\alpha)=\sigma_{\alpha}$ and for every $u\in P(\alpha)$ let $h(\alpha)(u)=otp(u)$.
\\
\\
$\textbf{Remark 1.4(B)}$: 1. On such forcing notions see [Sh628], [Sh872], [HwSh1067] for $\lambda=\aleph_0$ and [Sh1126] for inaccessible $\lambda$. In [Sh1126] we have $P(\alpha)=\{[\epsilon, \theta_{\alpha}] : \epsilon <\theta_{\alpha} \}$ with the reverse ordering, $g(\alpha)=\theta_{\alpha}$ which is regular $>|\alpha|$ and $h(\alpha)([\epsilon, \theta_{\alpha}])=\epsilon$.
\\
\\
2. The above example gives a justification for the (somewhat arbitrary) use of the assumption "$j<|Dom(\eta)|$" (rather than "$j<\lambda$") in Definition 1.4(D)(6).
\\
\\
Below is an additional example where $R$ is nontrivial:
\\
\\
$\textbf{Example 1.4(C)}$: Our next example is random real forcing with a modification needed to satisfy the requirement in Definition 1.4(D)(6). Let $(\eta_n : n<\omega)$ enumerate $2^{<\omega}$ without repetition and let $D_{\bold p}=D_{\aleph_0}^0$.
\\
\\
A. $p\in \mathbb Q$ iff $p=(tr(p), B_p)$ where:
\\
\\
a. $B_p\subseteq 2^{\omega}$ is Borel.
\\
\\
b. $\mu(B_p) >0$.
\\
\\
c. $\eta_p$ is the maximal element of $2^{<\omega}$ that is an initial segment of all members of $B_p$.
\\
\\
d. There is a natural $n(p)>1$ such that $2^{lg(\eta_p)} \mu(B_p) \in [1-\frac{1}{n(p)+1}, 1-\frac{1}{n(p)+2}]$ and $n(p) \leq lg(\eta_p)$.
\\
\\
e. $tr(p)$ is a constantly $1$ function with domain $\{ \eta_p\} \cup \{ \eta_p \restriction l : l<n(p)\}$.
\\
\\
B. For $p, q\in \mathbb Q$, $\mathbb Q \models p\leq q$ iff:
\\
\\
a. $B_q \subseteq B_p$.
\\
\\
b. $tr(p) \subseteq tr(q)$.
\\
\\
C. The generic will be the union of $\eta_p$ for every $p\in G$.
\\
\\
D. $\bold{T_p}=\{tr(p) : p \in \mathbb Q\}$.
\\
\\
E. $R_{\bold p}=\{ (\eta, \eta) : \eta \in \bold{T_p}\}$.
\\
$[$This gives an example where $R_{\bold p}$ is not the usual function compatibility. Note that as random real forcing is not $\sigma$-centered, we can't strengthen Definition 1.4(D)(6) to "$j<\lambda$".$]$
\\
\\
\textbf{Remark 1.4(D)}: The trunks will play a role in the definition of our iterations, where given a condition $p$ and $s\in Dom(p)$, $p(s)$ will be a name of a condition consisting of a trunk $tr(p(s))$ and a condition computed from names of other conditions of the form $p_{\iota}=\bold{B}_{p(s), \iota} (...,\underset{\sim}{{\eta}_{t_{\zeta}}}(a_{\zeta}),...)_{\zeta \in W_{p(s),\iota}}$ (this notation will be explained in due course) whose union of trunks is $tr(p(s))$. All of this will eventually play a role in the analysis of projections in Section 4.

\textbf{\large{}Iteration templates}{\large\par}

Similarly to forcing templates, iteration templates will contain the
information from which we shall construct our iterations. This information
will include a well-founded partial order along which we shall define
the iteration. For every element in the partial order, we shall assign
a forcing template and two types of memory: a strong memory which
will be used for the construction of the forcing conditions, and a
weak memory which will be used to define the necessary parameter for defining the forcing at the current stage. The parameters will then be computed
in a $\lambda$-Borel way from the previous generics. An additional complication in our memory apparatus (i.e. the $v_t \subseteq [u_t^0]^{\leq \lambda}$ in Definition 2.2.A) will then require a corresponding modification of our notion of strategic completeness in Definition 1.14.

\textbf{Definition 1.5: }A $(\lambda,D)$-iteration template $\mathbf q$
consists of the objects $\{L_{\mathbf q}, (\mathbf{p}_t : t\in L_{\mathbf q}), ((u_{t}^{0},\bar{u}_{t}^{1}) : t\in L_{\mathbf q}),((w_{t}^{0},\bar{w}_{t}^{1}) : t\in L_{\mathbf q}), D_{\mathbf q}, ((\mathbf{B}_{t,b}, (s_t(b,\zeta), a_{t,b,\zeta}) : \zeta<\xi(t,b)) : b\in I_{\mathbf{p}_t}^{0}) : t\in L_{\mathbf q}))\}$
such that:

A) $D_{\mathbf q}=D$, $L_{\mathbf q}$ is a well-founded partial order
with elements from $\mathbf{U}$.

B) For every $t\in L_{\mathbf q}$, $\mathbf{p}_t=\mathbf{p}_{\mathbf{q},t}$
is a $(\lambda,D)$ forcing template. Note that $D$ is fixed filter that doesn't depend on $t$. 

C) For every $t\in L_{\mathbf q}$, $u_{\mathbf q,t}^{0}=u_{t}^{0}\subseteq L_{<t}=\{s\in L_{\bold q} : s<_{L_{\bold q}}t\}$
and $\bar{u}_{\mathbf q,t}^{1}=\bar{u}_{t}^{1}=(u_{t,s}^{1} : s\in u_{t}^{0})$
where $u_{t,s}^{1}\subseteq I_{s}^{1}=I_{\mathbf{p}_s}^{1}$. We shall
refer to $u_{\mathbf q,t}^{0}$ as strong memory.

D) For every $t\in L_{\mathbf q}$, $w_{t}^{0}\subseteq u_{t}^{0}$
and $\bar{w}_{t}^{1}=(w_{t,s}^{1} : s\in w_{t}^{0})$ where $w_{t,s}^{1}\subseteq u_{t,s}^{1} \subseteq I_{s}^{1}$.
We shall refer to $w_{t}^{0}$ as weak memory.

Remark: In many interesting cases, $w_t^0=\emptyset$ for all $t$
(this will correspond to an iteration where the definitions of the
forcing notions are without parameters).

E) For every $t\in L_{\mathbf q}$ and $b\in I_{\mathbf{p}_t}^{0}$, $\mathbf{B}_{t,b}$
is a $\lambda-$Borel $\xi(t,b)-$place function $(\xi(t,b)<\lambda^+)$
from $\lambda^{\xi(t,b)}$ to $\lambda$. For every $\zeta<\xi(t,b)$
we have $s_t(b,\zeta) \in w_{t}^{0}$ and $a_{t,b,\zeta} \in w_{t,s_t(b,\eta)}^{1}$ (if $w_t^0=\emptyset$ then $\xi(t,b)=0$).
\\
$[$This will be used to compute $h$ when applying Definition 1.4.$]$

F) $D_{\mathbf q}$ is a $\lambda$-complete filter as in Hypothesis
0 such that $D_{\mathbf{p}_t}=D_{\mathbf q}$ for every $t\in L_{\mathbf q}$.

\textbf{Definition 1.6(A): }Given an iteration template $\mathbf q$
and $L\subseteq L_{\mathbf q}$, let $cl(L)=cl_{\mathbf q}(L)$ be the
minimal $L'$ such that $L\subseteq L' \subseteq L_{\mathbf q}$ and
$t\in L' \rightarrow w_{\mathbf q,t}^0 \subseteq L'$.
\\
\\
\textbf{Example 1.6(B)}: We shall briefly illustrate how to construct a concrete iteration within our general framework continued below. Let $\lambda$ be either $\aleph_0$ or inaccessible with $\bar \theta=(\theta_i : i<\lambda)$ a sufficiently fast increasing sequence such that $\theta_i=cf(\theta_i)>i$. Fix an ordinal $\alpha_*$ and let $(U_1, U_2, U_3)$ be a partition of $\alpha_*$. For $\alpha<\alpha_*$, let $\bar{\varphi_{\alpha}}$ define:
\\
a. Random real forcing (as in Example 1.4(C)) if $\alpha \in U_0$ and $\lambda=\aleph_0$.
\\
b. Random real forcing for inaccessible $\lambda$ (see [Sh:1004]) if $\alpha \in U_0$ and $\lambda$ is inaccessible.
\\
c. The forcing from Example 1.4(A) if $\alpha \in U_2$ and $\lambda=\aleph_0$.
\\
d. The forcing $\mathbb{Q}_{\bar \theta}$ from [Sh:945] if $\alpha \in U_2$ and $\lambda$ is inaccessible.
\\
e. Hechler forcing ($\lambda$-Hechler forcing) if $\alpha \in U_3$ and $\lambda=\aleph_0$ ($\lambda$ is inaccessible).
\\
\\
The filter $D$ will be $D^0_{\lambda}$ from Claim 1.3. If, for example, $\underset{\sim}{\mathbb{Q}_{t}}$ is $\mathbb{Q}_{\bar \theta}$ from [Sh945], then we might use a parameter $\bar \theta \in V$, but we might also want to use a parameter of the form $\bar \theta=\bold{B}(...,\underset{\sim}{\eta_{\zeta}}(a),...)$ where each $\zeta$ belongs to the weak memory $w^0_t$.
\\
\\
For every $\alpha < \alpha_*$, $u_{\alpha}^0$ will be a subset of $\alpha$. Note that if $\alpha_l \in U_2$ $(l=1,2,3)$, $\alpha_1 \in u_{\alpha_2}^0, \alpha_2 \in u_{\alpha_3}^0$ and $\alpha_1 \notin u_{\alpha_3}^0$, then it will still be forced that "$\underset{\sim}{\eta_{\alpha_1}} <_{bd} \underset{\sim}{\eta_{\alpha_3}}$". In [Sh:945] and [Sh:1126] the case $\alpha_*=U_2$ was used.

\textbf{Definition 1.7:} 1. Let $\mathbf{P}$ be a set of forcing templates,
we shall denote by $\mathbf{K}_{\mathbf P}$ the collection of iteration
templates $\mathbf q$ with forcing templates from $\mathbf{P}$ (i.e.
$\mathbf{p}_{\mathbf q,t} \in \mathbf{P}$ for every $t\in L_{\mathbf q}$).

2. For $\mathbf{q}_1,\mathbf{q}_2 \in \mathbf{K_P}$ we write $\mathbf{q}_1 \leq_{\mathbf{K_P}} \mathbf{q}_2$
if the following conditions hold:

a. $L_{\mathbf{q}_1} \subseteq L_{\mathbf{q}_2}$.

b. For every $t\in L_{\mathbf{q}_1}$, $\mathbf{p}_{\mathbf{q}_1,t}=\mathbf{p}_{\mathbf{q}_2,t}$
and $u_{\mathbf{q}_1}^0=u_{\mathbf{q}_2}^0 \cap L_{\mathbf{q}_1}$.

c. $(w_{\mathbf{q}_1,t}^0, \bar{w}_{\mathbf{q}_1,t}^1 : t\in L_{\mathbf{q}_1})=(w_{\mathbf{q}_2,t}^0, \bar{w}_{\mathbf{q}_2,t}^1 : t\in L_{\mathbf{q}_2}) \restriction L_{\mathbf{q}_1}$
and similarly for the other sequences appearing in definition 1.4.

\textbf{Definition 1.8: }Let $\mathbf q$ be an iteration template and
let $L\subseteq L_{\mathbf q}$, we shall say that $L$ is a closed
sub-partial order (or ``$L$ is closed with respect to weak memory'')
if $w_{t}^{0}\subseteq L$ for every $t\in L$.

\textbf{Definition 1.9: }1. Given $L\subseteq L_{\mathbf q}$, let $cl(L)=cl_{\mathbf q}(L)$
be the minimal set $L\subseteq L' \subseteq L_{\mathbf q}$ such that
$w_t^0 \subseteq L'$ for every $t\in L'$. 

\textbf{Convention 1.9(A): }Throughout this paper, whenever $\mathbf q$
is an iteration template, $L\subseteq L_{\mathbf q}$ and $\mathbf q \restriction L$
is defined or used (see definition 1.11), we shall assume that $L$ is closed
w.r.t. weak memory.

\textbf{Definition 1.10: }Let $\mathbf q$ be an iteration template,
we shall define for every $t\in L_{\mathbf q} \cup \{\infty\}$ a forcing
notion $\mathbb{P}_t=\mathbb{P}_{\mathbf q,t}$, a forcing notion $\mathbb{P}_{L}=\mathbb{P}_{\mathbf q,L}$
for any initial segment $L\subseteq L_{\mathbf q}$ and names $\underset{\sim}{\mathbb{Q}_t}=\underset{\sim}{\mathbb{Q}_{\mathbf q,t}}$,
$\underset{\sim}{\eta_t}$ (by the remark after Definition 1.4, this is always well-defined) by induction on $dp(t)$ (see definition
2.3):

A) $p\in \mathbb{P}_t$ $(\mathbb{P}_L)$ iff

1) $p$ is a function with domain $\subseteq$ $L_{<t}$ (or $\subseteq L$
in the case of $\mathbb{P}_L$) of cardinality $<\lambda$.

2) For every $s\in Dom(p)$, $p(s)=\mathbf{B}_{p(s)}(...,\underset{\sim}{\eta_{t_{\zeta}}}(a_{\zeta}),...)_{\zeta<\xi}$
(we may write $p(s)=(tr(p(s)), \mathbf{B}_{p(s)}(...,\underset{\sim}{\eta_{t_{\zeta}}}(a_{\zeta}),...))$, so it will be interpreted as a condition in $\underset{\sim}{\mathbb{Q}_s}$ that resulted from the respective computation by the $\lambda$-Borel function $\mathbf{B}_{p(s)}$)
for a $\lambda$-Borel function $\mathbf{B}_{p(s)}$ into $H_{\leq \lambda}(\mathbf{U} \cup \mathbf{I})$
and an object $tr(p(s))$ such that $tr(p(s))$ is computable from
$\mathbf{B}_{p(s)}$ (i.e. the range of $\mathbf{B}_{p(s)}$ consists
of conditions with trunk $tr(p(s))$), $\xi=\xi_{p(s)}\leq \lambda$,
$\{t_{\zeta} : \zeta<\xi\}\subseteq u_{s}^{0}$ and for every $\zeta$,
$a_{\zeta}\in u_{t_{\zeta}}^{1}$. Note that $\bold{B}_{p(s)}$ here is not the same function as $\bold{B}_{t, b}$ in Definition 1.5.
\\
$[$Remarks: a. The reader might wonder why not drop the $a_{\zeta}$ and use $\mathbf{B}_{p(s)}(...,\underset{\sim}{\eta_{t_{\zeta}}},...)_{\zeta<\xi}$ instead. The reason is that $Dom(\underset{\sim}{\eta_{t_{\zeta}}})=I_{t_{\zeta}}^1$ might be of cardinality $>\lambda$. Our choise allows $\bold{B}_{p(s)}$ to be a function with domain $H(\lambda)^{\xi}$.
\\
b. Note that if $p\leq q$ and $s\in Dom(p)$, then the corresponding set of $\{ t_{\zeta} : \zeta< \xi\}$ might increase. As a consequence, the number of input coordinates might increase between $\bold{B}_{p(s)}$ and $\bold{B}_{q(s)}$. $]$

3) For every $s\in Dom(p)$, $ \Vdash_{\mathbb{P}_s} "p(s)\in \underset{\sim}{\mathbb{Q}_s}"$.

B) $\mathbb{P}_t \models p\leq q$ iff $Dom(p) \subseteq Dom(q)$
and for every $s\in Dom(p)$, $q\restriction L_{<s} \Vdash_{\mathbb{P}_s} p(s) \leq_{\underset{\sim}{\mathbb{Q}_s}} q(s)$.

C) 1. Let $h_t : I_{\mathbf{p}_t}^{0}\rightarrow \lambda$ be the name
of a function defined by $h_t(b)=\mathbf{B}_{t,b}(...,\underset{\sim}{\eta_{s_t(b,\zeta)}}(a_{t,b,\zeta}),...)_{\zeta<\xi(t,b)}$. 

2. a. If $(\mathbf{p}_t,h_t)$ is active in $V^{\mathbb{P}_t}$ (see
Definition 1.4), we shall define $\underset{\sim}{\mathbb{Q}_t}$
as the $\mathbb{P}_t$-name of $\mathbb{Q}_{\mathbf{p}_t,h_t}^{V[\underset{\sim}{\eta_s} : s\in u_{t}^{0}]}$.

b. If $(\mathbf{p}_t,h_t)$ is not active in $V^{\mathbb{P}_t}$, we
shall define $\underset{\sim}{\mathbb{Q}_t}$ as the trivial forcing.

D) $\underset{\sim}{\eta_t}$ will be defined as the $\mathbb{P}_t \ast \underset{\sim}{\mathbb{Q}_t}$
name $\underset{\sim}{\eta_{\mathbf{p}_t,h_t}}$.

\textbf{Definition 1.11: }Given an iteration template $\mathbf q$ and
a sub partial order $L\subseteq L_{\mathbf q}$ we shall define the
iteration template $\mathbf q \restriction L$ as follows (recall that we assume that $L$ is closed under weak memory):

A) $L_{\mathbf q \restriction L}=L$.

B) For every $t\in L$, $\mathbf{p}_{\mathbf q\restriction L,t}=\mathbf{p}_{\mathbf q,t}$.

C) For every $t\in L$, $u_{\mathbf q \restriction L,t}^{0}=u_{\mathbf q,t}^{0}\cap L$
and $\bar{u}_{\mathbf{q}\restriction L,t}^{1}=\bar{u}_{\mathbf q,t}^{1}\restriction u_{\mathbf q\restriction L}^{0}$.

D) For every $t\in L$, $w_{\mathbf q\restriction L,t}^{0}=w_{\mathbf q,t}^{0}$
and $\bar{w}_{\mathbf q\restriction L,t}^{1}=\bar{w}_{\mathbf q,t}^{1}$.

E) For every $t\in L$ the other objects in the definition of $\mathbf q$
are not changed.

\textbf{Observation 1.12: }$\mathbf q\restriction L$ is an iteration
template (recall that $L$ is assumed to be closed under weak memory).

\textbf{Definition 1.13: }Let $\lambda$ be a regular cardinal, $\mathbb{P}$
a forcing notion and $Y\subseteq \mathbb{P}$. 

A) $\mathbb{L}_{\lambda^{+}}(Y)$ will be defined as the closure of
$Y$ under the operations $\neg$, $\underset{i<\alpha}{\wedge}$
for $\alpha < \lambda^+$.

B) For a generic set $G\subseteq \mathbb{P}$ and $\psi \in \mathbb{L}_{\lambda^+}(Y)$
the truth value of $\psi[G]$ will be defined naturally by induction
on the depth of $\psi$ (for example, for $p\in \mathbb{P}$, $p[G]=true$
iff $p\in G$).

C) The forcing $\mathbb{L}_{\lambda^+}(Y,\mathbb{P})$ will be defined
as follows: 

1) $\psi \in \mathbb{L}_{\lambda^+}(Y,\mathbb{P})$ iff $\psi \in \mathbb{L}_{\lambda^{+}}(Y)$
and $\nVdash_{\mathbb{P}} "\psi[\underset{\sim}{G}]=false"$.

2) $\psi_1 \leq \psi_2$ iff $\Vdash_{\mathbb{P}} "\psi_2[\underset{\sim}{G}]=true \rightarrow \psi_1[\underset{\sim}{G}]=true"$.

\textbf{\large{}More definitions and assumptions}{\large\par}

\textbf{Strategic completeness}

\textbf{Definition 1.14: } Let $\mathbb{P}$ be a forcing notion,
$\alpha \in Ord$ and $p\in \mathbb{P}$. 
\\
\\
1. The two player game $G_{\alpha}^0(p,\mathbb{P})$
will be defined as follows: 

A play in the game consists of $\alpha$ moves. In the $\beta$th move player
\textbf{I} chooses $p_{\beta}\in \mathbb{P}$ such that $p\leq p_{\beta} \wedge (\underset{\gamma<\beta}{\wedge}q_{\gamma}\leq p_{\beta})$,
player \textbf{II }responds with a condition $q_{\beta}$ such that
$p_{\beta}\leq q_{\beta}$.

Winning condition: Player \textbf{I }wins the play iff for each $\beta<\alpha$
there is a legal move for him.
\\
\\
2. Let $\mathbb P$ be a forcing notion, $tr=tr_{\mathbb P}$ a function from $\mathbb P$ into $\{ \eta : \eta$ is a function from a set of cardinality $<\lambda$ into $H(\lambda) \}$, $\alpha \in Ord$ and $p\in \mathbb P$. The game $G_{\alpha}^1(p, \mathbb P)$ will be defined as follows: The games consists of $\alpha$ moves. In the $\epsilon$th move the objects $j_{\epsilon}, \bar{q}_{\epsilon}, \eta_{\epsilon}, \nu_{\epsilon}$ are chosen such that:
\\
\\
a. $j_{\epsilon} < \lambda$ and $\xi \leq \epsilon \rightarrow j_{\xi} \leq j_{\epsilon}$.
\\
\\
b. $\bar{q}_{\epsilon}=(q^{\epsilon}_i : i<j_{\epsilon})$ is a sequence of members of $\mathbb P$ above $p$.
\\
\\
c. If $\xi < \epsilon$ and $i<j_{\xi}$ then:
\\
\\
c(1). $(q^{\zeta}_i : \zeta \in [\xi, \epsilon])$ is increasing.
\\
\\
c(2). $tr(q_i^{\epsilon})=\eta_{\epsilon}$.
\\
\\
c(3). $j_{\epsilon} \leq |Dom(\eta_{\epsilon})|$.
\\
\\
c(4). $com(q_i^{\epsilon}, \nu_{\epsilon})$ for every $i<j_{\epsilon}$.
\\
\\
c(5). $\eta_{\epsilon} \subseteq \nu_{\epsilon}$.
\\
\\
c(6). $\zeta < \epsilon \rightarrow \nu_{\zeta} \subseteq \eta_{\epsilon}$.
\\
\\
In the $\epsilon$th move, first INC chooses $j_{\epsilon}$, $\bar{q}_{\epsilon}$ and $\eta_{\epsilon}$, then COM chooses $\nu_{\epsilon}$. COM wins if he has a legal move at every stage during the play.
\\
\\
3. Let $\mathbb P$ be a forcing notion expanded by a function $tr=tr_{\mathbb P}$ as in (2). Let $\alpha \in Ord$ and let $\bar{F}_{\alpha}=(F_{\alpha, \epsilon} : \epsilon  < \alpha)$ be a winning strategy for $\textbf{I}$ in the game $G^0_{\alpha}(-, \mathbb P)$ that will naturally arise from the rest of the definition below. The game $G^2_{\alpha}(\mathbb P)$ will be defined as follows: In the $\epsilon$th move, the objects $j_{\epsilon}$, $\bar{p}_{\epsilon}$, $\bar{q}_{\epsilon}$, $\eta_{\epsilon}$ and $\nu_{\epsilon}$ such that:
\\
\\
a. In a preliminary move, $\textbf{II}$ chooses $\xi \in (0, \alpha)$ and $q_*$.
\\
\\
b. For $\zeta < \xi$, we let $j_{\zeta}=1$, $p_0^{\zeta}=q_0^{\zeta}=q_*$ (so $\bar{p}_{\zeta}=\bar{q}_{\zeta}=(q_*))$ and $\eta_{\zeta}=\nu_{\zeta}=tr(q_*)$.
\\
\\
c. $j_{\epsilon} < \lambda$ and $\xi < \epsilon \rightarrow j_{\xi} \leq j_{\epsilon}$.
\\
\\
d. Given $\epsilon$ and $i<j_{\epsilon}$:
\\
\\
d(1). $\bar{q}_{\epsilon}=(q^{\epsilon}_i : i<j_{\epsilon})$ and $\bar{p}_{\epsilon}=(p^{\epsilon}_i : i<j_{\epsilon})$ are sequences of members of $\mathbb P$.
\\
\\
d(2). $tr(q_i^{\epsilon})=\eta_{\epsilon}$.
\\
\\
d(3). $j_{\epsilon} \leq |Dom(\eta_{\epsilon})|$.
\\
\\
d(4). $com(q_i^{\epsilon}, \nu_{\epsilon})$ for every $i<j_{\epsilon}$.
\\
\\
d(5). $\eta_{\epsilon} \subseteq \nu_{\epsilon}$ and $\zeta < \epsilon \rightarrow \nu_{\zeta} \subseteq \eta_{\epsilon}$.
\\
\\
In the $\epsilon$th move for $\epsilon \geq \xi_*$, first COM chooses $(p_i^{\epsilon} : i< \cup \{j_{\zeta} : \zeta< \epsilon \})$ such that $p_i^{\epsilon}=F_{\alpha, \epsilon}((q_i^{\zeta} : \zeta < \epsilon))$. Next INC chooses $j_{\epsilon} = \cup \{ j_{\zeta} : \zeta < \epsilon\}$, $\bar{q}_{\epsilon}$ such that $p_i^{\epsilon} \leq q_i^{\epsilon}$ for all $i<j_{\epsilon}$ and $\eta_{\epsilon}$ as above. Finally, COM chooses $\nu_{\epsilon}$ as above. COM wins if at each stage there is a legal move for him.

4, Let $\mathbb{P}$ be a forcing notion and $\alpha \in Ord$, $\mathbb{P}$
is called $\alpha$-strategically $i$-complete ($i=0,1,2$) if for each $p\in \mathbb{P}$
player \textbf{I }has a winning strategy for $G_{\alpha}^i(p,\mathbb{P})$.

5. For a regular $\lambda$, we say that $\mathbb{P}$ is $(<\lambda)$-strategically $i$-complete ($i=0,1,2$)
if it's $\alpha$-strategically $i$-complete for every $\alpha<\lambda$.
\\
\\
6. Convention: We may omit the $i$ in $i$-completeness if $i=2$.
\\
\\
For discussion of various strategic completeness properties see [Sh:587].
\\
\\
We shall freely use the following fact:

\textbf{Fact 1.15: }$(<\lambda)$-strategic completeness is preserved
under $(<\lambda)$-support iterations.

\textbf{Absoluteness}

The following requirements will be assumed throughout the paper for
all $(\lambda,D)$-forcing templates $\mathbf p$:

\textbf{Requirement 1.16: }A $(\lambda,D)$ forcing template $\mathbf p$
is called $(\lambda,D)$-absolute when: If $\mathbb{P}_1$ and $\mathbb{P}_2$
are $(<\lambda)$-strategically complete forcing notions satisfying
$(\lambda,D)-cc$ (that is, $\{p_{\alpha} : \alpha<\lambda^+\} \subseteq \mathbb{P}_l \rightarrow \{(\alpha, \beta) : p_{\alpha}$ and $p_{\beta}$ are compatible$\} \in D$) such that $\mathbb{P}_1 \lessdot \mathbb{P}_2$,
$V_l=V^{\mathbb{P}_l}$ $(l=1,2)$ and $\mathbf p, h \in V_1$, then we
shall require that:

A) "$(\bold p, h)$ is active" and $"p\leq_{\mathbb{Q}_{\mathbf p,h}}q"$ is absolute between $V_1$
and $V_2$.

B) $"p\in \mathbb{Q}_{\mathbf p,h}"$ is absolute between $V_1$ and
$V_2$.

C) $"p$ and $q$ are incompatible in $\mathbb{Q}_{\mathbf p,h}"$ is
absolute between $V_1$ and $V_2$.

D) Similarly for the other formulas involved in the definition of
$\mathbf p$ (see definition 1.1).

\textbf{Definition 1.17: }Let $\mathbf p \in V_1$ be a forcing template
and let $\mathbf B$ be a $\lambda$-Borel function. We say that $\mathbf B$
is a $\lambda$-Borel function into $\mathbf p$ if for every $V_1 \subseteq V_2$
as above, the range of $\mathbf B$ is in $\mathbb{Q}_{\mathbf p,h}^{V_2}$
and the trunk of the members in the range is fixed. 
\\
\\
Remarks: The above definition is relevant in the context, e.g., of Definition 1.10(A)(2), where $(V_1, V_2)$ here stands for $(V, V^{\mathbb{P}_s})$ there.
\\
\\
\textbf{Requirement 1.18: }A) All $\lambda$-Borel functions will
be assumed to be into a relevant forcing template $\mathbf p$. That is, whenever a $\lambda$-Borel function $\bold B$ will be used, there will be an associated forcing template $\bold p$ such that $(\bold B, \bold p)$ are as in Definition 1.18, and $\bold p$ will be clear from the context.

B) $D_{\mathbf p}$ is fixed and is in $V$.

$\\$

\textbf{\large{}2. Iteration parameters and the corrected iteration}{\large\par}

\textbf{\large{}Iteration parameters}{\large\par}

We will be interested in iterations along a prescribed partial order
$M$. However, we will also have to consider iterations along a larger
partial order that $L$ that contains $M$. Therefore, we shall define
a binary relation $E'$ on $L$ such that $L \setminus M$ will consist
of equivalence classes that are only related via $M$. We shall require
that those equivalence classes will be preserved when we extend the
iteration, so extensions will be obtained by adding new equivalence
classes.

\textbf{Hypothesis 2.1: }We shall assume in this section that:

A) $\lambda=\lambda^{<\lambda}$ is a cardinal and $D$ is a filter
as in Hypothesis 0.

B) $\lambda \leq \lambda_1 \leq \lambda_2$ are cardinals such that
$\beth_3(\lambda_1)\leq \lambda_2$.

C) $\mathbf P$ is a set of $(\lambda, D)$-forcing templates that are
$(\lambda,D)$-absolute such that if $\mathbf p \in \mathbf P$ and $\mathbb{P}_1 \lessdot \mathbb{P}_2$ are $(<\lambda)$-strategically complete $(\lambda,D)$-cc
forcing notions, then $V^{\mathbb{P}_1} \models "(\mathbf p, h)$ is active$"$ implies that $V^{\mathbb{P}_2} \models "(\mathbf p, h)$ is active$"$ (with $(\lambda, D)$-cc as defined in Requirement 1.16).

D) $\mathbf I$ and $\mathbf U$ are disjoint sets such that $<_{\mathbf U}$
is a fixed well ordering of $\mathbf U$ and $\mathbf{I}\cup \mathbf{U}$
is $\lambda^+$.

E) $|\mathbf{P}|\leq 2^{\lambda_2}$.

\textbf{Definition 2.2.A: }Let $\mathbf M=\mathbf M[\lambda_1,\lambda_2]$
be the collection of triples $\mathbf m=(\mathbf{q_m},M_{\mathbf m}, E_{\mathbf m}')$
such that the following conditions hold (we may replace the index
$\mathbf m$ by $\mathbf{q_m}$ or omit it completely when the context
is clear):

A) $\mathbf{q_m}\in \mathbf{K_P}$.

B) $M=M_{\mathbf m}\subseteq L_{\mathbf{q_m}}$ is a sub partial order.

C) For every $t\in M$, $w_t^0 \subseteq M$.

D) $E'=E_{\mathbf m}'$ is a relation on $L=L_{\mathbf{q_m}}$ satisfying
the following properties:

1. $E''=E'\restriction (L\setminus M)$ is an equivalence relation
on $L\setminus M$.

2. For every non $E''$-equivalent $s,t\in L\setminus M$, $s<_{L}t$
iff there is $r\in M$ such that $s<_L r<_L t$.

3. If $sE't$ then $s\notin M$ or $t\notin M$.

4. If $t\in L\setminus M$ then $\{s\in L : sE't\}=\{s\in L : tE's\}$.
We shall denote this set by $t/E'$.

5. If $s,t\in L\setminus M$ are $E''$-equivalent, then $s/E'=t/E'$.

6. If $t\in L\setminus M$ then $u_t^0 \subseteq t/E'$.

7. If $t\in L\setminus M$ then $|t/E'|\leq \lambda_2$.

8. $||M||\leq \lambda_1$.

9. $|w_t^0|\leq \lambda$ for every $t$.

E) In addition to the objects mentioned in definition 1.5, $\mathbf{q_m}$
includes a sequence $\bar{v}_{\mathbf m}=(v_{\mathbf{m},t} : t\in L_{\mathbf m})=(v_t : t\in L_{\mathbf m})$
such that for every $t\in L_{\mathbf m}$ we have:

1. $v_t\subseteq [u_t^0]^{\leq \lambda}$, $w_t^0 \in v_t$ and for
every $u\in v_t$, $u\cup w_t^0 \in v_t$ (recall that the $u_t^0$
and $w_t^0$ are part of the definition of $\mathbf{q_m}$ mentioned
in 1.5).

2. $v_t$ is closed under subsets.

3. If $t\in L_{\mathbf m}\setminus M_{\mathbf m}$ then $|v_t|\leq \lambda_2$.
If $t\in M_{\mathbf m}$ and $s\in L\setminus M$ then $|\{u\in v_t : u\cap (s/E_{\mathbf m}'')\neq \emptyset\}|\leq \lambda_2$.

4. For every $u\in v_t$, if $u\nsubseteq M_{\mathbf m}$ then there
is $s\in L_{\mathbf m}\setminus M_{\mathbf m}$ such that $u\subseteq s/E'$.

We shall now supply the final definition of the forcing (recalling
definition 1.8).

\textbf{Definition 2.2.B: }For $\mathbf m \in \mathbf M$ and the corresponding
iteration template $\mathbf{q_m}$ we shall define $\mathbb{P}_t=\mathbb{P}_{\mathbf m,t}, \underset{\sim}{\mathbb{Q}_t}$
and $\underset{\sim}{\eta_t}$ in the same way as in 1.10, except
that we replace (A)(2) and (C) with the following definition:

For every $s\in Dom(p)$ there is $\iota(p(s))<\lambda$, a collection
of sets $W_{p(s),\iota}\subseteq \xi_{p(s)}\leq \lambda$ $(\iota<\iota(p(s)))$,
a collection of $\lambda$-Borel functions $\mathbf{B}_{p(s),\iota}$
$(\iota<\iota(p(s)))$, $\lambda$-Borel functions $\mathbf{C}_{p(s)}$
and $\mathbf{B}_{p(s)}$ and an object $tr(p(s))$ such that the following
conditions hold:

A) $\xi=\xi_{p(s)}=\underset{\iota<\iota(p(s))}{\cup}W_{p(s),\iota}$.

B) $\mathbf{B}_{p(s)}(...,\underset{\sim}{\eta_{t_{\zeta}}}(a_{\zeta}),...)_{\zeta<\xi}=\mathbf{C}_{p(s)}(...,\mathbf{B}_{p(s),\iota}(...,\underset{\sim}{\eta_{t_{\zeta}}}(a_{\zeta}),...)_{\zeta \in W_{p(s),\iota}},...)_{\iota<\iota(p(s))}$
such that $t_{\zeta} \in u_s^0$ and $a_{\zeta} \in u_{t_{\zeta}}^{1}$
for every $\zeta \in W_{p(s), \iota}$ (for $\bold{C}_{p(s)}$ recall Definition 1.4(D)(6)(b)).
\\
$[$Following Definition 1.4(6)(B), $\mathbf{C}_{p(s)}$ really has the form $\mathbf{C}_{\bold{p}_s, h_s, \iota(p(s))}$, but we shall abuse the notation and denote it $\mathbf{C}_{p(s)}$. In addition, the definition implies that $tr(\mathbf{B}_{p(s),\iota}(...,\underset{\sim}{\eta_{t_{\zeta}}}(a_{\zeta}),...)_{\zeta \in W_{p(s),\iota}})$ is constant for $\iota < \iota(p(s))$, say $\eta_{p(s)}$, and so $\iota(p(s)) \leq |Dom(\eta_{p(s)})| ]$.

C) For every $\iota<\iota(p(s))$ there is $u\in v_s$ such that $\{t_{\zeta} : \zeta \in W_{p(s),\iota}\}\subseteq u$.

D) $p(s)=\mathbf{B}_{p(s)}(...,\underset{\sim}{\eta_{t_{\zeta}}}(a_{\zeta}),...)_{\zeta<\xi}$.
We may write $p(s)=(tr(p(s)), \mathbf{B}_{p(s)}(...,\underset{\sim}{\eta_{t_{\zeta}}}(a_{\zeta}),...)_{\zeta<\xi})$.

E) Recall that the parameter $h_s$ was defined in Definition 1.10(C). $\underset{\sim}{\mathbb{Q}_s}$ will be defined as the $\mathbb{P}_s$-name
of the subforcing of $\mathbb{Q}_{\mathbf{p}_s,h_s}$ with elements
of the form $\mathbf{C}(...,p_i,...)_{i<i(*)}$
such that each $p_i$ belongs to $\mathbb{Q}_{\mathbf{p}_s,h_s}^{V[\underset{\sim}{\eta_r} : r\in u]}$
for some $u\in v_{\mathbf m,s}$ and the $\lambda$-Borel function $\mathbf{C}=\mathbf{C}(...,p_i,...)_{i<i(*)}$
is into $\mathbb{Q}_{\mathbf{p}_s,h_s}$. This can be seen as a refinement of the previous Definition 1.10. The way that $\bold C$ is defined (as a function of conditions $p_i$) will play a role in the analysis of projections in Section 4, where incompatibility with a condition $p(s)$ will be reduced to incompatibility with some $\mathbf{B}_{p(s),\iota}(...,\underset{\sim}{\eta_{t_{\zeta}}}(a_{\zeta}),...)_{\zeta \in W_{p(s),\iota}}$.

F) For each $q_{s,\iota}=\mathbf{B}_{p(s),\iota}(...,\underset{\sim}{\eta_{t_{\zeta}}}(a_{\zeta}),...)_{\zeta \in W_{p(s),\iota}}$
there is an object $tr(q_{s,\iota})$ such that the range of $\mathbf{B}_{p(s),\iota}$
consists of conditions with trunk $tr(q_{s,\iota})$.

G) $tr(p(s))=\underset{\iota}{\cup}tr(q_{s,\iota})$ (so in particular, the $tr(q_{s,\iota})$'s are compatible).

H) $\Vdash_{\mathbb{P}_s} "\mathbf{C}_{p(s)}(...,\mathbf{B}_{p(s),\iota}(...,\underset{\sim}{\eta_{t_{\zeta}}}(a_{\zeta}),...)_{\zeta \in W_{p(s),\iota}},...)_{\iota<\iota(p(s))} \in \underset{\sim}{G}$ 
\\
$\leftrightarrow (\forall \iota<\iota(p(s))) \mathbf{B}_{p(s),\iota}(...,\underset{\sim}{\eta_{t_{\zeta}}}(a_{\zeta}),...)_{\zeta \in W_{p(s),\iota}} \in \underset{\sim}{G}$.
\\
\\
\textbf{Remark 2.2(A)}: The reader might wonder about the difference between the above definition and 1.10. In the main case, we will really be interested in iterating $\mathbb{Q}_t$ for $t\in M_{\bold m}$, where $M_{\bold m}$ might be an ordinal. In order to obtain the parallel of [JuSh:292], we would like to correct the iteration in order to have enough saturation while maintaining the well-foundedness of the iteration's underlying partial order. For this we add the "pseudo coordinates" grouped in classes  of the form $t/E_{\bold m}$. For $t\in M_{\bold m}$, we have in the definition the new sets $v_{\bold m, t}$ giving us the following difference between the iteration here and the one in Definition 1.10: In 1.10, $\underset{\sim}{\mathbb{Q}_t}$ is computed via $(\bold{p}_t, h_t)$ in $V[\underset{\sim}{\eta} \restriction u_t^0]$, while here it is the closure of the union of the forcings computed via $(\bold{p}_t, h_t)$ in $V[\underset{\sim}{\eta} \restriction v]$ for every $v\in v_{\bold m, t}$.

\textbf{Definition 2.3: }Let $L$ be a well founded partial order,
we shall define the depth of an element of $L$ and the depth of $L$
by induction as follows:

A) $dp(t)=dp_L(t)=\cup \{dp_L(s)+1 : s<_{L}t\}$.

B) $dp(L)=\cup \{dp_L(t)+1 : t\in L\}$.

\textbf{Definition 2.4: }Let $\mathbf m \in \mathbf M$ and let $L\subseteq L_{\mathbf m}$
be a sub-partial order, we shall define $\mathbf{n}=\mathbf{m} \restriction L$
as follows:

A) $\mathbf{q_n}=\mathbf{q_m} \restriction L$.

B) $M_{\mathbf n}=M_{\mathbf m}\cap L$.

C) $E_{\mathbf n}'=E_{\mathbf m}'\cap L\times L$.

D) For every $t\in L$ we define $v_{\mathbf{q_n},t}$ as $\{u\cap L : u\in v_{\mathbf{q_m},t}\}$.

Remark: If $M_{\mathbf m}\subseteq L$ then $\mathbf n \in \mathbf{M}[\lambda_1,\lambda_2]$.

\textbf{Definition 2.5: }Let $\mathbf n, \mathbf m \in \mathbf M$, a function
$f: L_{\mathbf m} \rightarrow L_{\mathbf n}$ is an isomorphism of $\mathbf m$
and $\mathbf n$ if the following conditions hold:

A) $f$ is an isomorphism of the partial orders $L_{\mathbf m}$ and
$L_{\mathbf n}$.

B) For every $t\in L_{\mathbf m}$, $\mathbf{p}_{\mathbf{q_m},t}=\mathbf{p}_{\mathbf{q_n},f(t)}$.

C) For every $t\in L_{\mathbf m}$, $f(u_{\mathbf m, t}^{0})=u_{\mathbf n, f(t)}^{0}$
and $\bar{u}_{\mathbf m, t}^{1}=\bar{u}_{\mathbf n, f(t)}^{1}$.

D) For every $t\in L_{\mathbf m}$, $f(w_{\mathbf m, t}^{0})=w_{\mathbf n, f(t)}^{0}$
and $\bar{w}_{\mathbf m, t}^{1}=\bar{w}_{\mathbf n, f(t)}^{1}$.

E) $M_{\mathbf n}=f(M_{\mathbf m})$.

F) For every $s,t \in L_{\mathbf m}$, $sE_{\mathbf m}'t$ if and only
if $f(s)E_{\mathbf m}'f(t)$.

G) For every $t\in L_{\mathbf m}$, if $((\mathbf{B}_{\mathbf m,t,b}, (s_t(b,\zeta),a_{t,b,\zeta} : \zeta< \xi(t,b)) : b\in I_{\mathbf{p}_{\mathbf{q_m},t}}^{0}) : t\in L_{\mathbf{q_m}})$
is as in 1.4(F) for $\mathbf m$, then $((\mathbf{B}_{\mathbf m,t,b}, (f(s_t(b,\zeta)),a_{t,b,\zeta} : \zeta< \xi(t,b)) : b\in I_{\mathbf{p}_{\mathbf{q_n},f(t)}}^{0}) : t\in L_{\mathbf{q_m}})$
is as in 1.4(F) for $\mathbf n$ at $f(t)$.

H) For every $t\in L_{\mathbf m}$, $u\in v_{\mathbf{q_m},t}$ if and
only if $f(u)\in v_{\mathbf{q_n},t}$.

\textbf{Definition 2.6: }We say that $\mathbf m, \mathbf n \in \mathbf M$
are equivalent if $\mathbf{q_m}=\mathbf{q_n}$.

\textbf{Definition 2.7: }A) Let $L$ be a partial order, we shall
denote by $L^+$ the partial order obtained from $L$ by adding a
new element $\infty$ such that $t<\infty$ for every $t\in L$.

B) Given $\mathbf m \in \mathbf M$ we shall denote by $\mathbb{P}_{\mathbf m}$
the limit of $(\mathbb{P}_t, \underset{\sim}{\mathbb{Q}_t} : t\in L_{\mathbf m})$
with support $<\lambda$, i.e. $\mathbb{P}_{\mathbf m,\infty}$. We
shall denote $\mathbb{P}_t$ by $\mathbb{P}_{\mathbf m,t}$ and similarly
for $\underset{\sim}{\mathbb{Q}_t}$. 

C) $p,q \in \mathbb{P}_{\mathbf m}$ are strongly compatible if $tr(p(s))R_{\mathbf{p}_{\bold{q_m}, s}}tr(q(s))$
for every $s\in Dom(p)\cap Dom(q)$.

D) Given an initial segment $L\subseteq L_{\mathbf m}$, let $\mathbb{P}_{\mathbf m,L}=\mathbb{P}_{\mathbf m} \restriction \{p\in \mathbb{P}_{\mathbf m} : Dom(p) \subseteq L\}$.

\textbf{Claim 2.8: }Let $\mathbf m \in \mathbf M$ and $s<t\in L_{\mathbf m}^+$.

A) If $p\in \mathbb{P}_s$ then $p\in \mathbb{P}_t$ and $p\restriction L_{<s}=p$.

B) If $p,q \in \mathbb{P}_s$ then $\mathbb{P}_s \models p\leq q$
iff $\mathbb{P}_t \models p\leq q$.

C) If $p\in \mathbb{P}_t$ then $p\restriction L_{<s} \in \mathbb{P}_s$
and $\mathbb{P}_{s} \models "p\restriction L_{<s} \leq p"$.

D) If $\mathbb{P}_t \models p\leq q$ then $\mathbb{P}_s \models p\restriction L_{<s} \leq q\restriction L_{<s}$.

E) If $p\in \mathbb{P}_t$, $q\in \mathbb{P}_s$ and $p\restriction L_{<s} \leq q\in \mathbb{P}_s$
then $p,q \leq q\cup (p\restriction (L_{<t}\setminus L_{<s}))\in \mathbb{P}_t$.

F)\textbf{ }If $s<t \in L_{\mathbf m}^+$ then $\mathbb{P}_s \lessdot \mathbb{P}_t$.

\textbf{Proof: }Should be clear. $\square$

\textbf{Claim 2.8': }Suppose that $\mathbf m \in \mathbf M$ and $L_1 \subseteq L_2 \subseteq L_{\mathbf m}$
are initial segments.

A) If $p\in \mathbb{P}_{L_1}$ then $p\in \mathbb{P}_{L_2}$ and $p\restriction L_1=p$.

B) If $p,q \in \mathbb{P}_{L_1}$ then $\mathbb{P}_{L_1} \models p\leq q$
iff $\mathbb{P}_{L_2} \models p\leq q$.

C) If $p\in \mathbb{P}_{L_2}$ then $p\restriction L_1 \in \mathbb{P}_{L_1}$.

D) If $p,q \in \mathbb{P}_{L_2}$ and $\mathbb{P}_{L_2} \models p\leq q$
then $\mathbb{P}_{L_1} \models p\restriction L_1 \leq q\restriction L_1$.

E) If $p\in \mathbb{P}_{L_2}$, $q\in \mathbb{P}_{L_1}$ and $\mathbb{P}_{L_1} \models "p\restriction L_1 \leq q"$
then $\mathbb{P}_{L_2} \models "p,q\leq q\cup (p\restriction (L_2 \setminus L_1))"$.

F) $\mathbb{P}_{L_1} \lessdot \mathbb{P}_{L_2}$.

\textbf{Proof: }Should be clear. $\square$

\textbf{Claim 2.9: }If $\mathbf m \in \mathbf M$, $p\in \mathbb{P}_{\mathbf m}$
and $s\in Dom(p)$, then there is a $\lambda$-Borel name of the form
$\mathbf{B}(...,TV(\underset{\sim}{\eta_{s_{\zeta}}}(a_{\zeta})=j_{\zeta}),...)_{\zeta<\xi(p,s)}$
such that $\mathbf{B}(...,TV(\underset{\sim}{\eta_{s_{\zeta}}}(a_{\zeta})=j_{\zeta}),...)_{\zeta<\xi(p,s)}[G_{\underset{\sim}{\mathbb{Q}_s}}]=true$
iff $p(s) \in G_{\underset{\sim}{\mathbb{Q}_s}}$ (where $TV(\underset{\sim}{\eta_{s_{\zeta}}}(a_{\zeta})=j_{\zeta})$ stands for the truth value of the statement "$\underset{\sim}{\eta_{s_{\zeta}}}(a_{\zeta})=j_{\zeta}$", so it's either $0$ or $1$). That is, membership in the generic set can be computed in a $\lambda$-Borel way that depends on the (partial) values of the generics.

\textbf{Proof: }Follows from the definition of forcing templates and
the assumptions of the previous chapter using the $\lambda^+$-c.c..
$\square$
\\
\\
As promised earlier, the properties of forcing templates will play an important role in the proof of the following:

\textbf{Claim 2.10: }Let $\mathbf m \in \mathbf M$ and let $L\subseteq L_{\mathbf m}$
be an initial segment.

A) a. If $s\in L$ then $\Vdash_{\mathbb{P}_L} \underset{\sim}{\eta_s}\in \underset{r\in I_{\mathbf{p_s}}^{1}}{\Pi}X_r$
where $X_r=\{x\in H(\lambda) : \nVdash_{\underset{\sim}{\mathbb{Q}_s}} \underset{\sim}{\eta_s}(r) \neq x\} \subseteq H(\lambda)$
(we may take $H(\lambda)^{I_{\mathbf{p}_s}^{1}}$ instead of this product).
\\
\\
b. Moreover, if $p\in  \mathbb{P}_{\bold m}$ and $a\in I^1_s$, then for some $q \in \mathbb{P}_{\bold m}$ above $p$ we have $s\in Dom(q)$, $a\in Dom(tr(q(s)))$ and $s\in Dom(p) \rightarrow \iota(p(s))=\iota(q(s))$.
\\
\\
c. The set $\{ p \in \mathbb{P}_m :$ for every $s\in Dom(p)$, $|\iota(p(s))| \leq |tr(p(s))| \}$ is dense in $\mathbb{P}_{\bold m}$. 
\\
\\
d. If $\lambda=\aleph_0$ and $h\in \omega^{\omega}$, then the set $\{ p\in \mathbb{P}_{\bold m}: s\in Dom(p) \rightarrow h(\iota(p(s)))<|tr(p(s))| \}$ is dense in $\mathbb{P}_{\bold m}$.

B) $\mathbb{P}_{\mathbf m} \models (\lambda,D)-cc$ (hence $\mathbb{P}_{\mathbf m}\models \lambda^+-c.c.$).

C) a. $\mathbb{P}_{\mathbf m,L}$ is $(<\lambda)$-strategically $0$-complete.
\\
\\
b. If $p$ is a function with $Dom(p) \in [L]^{<\lambda}$ such that $s\in Dom(p) \rightarrow$  $\Vdash_{\mathbb{P}_{\bold m, L_{<s}}} "p(s) \in \underset{\sim}{\mathbb{Q}_s}"$, then there is $q\in \mathbb{P}_{\bold m, L}$ such that $Dom(p) \subseteq Dom(q)$ and $q\restriction L_{<s} \Vdash_{\mathbb{P}_{\bold m, L_{<s}}} "p(s) \leq q(s)"$ for every $s\in Dom(p)$.

D) Let $t\in L_{\mathbf m}$, if $\Vdash_{\mathbb{P}_t} "\underset{\sim}{y} \in \underset{\sim}{\mathbb{Q}_t}"$
then there is a $\lambda$-Borel function $\mathbf{B}$, $\xi \leq \lambda$
and a sequence $(r_{\zeta} : \zeta<\xi)$ of members of $u_t^0$ such
that $\Vdash_{\mathbb{P}_t} "\underset{\sim}{y}=\mathbf{B}(...,\underset{\sim}{\eta_{r_{\zeta}}}(a_{\zeta}),...)_{\zeta<\xi}"$
for some $a_{\zeta} \in u_{r_{\zeta}}^{1}$. 

E) $\Vdash_{\mathbb{P}_{\mathbf m}} V[\underset{\sim}{\eta_t} : t\in L_{\mathbf m}]=V[\underset{\sim}{G}]$.

F) If $\Vdash_{\mathbb{P}_L} "\underset{\sim}{\eta} \in V^{\zeta}"$
for some $\zeta<\lambda$, then there is a $\lambda$-Borel function
$\mathbf{B}$, $\xi \leq \lambda$ and a sequence $(r_{\zeta} : \zeta<\xi)$
of members of $u_t^0$ such that $\Vdash_{\mathbb{P}_L} "\underset{\sim}{\eta}=\mathbf{B}(...,\underset{\sim}{\eta_{r_{\zeta}}}(a_{\zeta}),...)_{\zeta<\xi}"$
for suitable $a_{\zeta} \in u_{r_{\zeta}}^{1}$. 

\textbf{Proof: }The proof is by induction on $dp(L)$, simultaneously for all clauses (though naturally this is not needed in all cases).

A)a) Let $p\in \mathbb{P}_L$ and $a\in I_{\mathbf{p_s}}^1$ and let $p_1=p\restriction L_{<s}$,
then $p_1 \in \mathbb{P}_{L_{<s}}$.

Case 1: $s\notin Dom(p)$. There is $f\in \mathbf{T_{p_s}}$ such that
$a\in Dom(f)$, and by absoluteness (and parts (D)(2) and (E)(1) of
Definition 1.4, together with the remark below it), $\Vdash_{\mathbb{P}_{L_{<s}}} " V[\underset{\sim}{\eta} \restriction u^0_s] \models$ There is $q\in \mathbb{Q}_{\mathbf{p}_s,h_s}$
such that $f=tr(q)"$  (so this holds whether $(\bold{p}_s, h_s)$ is active or not). By the induction hypothesis for clause (D),
there are $p_1 \leq p_2 \in \mathbb{P}_{L_<s}$, a $\lambda$-Borel
function $\mathbf B$, $\xi \leq \lambda$, a sequence $(r_{\zeta} : \zeta<\xi)$
of members of $u_s^0$ and $\{a_{\zeta} : \zeta<\xi\} \subseteq I_s^1$ such that $p_2 \Vdash_{\mathbb{P}_{L_{<s}}} "V[\underset{\sim}{\eta} \restriction u_s^0] \models f=tr(\mathbf{B}(...,\underset{\sim}{\eta_{r_{\zeta}}}(a_{\zeta}),...)_{\zeta<\xi})"$.
Now define a condition $p_3 \in \mathbb{P}_{L}$ as follows: $Dom(p_3)=Dom(p_2) \cup Dom(p) \cup \{s\}$,
$p_3 \restriction Dom(p_2)=p_2$, $p_3 \restriction (Dom(p) \setminus Dom(p_2))=p\restriction (Dom(p) \setminus Dom(p_2))$
and $p_3(s)=(f,\mathbf{B}(...,\underset{\sim}{\eta_{r_{\zeta}}}(a_{\zeta}),...)_{\zeta<\xi})$.
Then $p,p_2 \leq p_3$ by absoluteness, 2.8 and the definition of the partial order.

Case 2: $s\in Dom(p)$. $p(s)$ has the form $\mathbf{C}_{p(s)}(...,\mathbf{B}_{p(s),\iota}(...,\underset{\sim}{\eta_{t_{\zeta}}}(a_{\zeta}),...)_{\zeta \in W_{p(s),\iota}},...)_{\iota<\iota(p(s))}$
as in definition 2.2(B). In $V^{\mathbb{P}_{L_{<s}}}$, $V[...,\eta_{t_{\zeta}},...]_{\zeta<\xi_{p(s)}}$
(see definition 2.2(B) for $\xi_{p(s)}$) is a subuniverse, $\mathbb{Q}={\mathbb{Q}_{\mathbf{p_s}, \underset{\sim}{h_s}}}^{V[...,\eta_{t_{\zeta}},...]_{\zeta<\xi_{p(s)}}}$
is well-defined (recall Definitions 1.5(E) and 1.10(C)) and $p(s)[...,\eta_{t_{\zeta}},...]_{\zeta<\xi_{p(s)}}$
is a condition in $\mathbb Q$ with trunk $tr(p(s))$. Let $G\subseteq \mathbb{P}_{L_{<s}}$ be generic over $V$ such that $p_1 \in G$, so in $V[G]$, $\underset{\sim}{\mathbb{Q}_{\bold{p}_s, h_s}}[G]$ is well-defined and contains $p(s)$. Therefore, by Definition 1.4(E)(2), there is $q$ above $p(s)$ with trunk $\eta$ such that $a\in Dom(\eta)$ and $tr(p(s)) \subseteq \eta$. For every $\iota <\iota(p(s))$, by absoluteness we have $V[\underset{\sim}{\eta}[G] \restriction \{t_{\zeta} : \zeta \in W_{p(s), \iota}\}] \models$ "$p_{\iota}^1:=\bold{B}_{p(s),\iota}(...,\underset{\sim}{\eta_{t_{\zeta}}}(a_{\zeta}),...)_{\zeta \in W_{p(s),\iota}}[G]$ and $\eta$ are compatible". Therefore, for every $\iota<\iota(p(s))$ there is some $p_{\iota}^2$ above $p_{\iota}^1$ with trunk $\eta$. Now let $p_2 \in \mathbb{P}_{L_{<s}}$ be a condition above $p_1$ forcing the above statements, and using $p_2$ and the $p^2_{\iota}$ we can get an extension of $p$ as required.  
\\
\\
A)b) By the proof of clause (a).
\\
\\
A)c) By the previous clause and by clause (C) (whose proof doesn't depend on the current clause).
\\
\\
A)d) By clause (b).

B) First we shall introduce a new definition: Let $L\subseteq L_{\mathbf m}$
be an initial segment, $\zeta$ an ordinal, $\gamma<\lambda$ and
let $L[<\zeta]=\{t\in L : dp(t)<\zeta\}$. 

Now suppose that $\{p_{\alpha} : \alpha<\lambda^+\} \subseteq \mathbb{P}_{L[<\zeta]}$. By clause (A)(c), wlog $\alpha<\lambda^+ \wedge (s\in Dom(p_{\alpha})) \rightarrow |\iota(p(s))|\leq |tr(p(s))|$, with strict inequality in case that $\lambda=\aleph_0$. Fix an enumeration $(s_{\epsilon} : \epsilon<\epsilon_*)$ of $L[<\zeta]$.
For every $\alpha<\lambda^+$, let $u_{\alpha}=\{\epsilon : s_{\epsilon} \in Dom(p_{\alpha})\}$.
For $s\in Dom(p_{\alpha})$, let $h_{s,\alpha}=tr(p_{\alpha}(s))$.
By 1.4(D)(7), there is $X_s \in D$ such that $(\alpha,\beta)\in X_s \rightarrow h_{s,\alpha}R_{\mathbf{p_s}}h_{s,\beta}$ (unless $\{ \alpha : s\in Dom(p_{\alpha})\}$ is bounded by some $\gamma<\lambda^+$, in which case we choose $X_s$ to be $(\lambda^+ \setminus \gamma) \times (\lambda^+ \setminus \gamma)$).
For every $\alpha<\lambda^+$, $|u_{\alpha}| = |Dom(p_{\alpha})|<\lambda$.
For every $\alpha<\lambda^+$, define $f_{\alpha}: u_{\alpha} \rightarrow \lambda$
by $f_{\alpha}(\zeta)=otp(u_{\alpha} \cap \zeta)$, and define $g: \underset{\alpha<\lambda^+}{\cup}u_{\alpha} \rightarrow D$
by $g(\xi)=X_{s_{\xi}}$. Let $X\in D$ be the set described in Hypothesis 0(b)(2)
for $(g, (f_{\alpha}, u_{\alpha} : \alpha<\lambda^+))$, we shall
prove that for $(\alpha,\beta) \in X$, $s\in Dom(p_{\alpha}) \cap Dom(p_{\beta}) \rightarrow tr(p_{\alpha}(s))R_{\mathbf{p_s}}tr(p_{\beta}(s))$.
Given $s\in Dom(p_{\alpha}) \cap Dom(p_{\beta})$, $s=s_{\xi}$ for
some $\xi \in u_{\alpha} \cap u_{\beta}$, so $(\alpha,\beta) \in g(\xi)=X_{s_{\xi}}$.
It follows that $tr(p_{\alpha}(s))R_{\mathbf{p_s}}tr(p_{\beta}(s))$.
For such $\alpha$ and $\beta$, it will suffice to find a common upper bound $p$. This will be done as
follows: Let $(s_{\epsilon} : \epsilon<\zeta)$ list $Dom(p_{\alpha}) \cap Dom(p_{\beta})$ in increasing order. For $\epsilon \leq \zeta$ let $L_{\epsilon}:=\{s : s<_L s_{\xi}$ for some $\xi < \epsilon \}$. We shall now choose $(p_{\epsilon}^*, q_{\epsilon}^*)$ by induction on $\epsilon$ such that: 
\\
\\
a. $\mathbb{P}_{\bold m, L_{\epsilon}} \models "p_{\epsilon}^* \leq q_{\epsilon}^*$.
\\
\\
b. $\mathbb{P}_{\bold m, L_{\epsilon}} \models "q_{\xi}^* \leq p_{\epsilon}^*$ for every $\xi<\epsilon"$.
\\
\\
c. $\mathbb{P}_{\bold m, L_{\epsilon}} \models "p_{\alpha} \restriction L_{\epsilon}, p_{\beta} \restriction L_{\epsilon}$ are below $p_{\epsilon}^*"$.
\\
\\
d. If $\xi < \epsilon$ and $s\in Dom(q_{\xi}^*) \setminus \underset{\iota < \xi}{\cup} Dom(q_{\iota}^*)$, then $(p_{\iota}^*(s), q_{\iota}^*(s) : \iota \in [\xi +1, \epsilon])$ is an initial segment of a play in the game $G_{\zeta+1}(q^*_{\xi}(s), \mathbb{Q}_{\bold{p}_s, h_s})$ according to a winning strategy of play $\bold I$.
\\
\\
There is a subtle issue that needs to be addressed: Recall that in Definition 1.4(D)(5) we didn't require $tr(q)=tr(p_1) \cup tr(p_2)$. However, this is not a problem. Arriving at $\epsilon$, let $u_0= \cup \{Dom(q_{\xi}^*) : \xi< \epsilon)\}$, so we can choose a function $p^1_{\epsilon}$ with domain $u_0$ such that, for every $s \in u_0$, $p^1_{\epsilon}(s)$ is a $\mathbb{P}_{\bold m, L_{<s}}$-name as required in clause (d). Note that by the definition of the strategic completeness game, if $G\subseteq \mathbb{P}_{\bold m, L_{<s}}$ is generic over $V$ and $V[G] \models "p^1_{\epsilon}(s) \leq r"$, then in $V[G]$, $r$ can be chosen by player $\bold I$ according to the winning strategy. Let $L_{<\epsilon}:=\underset{\xi < \epsilon}{\cup}L_{\xi}$, then by clause (C)(b) of the theorem, there is $p^2_{\epsilon} \in \mathbb{P}_{\bold m, L_{<\epsilon}}$ such that if $s\in Dom(p^1_{\epsilon})$ then $s\in Dom(p^2_{\epsilon})$ and $p^2_{\epsilon} \restriction L_{<s} \Vdash "p^1_{\epsilon}(s) \leq p^2_{\epsilon}(s)"$. The choice of $p^*_{\epsilon}$ is now split to cases:
\\
\\
1. $\epsilon=0$: Trivial.
\\
\\
2. $\epsilon$ is limit: In this case, we choose $p^*_{\epsilon}=p^2_{\epsilon}$. In order to show that $p^2_{\epsilon}$ satisfies clause (b), one can show by induction on $\xi \leq \epsilon$ that $q^*_{\xi} \restriction L_{<\xi} \leq p^2_{\epsilon} \restriction L_{<\xi}$, using at each step the choice of $p^1_{\epsilon}(s)$. Cases (c) and (d) then follow by the induction hypothesis and the choice of $p^2_{\epsilon}(s)$.
\\
\\
3. $\epsilon=\zeta+1$: In this case $p^2_{\epsilon} \in \mathbb{P}_{L_{\zeta}}$. If $s_{\zeta} \in Dom(p_{\alpha}) \cap Dom(p_{\beta})$, then we know that $\Vdash "p_{\alpha}(s_{\zeta}), p_{\beta}(s_{\zeta})$ have a common upper bound $\underset{\sim}{r_{\zeta}}"$. Let $p^3_{\epsilon} \in \mathbb{P}_{\bold m, L_{\zeta}}$ be a condition above $p^2_{\epsilon}$ that forces a value for $tr(\underset{\sim}{r_{\zeta}})$, and we can now choose a $p^*_{\epsilon}$ as required.

Finally, given $p_{\zeta}^*$ constructed above, the existence of a common upper bound for $p_{\alpha}$ and $p_{\beta}$ follows.

C) See, e.g., {[}Sh:587{]} for the preservation of $(<\lambda)$-strategic
completeness under $(<\lambda)$-support iterations, or just work as in the proof of clause (B) (but we rely neither on clause (A) nor on clause (B)). Note that we use 1.4(D)(8). As for (C)(b), it follows from strategic completeness for $\mathbb{P}_{L_{<s}}$ where $s<_L t$.

D) In order to avoid awkward notation, we shall write $\mathbf{B}(...,\underset{\sim}{\eta_{\zeta}},...)_{\zeta<\xi}$
instead of $\mathbf{B}(...,\underset{\sim}{\eta_{\zeta}}(a_{\zeta}),...)_{\zeta<\xi}$
for suitable $a_{\zeta}\in u_{\zeta}^{1}$.

The proof of the claim is by induction on $dp(t)$. Given $t\in L_{\mathbf m}$,
we shall prove the following claim by induction on $\zeta<\lambda^+$:

1. For every $p\in \mathbb{P}_t$ and $\zeta<\lambda^+$ such that
$p\Vdash_{\mathbb{P}_t} "\underset{\sim}{y} \in H_{\leq \lambda}(\mathbf{I} \cup \mathbf{U})\wedge rk(\underset{\sim}{y})<\zeta"$
there is a $\lambda$-Borel function $\mathbf{B}_p$ such that $p \Vdash_{\mathbb{P}_t} "\underset{\sim}{y}=\mathbf{B}_p(...,\underset{\sim}{\eta_{r_{\zeta}}},...)_{\zeta<\xi(p)}"$ 
with $r_{\zeta} \in u_t^0$ (for some $\xi(p)$ which is the length of the inputs for the function).

By a standard argument of definition by cases, this claim is equivalent
to:

2. For every antichain $I=\{p_i : i<i(*)\leq \lambda\}$ such that
$p_i \Vdash_{\mathbb{P}_t} "\underset{\sim}{y}\in H_{\leq \lambda}(\mathbf{I} \cup \mathbf{U})\wedge rk(\underset{\sim}{y})<\zeta"$
for every $i$, there is a $\lambda$-Borel function $\mathbf{B}_I$
such that for every $i<i(*)$, $p_i \Vdash_{\mathbb{P}_t} "\underset{\sim}{y}=\mathbf{B}_I(...,\underset{\sim}{\eta_{r_{\zeta}}},...)_{\zeta<\xi(p)}"$.

\textbf{Clause I: }$\zeta=0$.

There is nothing to prove in this case.

\textbf{Clause II: }$\zeta$ is a limit ordinal.

We shall prove the second version of the claim. For every $i<i(*)$,
let $\{p_{i,j} : j<j(i)\}$ be a maximal antichain above $p_i$ such
that every $p_{i,j}$ forces a value $\zeta_{i,j}$ to $rk(\underset{\sim}{y})$.
As $p\Vdash rk(\underset{\sim}{y})<\zeta$, for every $i, j$ we have
$\zeta_{i,j}<\zeta$. Hence, by the induction, for every $i,j$ there
is $\mathbf{B}_{i,j}(...,\underset{\sim}{\eta_{r_{\zeta,i,j}}},...)_{\zeta<\xi(i,j)}$
as required. For every $i<i(*)$ define a name $\underset{\sim}{\mathbf{B}_i}$
such that $\underset{\sim}{\mathbf{B}_i}[G]=\mathbf{B}_{i,j}(...,\underset{\sim}{\eta_{r_{\zeta,i,j}}},...)_{\zeta<\xi(i,j)}[G]$
iff $p_{i,j}\in G$ and $p_{i,j'}\notin G$ for every $j'<j$. Finally
define a name $\underset{\sim}{\mathbf{B}}$ such that $\underset{\sim}{\mathbf{B}}[G]=\underset{\sim}{\mathbf{B}_i}[G]$
iff $p_i \in G$ and for every $j<i$, $p_j \notin G$. Now let $i<i(*)$,
let $G$ be a generic set such that $p_i \in G$, then there is a
unique $j<j(i)$ such that $p_{i,j}\in G$. Therefore, $\underset{\sim}{\mathbf{B}}[G]=\underset{\sim}{\mathbf{B}_i}[G]=\mathbf{B}_{i,j}(...,\underset{\sim}{\eta_{r_{\zeta,i,j}}},...)_{\zeta<\xi(i,j)}[G]=\underset{\sim}{y}[G]$,
hence $p_i \Vdash_{\mathbb{P}_t} "\underset{\sim}{y}=\underset{\sim}{\mathbf{B}}"$.

\textbf{Clause III: }$\zeta=\epsilon+1$.

We shall prove the first version of the claim. Let $\{p_i : i<i(*)\}$
be a aximal antichain above $p$ such that for every $i$, $p_i \Vdash_{\mathbb{P}_t} "|\underset{\sim}{y}|=\mu_i"$
for some $\mu_i$. Therefore for every $i<i(*)$ there is a sequence
$(\underset{\sim}{y_{i,\alpha}} : \alpha<\mu_i)$ such that $p_i \Vdash_{\mathbb{P}_t} "\underset{\sim}{y}=\{\underset{\sim}{y_{i,\alpha}} : \alpha< \mu_i\}"$.
By the assumption, $p_i \Vdash_{\mathbb{P}_t} "rk(\underset{\sim}{y_{i,\alpha}})<\epsilon"$
for every $i$ and $\alpha$. By the induction hypothesis, for every
such $i$ and $\alpha$ there is $\mathbf{B}_{i,\alpha}(...,\eta_{r(\zeta,i,\alpha)},...)_{\zeta<\xi(i,\alpha)}$
as required for $\underset{\sim}{y_{i,\alpha}}$ and $p_i$. Hence
for every $i$ there is a name $\underset{\sim}{\mathbf{B}_i}$ as required
such that $p_i \Vdash_{\mathbb{P}_t} "\underset{\sim}{y}=\underset{\sim}{\mathbf{B}_i}"$.
Now define a name $\underset{\sim}{\mathbf{B}}$ such that $\underset{\sim}{\mathbf{B}}[G]=\underset{\sim}{\mathbf{B}_i}[G]$
iff $p_i \in G$ and as before we have $p\Vdash_{\mathbb{P}_t} "\underset{\sim}{y}=\underset{\sim}{\mathbf{B}}"$.

Remark: For $\zeta=1$, let $\{p_i : i<i(*)\}$ be a maximal antichain
above $p$ of elements that force a value for $\underset{\sim}{y}$
from $\mathbf{I} \cup \mathbf{U}$. Let $Y\subseteq \mathbf{I} \cup \mathbf{U}$
be the set of all such values (so $|Y| \leq \lambda$) and denote
by $a_i$ the value that $p_i$ forces to $p_i$. For every generic
$G$ that conatians $p$, $\underset{\sim}{y}[G]=a_i$ iff $p_i \in G$.
Therefore it's enough to show that for every $p_i$ there is a name
$\underset{\sim}{\mathbf{B}_i}$ of the right form such that $\underset{\sim}{\mathbf{B}_i}[G]=true$
iff $p_i \in G$. Therefore it's enough to show that the truth value
of $"p\in G"$ can be computed by a $\lambda-$Borel function as above,
so it's enough to compute the truth value $p\restriction \mathbb{P}_s \in G\cap \mathbb{P}_s$
for every $s<t$, which follows from the induction hypothesis.

E) By the assumption, for every $p\in \mathbb{P}_{\mathbf m}$ and $t\in Dom(p)$
there is a $\lambda-$Borel function $\mathbf{B}_{p,t}$ and a sequence
$(s_{\zeta} : \zeta<\xi(p,t))$ of members of $u_t^0$ such that for
every generic $G\subseteq \mathbb{P}_{\mathbf m}$ we have $\mathbf{B}_{p,t}(...,TV(\underset{\sim}{\eta_{s_{\zeta}}}(a_{\zeta})=j_{\zeta}),...)_{\zeta<\xi(p,t)}[G]=true$
if and only if $p(t) \in G_{\underset{\sim}{\mathbb{Q}_t}}$ (for
suitable $a_{\zeta}$ and $j_{\zeta}$). Therefore $p\in G$ iff $(\underset{t\in Dom(p)}{\wedge}\mathbf{B}_{p,t}(...,TV(\underset{\sim}{\eta_{s_{\zeta}}}(a_{\zeta})=j_{\zeta}),...)_{\zeta<\xi(p,i)})[G]=true$,
hence we can compute $G$ from $(\underset{\sim}{\eta_t} : t\in L_{\mathbf m})$.

F) Similar to the proof of (D). $\square$

\textbf{\large{}Properties of the $\mathbb{L}_{\lambda^+}-$closure}{\large\par}

\textbf{Definition 2.11: }A) Let $p\in \mathbb{P}_{\mathbf m}$, the
full support of $p$ will be defined as follows: for every $s\in Dom(p)$,
if $p(s)=(tr(p(s)), \mathbf{B}_{p(s)}(...,\eta_{t(s,\zeta)}(a_{\zeta}),...)_{\zeta<\xi(s)})$,
then the full support of $p$ will be defined as $fsupp(p):= \underset{s\in Dom(p)}{\cup}\{t(s,\zeta) : \zeta<\xi(s)\} \cup \{s\}$.

B) For $L\subseteq L_{\mathbf m}$ define $\mathbb{P}_{\mathbf m}(L):= \mathbb{P}_{\mathbf m} \restriction \{p\in \mathbb{P}_{\mathbf m} : fsupp(p) \subseteq L\}$
with the order inherited from $\mathbb{P}_{\mathbf m}$.

C) Let $L\subseteq L_{\mathbf m}$, for every $s\in L$, $j<\lambda$
and $a\in I_{\mathbf{p}_s}^{1}$ let $p_{s,a,j}\in \mathbb{P}_{\mathbf m}$
be a condition that represents $\underset{\sim}{\eta_s}(a)=j$ such
that $Dom(p_{s,a,j})=s$ and let $X_L:= \{p_{s,a,j} : s\in L, a\in I_{\mathbf{p}_s}^{1}, j<\lambda\}$.
\\
$[$Note that such $p_{s,a,j}$ exist by Definition 1.4(D)(9). It is not necessarily unique, but it can be chosen in $\mathbb{P}_{L_*}$ if $L_*$ is a minimal closed subset of $L_{\bold m}$ that contains $s$.$]$ 

D) For $L\subseteq L_{\mathbf m}$ define $\mathbb{P}_{\mathbf m}[L]:= \mathbb{L}_{\lambda^+}(X_L, \mathbb{P}_{\mathbf m})$
(see definition 1.13).

Remark: For $\mathbf m \in \mathbf M$ we may define the partial order
$\leq^*$ on $\mathbb{P}_{\mathbf m}$ by $p\leq^* q$ if and only if
$q\Vdash_{\mathbb{P}_{\mathbf m}} "p\in \underset{\sim}{G}"$. As $(\mathbb{P}_{\mathbf m} ,\leq^*)$
is equivalent to $(\mathbb{P}_{\mathbf m}, \leq)$, it's $(<\lambda)$-strategically
complete and satisfies $(\lambda,D)-cc$ and we may replace $(\mathbb{P}_{\mathbf m},\leq)$
by $(\mathbb{P}_{\mathbf m}, \leq^*)$.

\textbf{Claim 2.12: }Let $\mathbf m \in \mathbf M$ and $L\subseteq L_{\mathbf m}$.

A) $\mathbb{P}_{\mathbf m}\subseteq \mathbb{P}_{\mathbf m}[L_{\mathbf m}]$
is dense and $\mathbb{P}_{\mathbf m} \lessdot \mathbb{P}_{\mathbf m}[L_{\mathbf m}]$,
therefore they're equivalent.

B) $\mathbb{P}_{\mathbf m}[L_{\mathbf m}]$ is $(<\lambda)$ strategically
complete and satisfies $\lambda^+-cc$.

C) $\mathbb{P}_{\mathbf m}(L)\subseteq \mathbb{P}_{\mathbf m}$ and $\mathbb{P}_{\mathbf m}[L] \lessdot \mathbb{P}_{\mathbf m}[L_{\mathbf m}]$.

D) $\mathbb{P}_{\mathbf m}[L]$ is $(<\lambda)$-strategically complete
and satisfies $\lambda^+-cc$.

E) Let $G\subseteq \mathbb{P}_{\mathbf m}$ be generic, for each $t\in L$
let $\eta_t:= \underset{\sim}{\eta_t}[G]$ and let $G_L^+:= \{\psi \in \mathbb{P}_{\mathbf m}[L] : \psi[G]=true\}$,
then $G_L^+$ is $\mathbb{P}_{\mathbf m}[L]$-generic over $V$ and
$V[G_L^+]=V[\eta_t : t\in L]$.

F) For $L_1 \subseteq L_2 \subseteq L_{\mathbf m}$ we have $\mathbb{P}_{\mathbf m}(L_1)\subseteq \mathbb{P}_{\mathbf m}(L_2)$
(as partial orders) and $\mathbb{P}_{\mathbf m}[L_1] \lessdot \mathbb{P}_{\mathbf m}[L_2]$.

G) If $\mathbf m, \mathbf n \in \mathbf M$ are equivalent (recall Definition
2.6), then $\mathbb{P}_{\mathbf m}(L)=\mathbb{P}_{\mathbf n}(L)$ and
$\mathbb{P}_{\mathbf m}[L]=\mathbb{P}_{\mathbf n}[L]$.

H) Let $I$ be a $\lambda_2^+$-directed partial order and let $\{L_t : t\in I\}$
be a collection of subsets of $L_{\mathbf m}$ such that $s<_I t \rightarrow L_s \subseteq L_t$.
Let $L:= \underset{t\in I}{\cup}L_t$, then $\mathbb{P}_{\mathbf m}[L]=\underset{t\in I}{\cup}\mathbb{P}_{\mathbf m}[L_t]$.

\textbf{Proof: }A) By claim 2.9, there is a natural embedding of $\mathbb{P}_{\mathbf m}$
in $\mathbb{P}_{\mathbf m}[L_{\mathbf m}]$. For $p \in \mathbb{P}_{\mathbf m}$,
denote by $p^*$ its image under the embedding. Now let $\psi \in \mathbb{P}_{\mathbf m}[L_{\mathbf m}]$,
there is $p\in \mathbb{P}_{\mathbf m}$ such that $p\Vdash_{\mathbb{P}_{\mathbf m}} \psi[\underset{\sim}{G}]=true$,
therefore for every generic $G\subseteq \mathbb{P}_{\mathbf m}$, if
$p^*[G]=true$ then $p\in G$ and $\psi[G]=true$, hence $\mathbb{P}_{\mathbf m}[L_{\mathbf m}]\models \psi \leq p^*$
and $\mathbb{P}_{\mathbf m}$ is dense in $\mathbb{P}_{\mathbf m}[L_{\mathbf m}]$.

B) By 2.10 (B+C), $\mathbb{P}_{\mathbf m}$ has these properties, and
by the clause (A), $\mathbb{P}_{\mathbf m}[L_{\mathbf m}]$ has these
properties too.

C) The first part is by the definition of $\mathbb{P}_{\mathbf m}(L)$.
For the second part, first note that, by definition, $\mathbb{P}_{\mathbf m}[L] \subseteq \mathbb{P}_{\mathbf m}[L_{\mathbf m}]$
as partial orders. Now note that if $\psi, \phi \in \mathbb{P}_{\mathbf m}[L]$
are compatible in $\mathbb{P}_{\mathbf m}[L_{\mathbf m}]$, then $\psi \wedge \phi \in \mathbb{P}_{\mathbf m}[L]$
is a common upper bound, so $\phi$ and $\psi$ are compatible in
$\mathbb{P}_{\mathbf m}[L]$ iff they're compatible in $\mathbb{P}_{\mathbf m}[L_{\mathbf m}]$.
Therefore if $I\subseteq \mathbb{P}_{\mathbf m}[L]$ is a maximal antichain,
then $I$ remains an antichain in $\mathbb{P}_{\mathbf m}[L_{\mathbf m}]$.
Furthermore, it's a maximal antichain in $\mathbb{P}_{\mathbf m}[L_{\mathbf m}]$:
Suppose towards contradiction that $\phi \in \mathbb{P}_{\mathbf m}[L_{\mathbf m}]$
is incompatible with all members of $I$. Let $\psi= \underset{\theta \in I}{\wedge} \neg \theta$.
As $I$ is an antichain in $\mathbb{P}_{\mathbf m}[L_{\mathbf m}]$ which
satisfies the $\lambda^+-c.c.$, we have that $|I| \leq \lambda$.
As $\phi \in \mathbb{P}_{\mathbf m}[L_{\mathbf m}]$, there is a generic
$G\subseteq \mathbb{P}_{\mathbf m}$ such that $\phi[G]=true$. As $\phi$
is incompatible with all elements of $I$, it follows that $\theta[G]=false$
for all $\theta \in I$. Therefore, $\psi \in \mathbb{P}_{\mathbf m}[L]$.
But $\psi$ is clearly incompatible with all members of $I$, a contradiction.
Therefore, $\mathbb{P}_{\mathbf m}[L]\lessdot \mathbb{P}_{\mathbf m}[L_{\mathbf m}]$.

D) By (B) and (C).

E) We shall first show that $G_{L_{\mathbf m}}^+$ is $\mathbb{P}_{\mathbf m}[L_{\mathbf m}]$-generic.
$G_{L_{\mathbf m}}^{+}$ is downward-closed, by the definition of $G_{L_{\mathbf m}}^+$
and of the order of $\mathbb{P}_{\mathbf m}[L_{\mathbf m}]$. If $\psi, \phi \in G_{L_{\mathbf m}}^{+}$
then $(\psi \wedge \phi)[G]=true$, hence $\psi \wedge \phi \in G_{L_{\mathbf m}}^{+}$,
so $G_{L_{\mathbf m}}^{+}$ is directed. Now let $I=\{ \psi _i : i<i(*)\} \subseteq \mathbb{P}_{\mathbf m}[L_{\mathbf m}]$
be a maximal antichain and let $J=\{p \in \mathbb{P}_{\mathbf m} : (\exists i<i(*))(p \Vdash "\psi_i[\underset{\sim}{G}]=true")\}$.
If $J$ is predense in $\mathbb{P}_{\mathbf m}$, then there is $q\in J\cap G$.
Let $i<i(*)$ such that $q\Vdash_{\mathbb{P}_{\mathbf m}} "\psi_i[\underset{\sim}{G}]=true"$,
then $\psi_i[G]=true$ hence $\psi_i \in G_{L_{\mathbf m}}^{+} \cap I$.
Suppose towards contradiction that $J$ is not predense and let $q\in \mathbb{P}_{\mathbf m}$
be incompatible with all members of $J$, so $q\Vdash_{\mathbb{P}_{\mathbf m}} "\psi_i[\underset{\sim}{G}]=false"$
for every $i<i(*)$. $i(*) \leq \lambda$ (as $\mathbb{P}_{\mathbf m} \models \lambda^+-c.c.$),
hence $\psi_* := \underset{i<i(*)}{\wedge}(\neg \psi_i) \in \mathbb{L}_{\lambda}(X_{L_{\mathbf m}})$
and $\psi_* \in \mathbb{L}_{\lambda}(X_{L_{\mathbf m}}, \mathbb{P}_{\mathbf m})$.
Obviously, $\psi_*$ is incompatible with the members of $I$, contradicting
our maximality assumption. Therefore we proved that $G_{L_{\mathbf m}}^{+}$
is $\mathbb{P}_{\mathbf m}[L_{\mathbf m}]$-generic.

Now let $L\subseteq L_{\mathbf m}$, then $G_{L_{\mathbf m}}^{+} \cap \mathbb{P}_{\mathbf m}[L]$
is $\mathbb{P}_{\mathbf m}[L]$-generic and $G_{L_{\mathbf m}}^{+} \cap \mathbb{P}_{\mathbf m}[L]=G_L^+$.

We shall now prove that $V[G_L^+]=V[\eta_t : t\in L]$.
We need to show that $G_L^+$ can be computed from $\{\eta_t : t\in L\}$.
Let $p_{s,a,j} \in X_L$, then $p_{s,a,j} \in G_L^+$ iff $p_{s,a,j}[G]=true$
iff $\underset{\sim}{\eta_s}[G](a)=j$. Therefore we can compute $G_L^+ \cap X_L$
and $G_L^+$ from $\{\underset{\sim}{\eta_s}[G] : s\in L\}$. As $\underset{\sim}{\eta_s}[G](a)=j$
iff $p_{s,a,j} \in G_L^+$, we can compute $\{\underset{\sim}{\eta_s}[G] : s\in L\}$
in $V[G_L^+]$, therefore $V[G_L^+]=V[\underset{\sim}{\eta_s} : s\in L]$.

F) If $fsupp(p) \subseteq L_1$ then $fsupp(p) \subseteq L_2$, hence
$p\in \mathbb{P}_{\mathbf m}(L_1) \rightarrow p\in \mathbb{P}_{\mathbf m}(L_2)$,
and by the definition of the order, $\mathbb{P}_{\mathbf m}(L_1) \subseteq \mathbb{P}_{\mathbf m}(L_2)$
as partial orders. For the second claim, first note that $\mathbb{P}_{\mathbf m}[L_1] \subseteq \mathbb{P}_{\mathbf m}[L_2]$
as partial orders. Now assume that $I \subseteq \mathbb{P}_{\mathbf m}[L_1]$
is a maximal antichain. By (C), $I$ is a maximal antichain in $\mathbb{P}_{\mathbf m}[L_{\mathbf m}]$,
hence in $\mathbb{P}_{\mathbf m}[L_2]$. Therefore $\mathbb{P}_{\mathbf m}[L_1] \lessdot \mathbb{P}_{\mathbf m}[L_2]$.

G) If $\mathbf m$ and $\mathbf n$ are equivalent, then $\mathbf{q_n}=\mathbf{q_m}$,
hence $\mathbb{P}_{\mathbf m}=\mathbb{P}_{\mathbf n}$, $\mathbb{P}_{\mathbf n}(L)=\mathbb{P}_{\mathbf m}(L)$
and $\mathbb{P}_{\mathbf m}[L]=\mathbb{P}_{\mathbf n}[L]$ for every $L$.

H) For every $t\in I$, $L_t \subseteq L$, therefore $\mathbb{P}_{\mathbf m}[L_t]\subseteq \mathbb{P}_{\mathbf m}[L]$,
so $\underset{t\in I}{\cup}\mathbb{P}_{\mathbf m}[L_t]\subseteq \mathbb{P}_{\mathbf m}[L]$.
In the other direction, suppose that $\psi \in \mathbb{P}_{\mathbf m}[L]$
is generated by the atoms $\{p_{s(i),a(i),j(i)} : s(i)\in L, a(i)\in I_{\mathbf{p}_{s(i)}}^{1}, j(i), i<\lambda\}$.
Recall that $\lambda \leq \lambda_2 \leq \lambda_2^+$, hence there
is $i(*)\in I$ such that $\{s(i) : i<\lambda\}\subseteq L_{i(*)}$,
therefore $\psi \in \mathbb{P}_{\mathbf m}[L_{i(*)}]$, so $\mathbb{P}_{\mathbf m}[L] \subseteq \underset{i\in I}{\cup} \mathbb{P}_{\mathbf m}[L_i]$.
$\square$

\textbf{\large{}Operations on members of $\mathbf M$}{\large\par}

We shall define a partial order $\leq_{\mathbf M}=\leq$ on $\mathbf M$
as follows:

\textbf{Definition 2.13: }Let $\mathbf m, \mathbf n \in \mathbf M$, we
shall write $\mathbf m \leq \mathbf n$ if:

A) $L_{\mathbf m}\subseteq L_{\mathbf n}$.

B) $M_{\mathbf m}=M_{\mathbf n}$ (yes, equal).

C) $\mathbf{q_m} \leq_{\mathbf{K_P}} \mathbf{q_n}$.

D) $u_{\mathbf{q_m},t}^{0}=u_{\mathbf{q_n},t}^{0}$ for every $t\in L_{\mathbf m}\setminus M_{\mathbf m}$.

E) $t/E_{\mathbf n}'=t/E_{\mathbf m}'$ for every $t\in L_{\mathbf m} \setminus M_{\mathbf m}$.

F) If $t\in M_{\mathbf m}$ then $v_{\mathbf{q_m},t}=\{u\cap L_{\mathbf m} : u\in v_{\mathbf{q_n},t}\}$,
if $t\in L_{\mathbf m}\setminus M_{\mathbf m}$ then $v_{\mathbf{q_n},t}=v_{\mathbf{q_m},t}$.

G) If $t\in M_{\mathbf m}$ then $\{u\in v_{\mathbf m,t} : u\subseteq M_{\mathbf m}\}=\{u\in v_{\mathbf n,t} : u\subseteq M_{\mathbf m}\}$.

H) If $t\in M_{\mathbf m}$ and $s\in L_{\mathbf m}\setminus M_{\mathbf m}$
then $\{u\in v_{\mathbf m,t} : u\subseteq s/E_{\mathbf m}'\}=\{u\in v_{\mathbf n,t} : u\subseteq s/E_{\mathbf n}'\}$.

\textbf{Definition 2.14: }Let $(\mathbf{m}_{\alpha} : \alpha<\delta)$
be an increasing sequence of elements of $\mathbf M$ with respect to
$\leq_{\mathbf M}$, we shall define the union $\mathbf n= \underset{\alpha<\delta}{\cup}\mathbf{m}_{\alpha}$
as follows:

A) $M_{\mathbf n}=M_{\mathbf{m}_{\alpha}}$ $(\alpha<\delta)$.

B) $E_{\mathbf n}'=\underset{\alpha<\delta}{\cup}E_{\mathbf{m}_{\alpha}}'$.

C) $\mathbf{q_n}$ will be defined as follows: 

1. $L_{\mathbf n}=\underset{\alpha<\delta}{\cup}L_{\mathbf{m}_{\alpha}}$.

2. For every $t\in L_{\mathbf{q_n}}$, $\mathbf{p}_{\mathbf{q_n},t}=\mathbf{p}_{\mathbf{q_{m_{\alpha}}},t}$
(for $\alpha<\delta$ such that $t\in L_{\mathbf{m}_{\alpha}}$).

3. For every $t\in L_{\mathbf n}$, $u_{\mathbf{q_n},t}^{0}=\cup \{u_{\mathbf{q_{m_{\alpha}}},t}^{0} : \alpha<\delta \wedge t\in L_{\mathbf{m}_{\alpha}} \}$
and $\bar{u}_{\mathbf{q_n},t}^{1}=\underset{\alpha<\delta}{\cup}\bar{u}_{\mathbf{q}_{\mathbf{m}_{\alpha}},t}^{1}$.

4. For every $t\in L_{\mathbf n}$, $w_{\mathbf{q_n},t}^{0}=\cup \{w_{\mathbf{q_{m_{\alpha}}},t}^{0} : \alpha<\delta \wedge t\in L_{\mathbf{m}_{\alpha}}\}$
and $\bar{w}_{\mathbf{q_n},t}^{1}=\underset{\alpha<\delta}{\cup}\bar{w}_{\mathbf{q}_{\mathbf{m}_{\alpha}},t}^{1}$.

5. $((\mathbf{B}_{t,b}, (s_t(b,\zeta), a_{t,b,\zeta}) : \zeta<\xi(t,b)) : b\in I_{\mathbf{p}_t}^{0}) : t\in L_{\mathbf{q_n}}))$
will be defined naturally as the union of the sequences corresponding
to the sequence of the $\mathbf{m}_{\alpha}$'s.

6. $v_{\mathbf{q_n},t}=\underset{\alpha<\delta}{\cup} v_{\mathbf{q}_{\mathbf{m}_{\alpha}},t}$
for every $t\in L_{\mathbf n}$.

It's easy to see that the union is a well defined member of $\mathbf M$.

\textbf{Claim 2.15: }Let $(\mathbf{m}_{\alpha} : \alpha<\delta)$ and
$\mathbf n$ be as above, then $\mathbf n \in \mathbf M$ and $\mathbf{m}_{\alpha}\leq \mathbf n$
for every $\alpha<\delta$.

\textbf{Proof: }It's straightforward to verify that $\mathbf{m}_{\alpha}\leq \mathbf n$
for every $\alpha<\delta$. $\square$

\textbf{Defintion and claim 2.16 (Amalgamation): }Suppose that

A) $\mathbf{m}_0, \mathbf{m}_1, \mathbf{m}_2 \in \mathbf M$.

B) $\mathbf{m}_0 \leq \mathbf{m}_l$ $(l=1,2)$.

C) $L_{\mathbf{m}_1}\cap L_{\mathbf{m}_2}=L_{\mathbf{m}_0}$.

We shall define the amalgamation $\mathbf m$ of $\mathbf{m}_1$ and $\mathbf{m}_2$
over $\mathbf{m}_0$ as follows:

1. $E_{\mathbf m}'=E_{\mathbf{m}_1}'\cup E_{\mathbf{m}_2}'$.

2. $M_{\mathbf m}=M_{\mathbf{m}_0}$.

$\mathbf{q_m}$ will be defined as follows: 

3. $L_{\mathbf m}$ is the minimal partial order containing $L_{\mathbf{m}_1}$
and $L_{\mathbf{m}_2}$.

4. For every $t\in L_{\mathbf m}$, $\mathbf{p}_{\mathbf{q_m},t}=\mathbf{p}_{\mathbf{q_{m_l}},t}$
provided that $t\in L_{\mathbf{m}_l}$.

5. $u_{\mathbf{q_m},t}^{0}=u_{\mathbf{q_{m_1}},t}^{0} \cup u_{\mathbf{q_{m_2}},t}^{0}$
(where $u_{\mathbf{q_{m_l}},t}^{0}=\emptyset$ if $t\notin L_{\mathbf{m}_l}$).

6. $w_{\mathbf{q_m},t}^{0}=w_{\mathbf{q_{m_1}},t}^{0} \cup w_{\mathbf{q_{m_2}},t}^{0}$
(where $w_{\mathbf{q_{m_l}},t}^{0}=\emptyset$ if $t\notin L_{\mathbf{m}_l}$).

7. $\bar{u}_{\mathbf{q_m},t}^{1}=\bar{u}_{\mathbf{q_{m_1}},t}^{1} \cup \bar{u}_{\mathbf{q_{m_2}},t}^{1}$,
$\bar{w}_{\mathbf{q_m},t}^{1}=\bar{w}_{\mathbf{q_{m_1}},t}^{1} \cup \bar{w}_{\mathbf{q_{m_2}},t}^{1}$,
i.e. coordinatewise union (similarly to 5+6, if $t\notin L_{\mathbf{m}_l}$,
the corresponding sequence will be defined as the empty sequence).

8. For $t\in L_{\mathbf{m}_1}\cup L_{\mathbf{m}_2}$, the $\lambda$-Borel
functions from 1.5(E) will be defined in the same way as in the case
of $\mathbf{m}_1$ and $\mathbf{m}_2$.

9. If $t\in L_{\mathbf{m}_0}$ then $v_{\mathbf{q_m},t}=v_{\mathbf{q_{m_1}},t} \cup v_{\mathbf{q_{m_2}},t}$.
If $t\in L_{\mathbf{m}_l}\setminus L_{\mathbf{m}_0}$ $(l=1,2)$ then
$v_{\mathbf{q_m},t}=v_{\mathbf{q_{m_l}},t}$.

\textbf{Claim 2.16: $\mathbf m$ }is well defined, $\mathbf m \in \mathbf M$
and $\mathbf{m}_1, \mathbf{m}_2\leq \mathbf m$.

\textbf{Proof: }Straightforward. $\square$

Remark: The amalgamation of a set $\{\mathbf{m}_i : 1\leq i<i(*)\}$
over $\mathbf{m}_0$ can be defined naturally as in 2.16.

\textbf{\large{}Existentially closed iteration parameters}{\large\par}

Given $\mathbf m \in \mathbf M$, we would like to construct extensions
$\mathbf m \leq \mathbf n$ which are, in a sense, existentially closed.

\textbf{Definition and Observation 2.17 }A) Let $\mathbf m \in \mathbf M$,
$L\subseteq L_{\mathbf m}$, we shall define the relative depth of $L$
as follows: $dp_{\mathbf m}^{*}(L):=\cup \{dp_{M_{\mathbf m}}(t)+1 : t\in L\cap M_{\mathbf m}\}$ (so this is $dp_{M_{\bold m}}(L\cap M_{\bold m})$).

B) For $\gamma \in Ord$ we shall define $\mathbf{M}_{\gamma}^{ec}$
as the set of elements $\mathbf m \in \mathbf M$ satisfying the following
property: Let $\mathbf m \leq \mathbf{m}_1 \leq \mathbf{m}_2$, $L_{\mathbf{m}_l, \gamma}^{dp}:=\{t\in L_{\mathbf{m}_l} : sup\{dp_{M_{\mathbf m}}(s) : s<t, s\in M_{\mathbf m}\}< \gamma \}$
$(l=1,2)$, then $\mathbb{P}_{\mathbf{m}_1}(L_{\mathbf{m}_1, \gamma}^{dp})\lessdot \mathbb{P}_{\mathbf{m}_2}(L_{\mathbf{m}_2,\gamma}^{dp})$. Note that in this case we have $\mathbb{P}_{\mathbf{m}_1}(L)=\mathbb{P}_{\mathbf{m}_2}(L)$ for
every $L\subseteq L_{\mathbf{m}_1,\gamma}^{dp}$.

C) $\mathbf{M}_{ec}$ will be defined as the collection of elements
$\mathbf m \in \mathbf M$ such that $\mathbf m \in \mathbf{M}_{\gamma}^{ec}$
for every $\gamma \in Ord$.

\textbf{Observation}: $\mathbf m \in \mathbf{M}_{ec}$ if and only if
$\mathbb{P}_{\mathbf{n}_1}\lessdot \mathbb{P}_{\mathbf{n}_2}$ for every
$\mathbf m \leq \mathbf{n}_1 \leq \mathbf{n}_2$.

Proof: Suppose that $\mathbf m \in \mathbf{M}_{\gamma}^{ec}$ for every
$\gamma$ and $\mathbf m \leq \mathbf{m}_1 \leq \mathbf{m}_2$. Choose some
$\gamma'$ such that $\gamma'>dp_{M_{\mathbf{m}_l}}(s)$ for every $s\in M_{\mathbf{m}_l}$
$(l=1,2)$ and let $\gamma=\gamma'+1$. Obviously $L_{\mathbf{m}_l}=L_{\mathbf{m}_l,\gamma}^{dp}$
$(l=1,2)$, so $\mathbb{P}_{\mathbf{m}_1}=\mathbb{P}_{\mathbf{m}_1}(L_{\mathbf{m}_1,\gamma}^{dp})\lessdot \mathbb{P}_{\mathbf{m}_2}(L_{\mathbf{m}_2,\gamma}^{dp})=\mathbb{P}_{\mathbf{m}_2}$.
In the other direction, suppose that $\mathbb{P}_{\mathbf{m}_1}\lessdot \mathbb{P}_{\mathbf{m}_2}$
for every $\mathbf{m} \leq \mathbf{m}_1 \leq \mathbf{m}_2$ and let $\gamma \in Ord$.
As $L_{\mathbf{m}_l, \gamma}^{dp}$ is an initial segment of $L_{\bold{m}_l}$, it follows that   $\mathbb{P}_{\mathbf{m}_l}(L_{\mathbf{m}_l, \gamma}^{dp})\lessdot \mathbb{P}_{\mathbf{m}_l}$
$(l=1,2)$, and we have $\mathbb{P}_{\mathbf{m}_1}(L_{\mathbf{m}_1,\gamma}^{dp}) \lessdot \mathbb{P}_{\mathbf{m}_1} \lessdot \mathbb{P}_{\mathbf{m}_2}$
and $\mathbb{P}_{\mathbf{m}_2}(L_{\mathbf{m}_2,\gamma}^{dp})\lessdot \mathbb{P}_{\mathbf{m}_2}$.
Note that $L_{\mathbf{m}_1,\gamma}^{dp} \subseteq L_{\mathbf{m}_2,\gamma}^{dp}$,
so $\mathbb{P}_{\mathbf{m}_1}(L_{\mathbf{m}_1,\gamma}^{dp})\subseteq \mathbb{P}_{\mathbf{m}_2}(L_{\mathbf{m}_2,\gamma}^{dp})$
and it follows that every maximal antichain in $\mathbb{P}_{\mathbf{m}_1}(L_{\mathbf{m}_1,\gamma}^{dp})$
is a maximal antichain in $\mathbb{P}_{\mathbf{m}_2}(L_{\mathbf{m}_2,\gamma}^{dp})$,
so $\mathbf{m} \in \mathbf{M}_{\gamma}^{ec}$. $\square$

\textbf{Definition 2.18: }Let $\chi$ be a cardinal, we shall denote
by $\mathbf{M}_{\chi} (\mathbf{M}_{\leq \chi})$ the collection of members
$\mathbf m \in \mathbf M$ such that $|L_{\mathbf m}|=\chi$ $(|L_{\mathbf m}| \leq \chi)$.

\textbf{Claim 2.19: }Let $2^{\lambda_2}\leq \chi$ and $\mathbf m \in \mathbf M_{\leq \chi}$,
then there is $\mathbf m \leq \mathbf n \in \mathbf{M}_{\chi}$ such that
$\mathbf n \in \mathbf{M}_{ec}$.

\textbf{Proof: }Denote by $C=C_{\mathbf m}$ the collection of elements
$\mathbf n \in \mathbf M$ such that:

1. $\mathbf m \restriction M_{\mathbf m} \leq \mathbf n$ (recall Definition 2.4).

2. $L_{\mathbf n}\setminus M_{\mathbf m}=t/E_{\mathbf n}''$ for some $t$.

Definition: Let $\mathbf{n}_1, \mathbf{n}_2 \in C$, a function $h: L_{\mathbf{n}_1}\rightarrow L_{\mathbf{n}_2}$
is called a strong isomorphism of $\mathbf{n}_1$ onto $\mathbf{n}_2$
If:

1. $h$ is an isomorphism of $\mathbf{n}_1$ onto $\mathbf{n}_2$.

2. $h$ is the identity on $M_{\mathbf m}$.

Definition: Let $R=R_{\mathbf m}$ be the following equivalence relation
on $C_{\mathbf m}$:

$\mathbf{n}_1 R 
\mathbf{n}_2$ iff there is a strong isomorphism of $\mathbf{n}_1$ onto $\mathbf{n}_2$.

We shall now estimate the number of $R$-equivalence relations:

1. As $|L_{\mathbf n}|\leq \lambda_2$ for every $\mathbf n \in C$, once
we fix $M_{\mathbf n}$ there are at most $2^{\lambda_2}$ possible
isomorphism types of $(L_{\mathbf n},\leq_{L_{\mathbf n}})$ over $M_{\mathbf n}$.

2. Given such $L_{\mathbf n}$, there are at most $2^{\lambda_2}$ possible
forcing templates from $\mathbf{P}$.

3. For every $\mathbf n \in C$ there is $t$ such that $|L_{\mathbf n}|=|L_{\mathbf n}\setminus M_{\mathbf m}|+|M_{\mathbf m}|=|t/E_{\mathbf n}''|+|M_{\mathbf m}|\leq \lambda_2$
(recalling definition 2.2.A), hence $|\mathcal{P}(L_{\mathbf n})|\leq 2^{\lambda_2}$
and for every $t\in L_{\mathbf n}$ there are at most $2^{\lambda_2}$
possible values for $u_{\mathbf{q_n},t}^{0}$ and $w_{\mathbf{q_n},t}^{0}$.

4. For every $t$, $\bar{u}_{\mathbf{q_n},t}^{1}$ is a function assigning
for each $s$ a member of $\mathcal{P}(I_{s}^{1})$, so we have at
most $(2^{|\mathbf I|})^{|L_{\mathbf n}|}\leq 2^{(|\mathbf{I}|+\lambda_2)}$
possible functions. Similar argument applies to $\bar{w}_{\mathbf{q_n},t}^{1}$
as well.

Therefore there are at most $2^{\lambda_2}$ $R-$equivalence classes.
Let $(\mathbf{n}_{\alpha} : \alpha<2^{\lambda_2})$ list all such classes.
For every $\alpha<2^{\lambda_2}$ we shall choose the sequence $(\mathbf{n}_{\alpha}^{i} : i<\chi)$
such that each $\mathbf{n}_{\alpha}^{i}$ is obtained from $\mathbf{n}_{\alpha}$
by the changing the names of the elements in $L_{\mathbf{n}_{\alpha}}\setminus M_{\mathbf m}$
such that the new sets are pairwise disjoint and also disjoint to
$L_{\mathbf m}$ (for $i<\chi$). For every $i$ there is $t_{\alpha,i}$
such that $t_{\alpha,i}/E_{\mathbf{n}_{\alpha}^{i}}''=L_{\mathbf{n}_{\alpha}^{i}}\setminus M_{\mathbf m}$
and $t_{\alpha,i}/E_{\mathbf{n}_{\alpha}^{i}}''\cap t_{\alpha,j}/E_{\mathbf{n}_{\alpha}^{j}}''=\emptyset$.
Now let $\mathbf n$ be the amalgamation of $\{\mathbf m\} \cup \{\mathbf{n}_{\alpha}^{i} : i<\chi, \alpha<2^{\lambda_2}\}$
over $\mathbf m \restriction M_{\mathbf m}$. Obviously, $\mathbf n \in \mathbf M_{\chi}$.

Suppose now that $\mathbf n \leq \mathbf{n}_1 \leq \mathbf{n}_2$. Let $\mathcal{F}$
be the collection of functions $f$ such that for some $L_1,L_2 \subseteq L_{\mathbf{n}_2}$:

a. $Dom(f)=L_1$, $Ran(f)=L_2$.

b. $M_{\mathbf m}=M_{\mathbf n}\subseteq L_1 \cap L_2$.

c. $|L_l\setminus M_{\mathbf m}|\leq \lambda_2$ $(l=1,2)$.

d. $t/E_{\mathbf{n}_2}\subseteq L_l$ for every $t\in L_l \setminus M_{\mathbf m}$.

e. $f$ is the identity on $M_{\mathbf m}$.

f. $f$ is an isomorphism of $\mathbf{n}_2 \restriction L_1$ onto $\mathbf{n}_2 \restriction L_2$.

Claim 1: Let $f\in \mathcal{F}$, $L'\subseteq L_{\mathbf{n}_1}$, $L''\subseteq L_{\mathbf{n}_2}$
such that $|L'|+|L''|\leq \lambda_2$, then there is $g\in \mathcal{F}$
such that $f\subseteq g$, $L'\subseteq Dom(g)$ and $L''\subseteq Ran(g)$.

Proof: WLOG $L'\cap Dom(f)=\emptyset=L'' \cap Ran(f)$ and $|L'|=|L''|=\lambda_2$.
Let $(a_i : i<\lambda_2)$ and $(b_j : j<\lambda_2)$ list $L'$ and
$L''$, respectively. For $b\in L_{\bold{n}_2} \setminus M_{\bold m}$, let $B_b:=(b/E_{\bold{n}_2}') \cup M_{\bold m}$, then $\bold m \restriction M_{\bold m} \leq \bold{n}_2 \restriction B_b$, $\bold{n}_2 \restriction B_b \in C$ and $\bold{n}_2 \restriction B_b \leq \bold{n}_2$. We shall construct by induction on $i<\lambda_2$
an increasing continuous sequence of functions $f_i \in \mathcal{F}$
such that $g:=\cup f_i$ will give the desired function of the claim.

I. $i=0$: $f_0:=f$.

II. $i$ is a limit ordinal: $f_i:=\underset{j<i}{\cup}f_j$.

III. $i=2j+1$: By the "WLOG" above,
$L'' \cap M_{\mathbf m}=\emptyset$, hence $b_j \in L_{\mathbf{n}_2}\setminus M_{\mathbf m}$.
Therefore it follows that $\mathbf{m} \restriction M_{\mathbf m} \leq \mathbf{n}_2 \restriction B_{b_j}$,
hence $\mathbf{n}_2 \restriction B_{b_j}\in C$. Let $\mathbf{n}_{\alpha}$
be the representative of the $R$-equivalence class of $\mathbf{n}_2 \restriction B_{b_j}$.
By $\mathcal{F}$'s definition, $|Dom(f_{2j})|\leq \lambda_2$. Since
$\mathbf{n}$ is the result of an amalgamation that includes $\mathbf{n}_{\alpha}^{i}$
$(i<\chi)$, each $\mathbf{n}_{\alpha}^{i}$ is $R$-equivalent to $\mathbf{n}_{\alpha}$
and $\lambda_2 <\chi$, it follows that for some $i<\chi$, $L_{\mathbf{n}_{\alpha}^{i}}\setminus M_{\mathbf m} \cap Dom(f_{2j})=\emptyset$.
Since $\mathbf{n}_2 \restriction B_{b_j}R \mathbf{n}_{\alpha}^{i}$,
there is a strong isomorphism $h$ from $\mathbf{n}_2 \restriction L_{\mathbf{n}_{\alpha}^{i}}=\mathbf{n}_{\alpha}^{i}$
onto $\mathbf{n}_2 \restriction B_{b_j}$. Therefore $f_i:=f_{2j}\cup h$
is a well defined function, $b_j \in Ran(f_i)$ and $f_{2j}\subseteq f_i$.
We shall now show that $f_i \in \mathcal{F}$: conditions a, b, c
and e are obviously satisfied. If $t\in L_{\mathbf{n}_{\alpha}^{i}}\setminus M_{\mathbf m}$,
then $t/E_{\mathbf n}=t/E_{\mathbf{n}_2}$ (as $\mathbf{n} \leq \mathbf{n}_2$)
and $t/E_{\mathbf n}=t/E_{\mathbf{n}_{\alpha}^{i}}$. Therefore $t/E_{\mathbf{n}_2}=t/E_{\mathbf{n}_{\alpha}^{i}}\subseteq L_{\mathbf{n}_{\alpha}^{i}}\subseteq Dom(f_i)$.
Similarly, if $t\in b_j/E_{\mathbf{n}_2}''$ then $t/E_{\mathbf{n}_2}=b_j/E_{\mathbf{n}_2}\subseteq Ran(f_i)$,
hence condition d is satisfied. It remains to show that $f_i$ is
an isomorphism of $\mathbf{n}_2 \restriction Dom(f_i)$ onto $\mathbf{n}_2 \restriction Ran(f_i)$.
Note that $b_j/E_{\mathbf{n}_2}'' \cap Ran(f_{2j})=\emptyset$ (as we may assume WLOG that $b_j \notin Ran(f_{2j})$), hence $f_i$ is an order
preserving bijection, as a union of two such functions (that are identified
on $M_{\mathbf m}$). It's easy to check that $f_i$ is as required.

$IV.$ $i=2j+2:$ Similar to the previous case, ensuring that $a_j \in Dom(f_{2j+1})$.

As $\mathcal{F}$ is closed to increasing unions of length $\lambda_2$,
$g:=\underset{i<\lambda_2}{\cup}f_i \in \mathcal{F}$ is as required,
hence we're done proving claim 1.

Denote $L_{\gamma}:=\{s\in L_{\mathbf{n}_2} : dp_{\mathbf{n}_2}(s)<\gamma\}$
(so $L_{\mathbf{n}_2}=L_{|L_{\mathbf{n}_2}|^+}$).

Claim 1(+): Let $f\in \mathcal{F}$, $L'\subseteq L_{\mathbf{n}_2}$
such that $|L'|\leq \lambda_2$ and $Ran(f)\subseteq L_{\mathbf{n}_1}$,
then there exists $g\in \mathcal{F}$ such that $f\subseteq g$, $L' \subseteq Dom(g)$
and $Ran(g) \subseteq L_{\mathbf{n}_1}$.

Proof: Repeat the proof of claim 1 (in particular, stage $2j+2$).
Note that at each stage we add a set of the form $L_{\mathbf{n}_{\alpha}^{i}}$
to the range. As $L_{\mathbf{n}_{\alpha}^{i}} \subseteq L_{\mathbf n} \subseteq L_{\mathbf{n}_1}$
and $Ran(f) \subseteq L_{\mathbf{n}_1}$, it follows that $Ran(g) \subseteq L_{\mathbf{n}_1}$.

Claim 2: Let $g\in \mathcal{F}$, then $g(Dom(g)\cap L_{\gamma})=Ran(g) \cap L_{\gamma}$.

Proof: By induction on $\gamma$.

Claim 3: Given $g\in \mathcal{F}$ and $\gamma<|L_{\mathbf{n}_2}|^+$,
the map $\hat{g}$ is an isomorphism of $\mathbb{P}_{\mathbf{n}_2}(Dom(g)\cap L_{\gamma})$
onto $\mathbb{P}_{\mathbf{n}_2}(Ran(g) \cap L_{\gamma})$ where $\hat{g}$
is defined as follows: Given $p\in \mathbb{P}_{\mathbf{n}_2}(Dom(g) \cap L_{\gamma})$,
$\hat{g}(p)=q$ has the domain $g(Dom(p))$, and for every $g(s)\in Dom(q)$,
$q(g(s))=(tr(p(s)), \mathbf{B}_{p(s)}(...,\eta_{g(t_{\zeta})}(a_{\zeta}),...)_{\zeta<\xi})$
where $p(s)=(tr(p(s)), \mathbf{B}_{p(s)}(...,\eta_{t_{\zeta}}(a_{\zeta}),...)_{\zeta<\xi})$.

Proof: Given $g\in \mathcal{F}$, by the previous claim $g$ is a
bijection from $Dom(g)\cap L_{\gamma}$ onto $Ran(g)\cap L_{\gamma}$.
As $g\in \mathcal{F}$, it's order preserving and the information
of $\mathbf{q_{n_2}}\restriction (Dom(g)\cap L_{\gamma})$ is preserved.
Hence clearly $\hat{g}$ is an isomorphism from $\mathbb{P}_{\mathbf{n}_2}(Dom(g)\cap L_{\gamma})$
onto $\mathbb{P}_{\mathbf{n}_2}(Ran(g)\cap L_{\gamma})$.

Claim 4: $\mathbb{P}_{\mathbf{n}_2}(L_{\gamma}\cap L_{\mathbf{n}_1})\lessdot \mathbb{P}_{\mathbf{n}_2}(L_{\gamma})$.

Proof: By induction on $\gamma$. Arriving at stage $\gamma$, note
that $\mathbb{P}_{\mathbf{n}_2}(L_{\gamma}\cap L_{\mathbf{n}_1})\subseteq \mathbb{P}_{\mathbf{n}_2}(L_{\gamma})$
(as partial orders). Suppose that $p_1,p_2 \in \mathbb{P}_{\mathbf{n}_2}(L_{\gamma}\cap L_{\mathbf{n}_1})$
are compatible in $\mathbb{P}_{\mathbf{n}_2}(L_{\gamma})$, and let
$q\in \mathbb{P}_{\mathbf{n}_2}(L_{\gamma})$ be a common uppper bound.
Since $|fsupp(p_1)|,|fsupp(p_2)|\leq \lambda$, there is $L'$ such
that $fsupp(p_1)\cup fsupp(p_2)\subseteq L'\subseteq (L_{\gamma}\cup L_{\mathbf{n}_1})$,
$|L'|\leq \lambda_2$ and $L'$ is $E_{\mathbf{n}_2}$-closed. Therefore
$p_1,p_2 \in \mathbb{P}_{\mathbf{n}_2}(L')$. Similarly, there is $L''\subseteq L_{\gamma}$
such that $|L''|\leq \lambda_2$, $fsupp(q)\cup L' \subseteq L''$
and $L''$ is $E_{\mathbf{n}_2}$-closed, hence $q\in \mathbb{P}_{\mathbf{n}_2}(L'')$.
Let $f$ be the identity function on $L_1=L_2=\cup \{t/E_{\mathbf{n}_2} : t\in L'\setminus M_{\mathbf m}\}$.
Note that $|L_i|\leq \lambda_2$ $(i=1,2)$ and $f\in \mathcal{F}$.
Let $L_1':=\cup \{t/E_{\mathbf{n}_2} : t\in L''\setminus M_{\mathbf m}\}$,
then $|L_1'|\leq \lambda_2$, hence by claim $1(+)$, there is $g\in \mathcal{F}$
such that $f\subseteq g$ such that $L_1'\subseteq Dom(g)$ and $Ran(g) \subseteq L_{\mathbf{n}_1}$.
As $fsupp(q)\cup fsupp(p_1)\cup fsupp(p_2)\subseteq Dom(g)\cap L_{\gamma}$,
we have $p_1, p_2, q\in \mathbb{P}_{\mathbf{n}_2}(Dom(g)\cap L_{\gamma})$,
hence $\hat{g}(p_1),\hat{g}(p_2),\hat{g}(q)\in \mathbb{P}_{\mathbf{n}_2}(Ran(g)\cap L_{\gamma})$
(in particular, $\hat{g}(q),\hat{g}(p_1),\hat{g}(p_2)$ are well defined).
By the choice of $g$, $\hat{g}(p_1)=p_1$ and $\hat{g}(p_2)=p_2$.
By claim $3$, $\mathbb{P}_{\mathbf{n}_2}(Ran(g)\cap L_{\gamma})\models p_1,p_2 \leq \hat{g}(q)$.
As $Ran(g)\subseteq L_{\mathbf{n}_1}$, $\hat{g}(q)\in \mathbb{P}_{\mathbf{n}_2}(L_{\gamma}\cap L_{\mathbf{n}_1})$,
hence $p_1$ and $p_2$ are compatible in $\mathbb{P}_{\mathbf{n}_2}(L_{\gamma}\cap L_{\mathbf{n}_1})$.
Therefore, if $I\subseteq \mathbb{P}_{\mathbf{n}_2}(L_{\gamma}\cap L_{\mathbf{n}_1})$,
then $I$ remains an antichaim in $\mathbb{P}_{\mathbf{n}_2}(L_{\gamma})$.

Suppose now that $I\subseteq \mathbb{P}_{\mathbf{n}_2}(L_{\gamma}\cap L_{\mathbf{n}_1})$
is a maximal antichain, and suppose towards contradiction that $q\in \mathbb{P}_{\mathbf{n}_2}(L_{\gamma})$
is incompatible with all members of $I$. We can show by induction on $\gamma$ that $\mathbb{P}_{\mathbf{n}_1}(L_{\gamma}\cap L_{\mathbf{n}_1})=\mathbb{P}_{\mathbf{n}_2}(L_{\gamma}\cap L_{\mathbf{n}_1})$.
Since $L_{\gamma}\cap L_{\mathbf{n}_1}$ is an initial segment of $L_{\mathbf{n}_1}$,
$\mathbb{P}_{\mathbf{n}_1}(L_{\gamma}\cap L_{\mathbf{n}_1})=\mathbb{P}_{\mathbf{n}_1 \restriction (L_{\gamma} \cap L_{\mathbf{n}_1})} \lessdot \mathbb{P}_{\mathbf{n}_1}$,
hence $\mathbb{P}_{\mathbf{n}_2}(L_{\gamma}\cap L_{\mathbf{n}_1})\models \lambda^+-c.c.$
and $|I|\leq \lambda \leq \lambda_2$. Let $(p_i : i<\lambda_2)$
enumerate $I$'s members, then there is $L'\subseteq L_{\gamma}\cap L_{\mathbf{n}_1}$
such that $|L'|\leq \lambda_2$ and $\underset{i<\lambda_2}{\cup}fsupp(p_i)\subseteq L'$,
hence $I\subseteq \mathbb{P}_{\mathbf{n}_2}(L')$. Define $L''$ and
choose $f$ and $g$ as before. Again, $\hat{g}: \mathbb{P}_{\mathbf{n}_2}(L_{\gamma}\cap Dom(g))\rightarrow \mathbb{P}_{\mathbf{n}_2}(L_{\gamma}\cap Ran(g))$
is an isomorphism, $I\cup \{q\} \subseteq Dom(\hat{g})$ and $\hat{g}$
is thee identity on $I$. Hence $\hat{g}(q)$ is incompatible in $\mathbb{P}_{\mathbf{n}_2}(L_{\gamma}\cap Ran(g))$
with all members of $I$. As before, $\hat{g}(q)\in \mathbb{P}_{\mathbf{n}_2}(L_{\gamma}\cap L_{\mathbf{n}_1})$,
therefore, in order to get a contradiction, it's enough to show that
$\hat{g}(q)$ is incompatible in $\mathbb{P}_{\mathbf{n}_2}(L_{\gamma}\cap L_{\mathbf{n}_1})$
with all members of $I$. Suppose that for some $p\in I$, $r\in \mathbb{P}_{\mathbf{n}_2}(L_{\gamma}\cap L_{\mathbf{n}_1})$
we have $p,\hat{g}(q)\leq r$. Since $g^{-1}\in \mathcal{F}$, as
in previous arguments, there is $g^{-1}\subseteq h\in \mathcal{F}$
such that $\hat{h}(r), \hat{h}(\hat{g}(q))$ are well-defined and
$\hat{h}(p)=p$, $\hat{h}(\hat{g}(q))=q$. Hence $p$ and $q$ are
compatible in $\mathbb{P}_{\mathbf{n}_2}(L_{\gamma}\cap Ran(h))$ and
therefore in $\mathbb{P}_{\mathbf{n}_2}(L_{\gamma})$, contradicting
the assumption. This proves claim $4$.

Claim $5$: $\mathbb{P}_{\mathbf{n}_1}\lessdot \mathbb{P}_{\mathbf{n}_2}$.

Proof: By the previous claim, for $\gamma=|L_{\mathbf{n}_2}|^+$ we
get $\mathbb{P}_{\mathbf{n}_2}(L_{\mathbf{n}_1})=\mathbb{P}_{\mathbf{n}_2}(L_{\gamma}\cap L_{\mathbf{n}_1})\lessdot \mathbb{P}_{\mathbf{n}_2}(L_{\gamma})=\mathbb{P}_{\mathbf{n}_2}$.
We can show by induction on $\delta$ that $\mathbb{P}_{\mathbf{n}_1}(L_{\delta}\cap L_{\mathbf{n}_1})=\mathbb{P}_{\mathbf{n}_2}(L_{\delta}\cap L_{\mathbf{n}_1})$,
hence for $\delta=\gamma$ we get $\mathbb{P}_{\mathbf{n}_1}\lessdot \mathbb{P}_{\mathbf{n}_2}$.
This proves claim $2.19$. $\square$

The following observation will be useful throughout the rest of this
paper:

\textbf{Observation 2.20: }Let $\mathbf n \in \mathbf{M}_{ec}$ and $\mathbf{n} \leq \mathbf{n}_1 \leq \mathbf{n}_2$,
then for every $L\subseteq L_{\mathbf{n}_1}$, $\mathbb{P}_{\mathbf{n}_1}[L]=\mathbb{P}_{\mathbf{n}_2}[L]$.

\textbf{Proof}: $\mathbf{n}_1 \leq \mathbf{n}_2$, hence for $L\subseteq L_{\mathbf{n}_1}$,
the set $X_L$ in definition $2.11(c)$ is the same for $\mathbf{n}_1$
and $\mathbf{n}_2$. Let $\psi \in \mathbb{L}_{\lambda}(X_L)$, since
$\mathbb{P}_{\mathbf{n}_1}\lessdot \mathbb{P}_{\mathbf{n}_2}$, there
is a generic set $G\subseteq \mathbb{P}_{\mathbf{n}_2}$ such that $\psi[G]=true$
iff there is a generic set $H\subseteq \mathbb{P}_{\mathbf{n}_1}$ such
that $\psi[H]=true$. Similarly, if an implication of the form $"\psi[G]=true \rightarrow \phi[G]=true"$
holds for every generic $G\subseteq \mathbb{P}_{\mathbf{n}_2}$, then it holds for every generic $H\subseteq \mathbb{P}_{\mathbf{n}_1}$, and vice
versa. Therefore, $\mathbb{P}_{\mathbf{n}_1}[L]=\mathbb{P}_{\mathbf{n}_2}[L]$.
$\square$

\textbf{Claim 2.21: }Suppose that

A) $\mathbf{m}_1, \mathbf{m}_2 \in \mathbf{M}_{ec}$.

B) $M_l=M_{\mathbf{m}_l}$ $(l=1,2)$.

C) $h: M_1 \rightarrow M_2$ is an isomorphism from $\mathbf{m}_1 \restriction M_1$
onto $\mathbf{m}_2 \restriction M_2$.

then $\mathbb{P}_{\mathbf{m}_1}[M_1]$ is isomorphic to $\mathbb{P}_{\mathbf{m}_2}[M_2]$.

\textbf{Proof}: WLOG $M_1=M_2$ (denote this set by $M$), $L_{\mathbf{m}_1}\cap L_{\mathbf{m}_2}=M$
and $h$ is the identity. Let $\mathbf{m}_0:= \mathbf{m}_1 \restriction M=\mathbf{m}_2 \restriction M$,
then $\mathbf{m}_0 \leq \mathbf{m}_1,\mathbf{m}_2$ and $L_{\mathbf{m}_0}=L_{\mathbf{m}_1}\cap L_{\mathbf{m}_2}$,
therefore, by 2.16, there is $\mathbf m \in \mathbf M$ such that $\mathbf m$
is the amalgamation of $\mathbf{m}_1$ and $\mathbf{m}_2$ over $\mathbf{m}_0$
and $\mathbf{m}_1, \mathbf{m}_2 \leq \mathbf m$. By the definition of $\mathbf{M}_{ec}$,
as $\mathbf{m}_l \in \mathbf{M}_{ec}$, $\ \mathbf{m}_l \leq \mathbf m$ $(l=1,2)$ and $M\subseteq L_{\mathbf{m}_l}$ $(l=1,2)$, it follows that $\mathbb{P}_{\mathbf{m}_1}[M]=\mathbb{P}_{\mathbf m}[M]=\mathbb{P}_{\mathbf{m}_2}[M]$.
$\square$

\textbf{\large{}The Corrected Iteration}{\large\par}

We shall now describe how to correct an iteration $\mathbb{P}_{\mathbf m}$
in order to obtain the desired iteration for the main result.

\textbf{Definition 2.22: }Let $\mathbf m \in \mathbf M$, we shall define
the corrected iteration $\mathbb{P}_{\mathbf m}^{cr}$ as $\mathbb{P}_{\mathbf n}[L_{\mathbf m}]$
for $\mathbf m \leq \mathbf n \in \mathbf{M}_{ec}$ (we'll show that $\mathbb{P}_{\mathbf m}^{cr}$
is indeed well-defined). For $L\subseteq L_{\mathbf m}$, define $\mathbb{P}_{\mathbf m}^{cr}[L]:=\mathbb{P}_{\mathbf n}[L]$
for $\mathbf n$ as above.

\textbf{Claim 2.23 }A) $\mathbb{P}_{\mathbf m}^{cr}[L]$ is well-defined
for every $\mathbf m \in \mathbf M$ and $L\subseteq L_{\mathbf m}$.

B) $\mathbb{P}_{\mathbf m}^{cr}[M_{\mathbf m}]$ is well-defined for every
$\mathbf m \in \mathbf M$ and depends only on $\mathbf m \restriction M_{\mathbf m}$.

C) If $\mathbf m \leq \mathbf n$ then $\mathbb{P}_{\mathbf m}^{cr} \lessdot \mathbb{P}_{\mathbf n}^{cr}$.

D) If $\mathbf m \leq \mathbf n$ and $L\subseteq L_{\mathbf m}$, then
$\mathbb{P}_{\mathbf m}^{cr}[L]=\mathbb{P}_{\mathbf n}^{cr}[L]$.

\textbf{Proof: }A) By claim 2.19, there is $\mathbf m \leq \mathbf n \in \mathbf{M}_{ec}$,
so it's enough to show that the definition does not depend on the
choice of $\mathbf n$. Given $\mathbf{n}_1, \mathbf{n}_2 \in \mathbf{M}_{ec}$
such that $\mathbf m \leq \mathbf{n}_l$, we have to show that $\mathbb{P}_{\mathbf{n}_1}[L_{\mathbf m}]=\mathbb{P}_{\mathbf{n}_2}[L_{\mathbf m}]$.
WLOG $L_{\mathbf{n}_1}\cap L_{\mathbf{n}_2}=L_{\mathbf m}$. Let $\mathbf n$
be the amalgamation of $\mathbf{n}_1, \mathbf{n}_2$ over $\mathbf m$.
Since $\mathbf{n}_1 \in \mathbf{M}_{ec}$, $\mathbf{n}_1 \leq \mathbf{n}_1 \leq \mathbf n$
and $L_{\mathbf m}\subseteq L_{\mathbf{n}_1}$, we get $\mathbb{P}_{\mathbf{n}_1}[L_{\mathbf m}]=\mathbb{P}_{\mathbf n}[L_{\mathbf m}]$.
Similarly, $\mathbb{P}_{\mathbf{n}_2}[L_{\mathbf m}]=\mathbb{P}_{\mathbf n}[L_{\mathbf m}]$,
therefore, $\mathbb{P}_{\mathbf{n}_1}[L_{\mathbf m}]=\mathbb{P}_{\mathbf{n}_2}[L_{\mathbf m}]$.
The argument for $\mathbb{P}_{\mathbf{m}}^{cr}[L]$ is similar.

B) Suppose that $\mathbf{m}_1 \restriction M_{\mathbf{m}_1}$ is isomorphic
to $\mathbf{m}_2 \restriction M_{\mathbf{m}_2}$ and choose $\mathbf{n}_l$
$(l=1,2)$ such that $\mathbf{m}_l \leq \mathbf{n}_l \in \mathbf{M}_{ec}$.
Now, $\mathbf{m}_1 \restriction M_{\mathbf{m}_1}=\mathbf{n}_1 \restriction M_{\mathbf{m}_1}$
is isomorphic to $\mathbf{n}_2 \restriction M_{\mathbf{m}_2}=\mathbf{m}_2 \restriction M_{\mathbf{m}_2}$,
hence by claim 2.21, $\mathbb{P}_{\mathbf{n}_1}[M_{\mathbf{m}_1}]$ is
isomorphic to $\mathbb{P}_{\mathbf{n}_2}[M_{\mathbf{m}_2}]$. Moreover,
the proof of 2.21 shows that if $\mathbf{m}_1 \restriction M_{\mathbf{m}_1}=\mathbf{m}_2 \restriction M_{\mathbf{m}_2}$,
then $\mathbb{P}_{\mathbf{n}_1}[M_{\mathbf{m}_1}]=\mathbb{P}_{\mathbf{n}_2}[M_{\mathbf{m}_2}]$,
therefore $\mathbb{P}_{\mathbf{m}_1}^{cr}[M_{\mathbf{m}_1}]=\mathbb{P}_{\mathbf{m}_2}^{cr}[M_{\mathbf{m}_2}]$.

C) Choose $\mathbf n \leq \mathbf{n}_*$ such that $\mathbf{n}_* \in \mathbf{M}_{ec}$,
then $\mathbb{P}_{\mathbf n}^{cr}=\mathbb{P}_{\mathbf{n}_*}[L_{\mathbf n}]$.
As $\mathbf m \leq \mathbf{n}_*$, it follows that $\mathbb{P}_{\mathbf m}^{cr}=\mathbb{P}_{\mathbf{n}_*}[L_{\mathbf m}]$.
By $2.12(F)$, $\mathbb{P}_{\mathbf m}^{cr}=\mathbb{P}_{\mathbf{n}_*}[L_{\mathbf m}]\lessdot \mathbb{P}_{\mathbf{n}_*}[L_{\mathbf n}]=\mathbb{P}_{\mathbf n}^{cr}$.

D) Choose $(\mathbf m \leq)\mathbf n \leq \mathbf{n}_* \in \mathbf{M}_{ec}$,
then by definition we get $\mathbb{P}_{\mathbf m}^{cr}[L]=\mathbb{P}_{\mathbf{n}_*}[L]=\mathbb{P}_{\mathbf n}^{cr}[L]$.
$\square$

\textbf{\large{}The main result}{\large\par}

\textbf{Definition 2.24: }Let $\mathbf q$ be a $(\lambda,D)$-iteration
template such that $|L_{\mathbf q}| \leq \lambda_1$ and $|w_t^0| \leq \lambda$
for every $t\in L_{\mathbf q}$. 

We call $\mathbf m=\mathbf{m_q} \in \mathbf M$ the iteration parameter
derived from $\mathbf q$ if:

a. $\mathbf{q_m}=\mathbf q$.

b. $M_{\mathbf m}=L_{\mathbf q}$.

c. $E_{\mathbf m}'=\emptyset$.

d. For every $t\in L_{\mathbf q}$, $v_t=[u_t^0]^{\leq \lambda}$.

\textbf{Definition 2.25:} Given $\mathbf m \in \mathbf M$, we define
the forcing notions $(\mathbb{P}_t' : t\in L_{\mathbf m} \cup \{\infty\})=(\mathbb{P}_{\mathbf m, t}' : t\in L_{\mathbf m} \cup \{\infty\})$
as follows: Fix $\mathbf m \leq \mathbf n \in \mathbf{M}_{ec}$ and let
$\mathbb{P}_t':=\mathbb{P}_{\mathbf n}[\{s\in L_{\mathbf m} : s<t\}]$
(so $\mathbb{P}_t'=\mathbb{P}_{\mathbf m}^{cr}[\{s\in L_{\mathbf m} : s<t\}]$
for $t\in L_{\mathbf m}$ and $\mathbb{P}_{\infty}'=\mathbb{P}_{\mathbf m}^{cr}$).
Similarly, let $\mathbb{P}_t'':=\mathbb{P}_{\mathbf n}[\{s\in L_{\mathbf m} : s\leq t\}]$.

\textbf{Main conclusion 2.26: }Let $\mathbf q$ be a $(\lambda,D)$-iteration
template. The sequence of forcing notions $(\mathbb{P}_t' : t\in L_{\mathbf q} \cup \{\infty\})$
from 2.25 has the following properties:

A) $(\mathbb{P}_t' : t\in L_{\mathbf q} \cup \{\infty\})$ is $\lessdot$-increasing,
and $s<t \in L_{\mathbf q}^+ \rightarrow \mathbb{P}_s' \lessdot \mathbb{P}_s'' \lessdot \mathbb{P}_t'$.

B) $\underset{\sim}{\eta_{t}}$ is a $\mathbb{P}_{t}''$-name of a
function from $I_{\mathbf{p}_{t}}^{1}$ to $\lambda$.

C) $(\underset{\sim}{\eta_{s}} : s<t)$ is generic for $\mathbb{P}_{t}'$.

D) $\mathbb{P}_{t}'$ is $(<\lambda)$-strategically complete and satisfies
$(\lambda,D)$-cc.

E) If $t\in L_{\mathbf q} \cup \{\infty\}$ and every set of $\leq \lambda$
elements below $t$ has a common upper bound $s<t$, then $\mathbb{P}_{t}'=\underset{s<t}{\cup}\mathbb{P}_{s}'$.

F) $|\mathbb{P}_{\infty}'|\leq (\underset{t\in L_{\mathbf q}}{\Sigma}(|I_{t}^{1}|+\lambda))^{\lambda}$.

G) If $U_1,U_2 \subseteq L_{\mathbf q}$ and $\mathbf n \restriction U_1$
is isomorphic to $\mathbf n \restriction U_2$, then $\mathbb{P}_{\mathbf m}^{cr}[U_1]=\mathbb{P}_{\mathbf n}[U_1]$
is isomorphic to $\mathbb{P}_{\mathbf m}^{cr}[U_2]=\mathbb{P}_{\mathbf n}[U_2]$.
Moreover, if $U\subseteq L_{\mathbf q}$ is closed under weak memory
(as is always the case), then $\mathbb{P}_{\mathbf m \restriction U}^{cr}$
is isomorphic to $\mathbb{P}_{\mathbf m}^{cr}[U]$. It follows that
for every $t\in L_{\mathbf q}$, $\mathbb{P}_{\mathbf m \restriction L_{<t}}^{cr}$
is isomorphic to $\mathbb{P}_{\mathbf m}^{cr}[L_{<t}]=\mathbb{P}_{t}'$.
\\
\\
H) For each $t\in L_{\bold q}$, let $V^t:=V[...,\underset{\sim}{\eta_s},...]_{s\in u_{\bold q, t}^0}$, then $\underset{\sim}{\eta_t}$ is "somewhat generic" for $\underset{\sim}{\mathbb{Q}_t^{V^t}}$ in the following sense: If $I$ is an antichain in $\underset{\sim}{\mathbb{Q}_t^{V^t}}$ that remains maximal in $V^{\mathbb{P}_{\bold n}}$ for every $\bold n$ such that $\bold m \leq \bold n \in \bold{M}_{ec}$, then $\underset{\sim}{\eta_t}$ satisfies some $p\in I$.
\\
$[$This means that if $I=\{p_{\epsilon} : \epsilon<\epsilon(*)\}$ where each $p_{\epsilon}$ has the form $(tr(p_{\epsilon}), \bold{B}_{p_{\epsilon}}(...,\underset{\sim}{\eta_{t_{\zeta}}}(a_{\zeta}),...)_{\zeta<\xi})$, then $\Vdash_{\mathbb{P}_{\bold m}^{cr}} "$There is some $\epsilon<\epsilon(*)$ such that $\underset{\sim}{\eta_t}$ extends $tr(p_{\epsilon})$ and belongs to $\bold{B}_{p_{\epsilon}}(...,\underset{\sim}{\eta_{t_{\zeta}}}(a_{\zeta}),...)_{\zeta<\xi}$".$]$
\\
$[$The reason for the absoluteness requirement is that in Requirement 1.16 we didn't demand the property of being a maximal antichain to be absolute (this would seriously restrict the range of forcing notions covered).$]$

\textbf{Proof: }A) By $2.12(F)$.

B) By the definition of $\underset{\sim}{\eta_{\alpha}}$.

C) By the definition of $\mathbb{P}_{\mathbf n}[\{i : i<\alpha\}]$.
More generally, this is true by the definition of the $\mathbb{L}_{\lambda^+}$-closure,
as $(\underset{\sim}{\eta_{\alpha}} : \alpha \in L)$ is generic for
$\mathbb{P}_{\mathbf n}[L]$ for every $L\subseteq \delta_*$.

D) By $2.12(D)$.

E) By $2.12(F)$, $\underset{s<t}{\cup}\mathbb{P}_{s}' \subseteq \mathbb{P}_{t}'$.
In the other direction, suppose that $\psi \in \mathbb{P}_{t}'=\mathbb{P}_{\mathbf n}[\{s : s<t\}]$
and let $\{p_{s(i),a(i),j(i)} : i<\lambda\}\subseteq X_{L_{<t}}$
be the set that $\mathbb{L}_{\lambda^+}$-generates $\psi$. By our
assumption, the set $\{s(i) : i<\lambda\}$ has a common upper bound
$s'<t$. Hence $\{p_{s(i),a(i),j(i)} : i<\lambda\}\subseteq X_{L_{<s'}}$,
so $\psi \in \mathbb{P}_{\mathbf n}[\{s : s<s'\}]=\mathbb{P}_{s'}'$
and equality follows.

F) As $\mathbb{P}_{\infty}'=\mathbb{P}_{\mathbf n}[L_{\mathbf q}]=\mathbb{L}_{\lambda^+}(X_{L_{\mathbf q}},\mathbb{P}_{\mathbf n})$
(recall definition $2.11$), the claim follows by the definition of
$X_{L_{\mathbf q}}$ and the definition of the $\mathbb{L}_{\lambda^+}$-closure.

G) Choose $\mathbf n \geq \mathbf m$ such that $\mathbf n \in \mathbf{M}_{ec}$
and $M_{\mathbf n}=L_{\mathbf q}$, therefore, by claim 3.12 in the next
section (the proof of which does not rely on the current claim), $\mathbb{P}_{\mathbf n}[U_1]$ is isomorphic to $\mathbb{P}_{\mathbf n}[U_2]$
where $(\mathbf n, \mathbf n, U_1, U_2)$ here stands for $(\mathbf{m}_1,\mathbf{m}_2,M_1,M_2)$
there. For the second part of the claim, choose $\mathbf m \restriction U \leq \mathbf{n}' \in \mathbf{M}_{ec}$,
then $\mathbf{n}'\restriction U=\mathbf m \restriction U=\mathbf n \restriction U$,
and as before, $\mathbb{P}_{\mathbf m}^{cr}[U]=\mathbb{P}_{\mathbf n}[U]$
is isomorphic to $\mathbb{P}_{\mathbf{n}'}[U]=\mathbb{P}_{\mathbf{m} \restriction U}^{cr}$.
\\
\\
H) Follows from the definition and the absoluteness requirement.

\textbf{\large{}3. Proving the main claim}{\large\par}

\textbf{\large{}Existence of an existentially closed extension of
adequate cardinality for a given $\mathbf m \in \mathbf M$}{\large\par}

Our goal will be to show that for every $\mathbf m \in \mathbf M$, if
$L_{\mathbf m}=M_{\mathbf m}$ and $\mathbf n = \mathbf m \restriction M$
where $M\subseteq M_{\mathbf m}$, then $\mathbb{P}_{\mathbf n}^{cr} \lessdot \mathbb{P}_{\mathbf m}^{cr}$.
In particular, in Conclusion 3.13 we get that for every $U\subseteq \delta_*$
closed under weak memory, $\mathbb{P}_{\mathbf m \restriction U}^{cr} \lessdot \mathbb{P}_{\mathbf m}^{cr}=\mathbb{P}_{\delta_*}$.
\\
\\
Remark: Note that we don't rely in this section on 2.26.

\textbf{Definition 3.1: }A) $\mathbf m \in \mathbf M$ is wide if for
every $t\in L_{\mathbf m}\setminus M_{\mathbf m}$ there are $t_{\alpha} \in L_{\mathbf m} \setminus M_{\mathbf m}$
$(\alpha< \lambda^+)$ such that:

1. $\mathbf{m} \restriction (t_{\alpha}/E_{\mathbf m}')$ is isomorphic
to $\mathbf{m} \restriction (t/E_{\mathbf m}')$ over $M_{\mathbf m}$.

2. $t_{\alpha}/E_{\mathbf m}'' \neq t_{\beta}/E_{\mathbf m}''$ for every
$\alpha<\beta<\lambda^+$.

B) $\mathbf m \in \mathbf M$ is very wide if $\mathbf m$ satisfies the
above requirements with $\lambda^+$ replaced by $|L_{\mathbf m}|$.

C) $\mathbf m \in \mathbf M$ is full if for every $\mathbf m \restriction M_{\mathbf m} \leq \mathbf n$
such that $E_{\mathbf n}''$ consists of one equivalence class, there
is $t\in L_{\mathbf m}\setminus M_{\mathbf m}$ such that $\mathbf n$ is
isomorphic to $\mathbf m \restriction (t/E_{\mathbf m}')$ over $M_{\mathbf m}$.

Remark: In the proof of theorem $2.19$, we constructeed $\mathbf n \in \mathbf{M}_{ec}$
by amalgamating $(\mathbf{n}_{\alpha}^{i} : i<\chi, \alpha<2^{\lambda_2})$.
Therefore, for every $t\in L_{\mathbf n}\setminus M_{\mathbf n}$ there
are $i$ and $\alpha$ such that $t$ belongs to $\mathbf n \restriction t/E_{\mathbf n}=\mathbf{n}_{\alpha}^{i}$.
As $\mathbf n$ includes $(\mathbf{n}_{\alpha}^{i} : i<\chi)$, by choosing
representatives $t_i \in L_{\mathbf{n}_{\alpha}^{i}}\setminus M_{\mathbf n}$
$(i<\chi)$ we get that $\mathbf{n} \restriction (t/E_{\mathbf n}')$ is
isomorphic to $\mathbf n \restriction (t_i/E_{\mathbf n}')$ for every
$i<\chi$. Since $t_i/E_{\mathbf n}'\neq t_j/E_{\mathbf n}'$ for every
$i<j<\chi$ and $|L_{\mathbf n}|=\chi$, it follows that $\mathbf n$
is very wide. By the construction of $\mathbf n$, it's also easy to
see that $\mathbf n$ is full.

\textbf{Definition 3.2: }Let $L\subseteq L_{\mathbf m}$ and $q\in \mathbb{P}_{\mathbf m}$,
we say that $p$ is the projection of $q$ to $L$ and write $p=\pi_L(q)$
if the following conditions hold:

a. $Dom(p)=Dom(q)\cap L$.

b. If $s\in Dom(p)$ then:

1. $\{\mathbf{B}_{p(s),\iota}(...,\underset{\sim}{\eta_{t_{\zeta}}}(a_{\zeta}),...)_{\zeta \in W_{p(s),\iota}} : \iota<\iota(p(s))\}=\{\mathbf{B}_{q(s),\iota}(...,\underset{\sim}{\eta_{t_{\zeta}}}(a_{\zeta}),...)_{\zeta \in W_{q(s),\iota}} : \iota< \iota(q(s)) \wedge \{t_{\zeta} : \zeta \in W_{q(s),\iota}\} \subseteq L\}$.

2. $tr(p(s))=\underset{\iota}{\cup}tr(\mathbf{B}_{q(s),\iota}(...,\underset{\sim}{\eta_{t_{\zeta}}}(a_{\zeta}),...)_{\zeta \in W_{q(s),\iota}})$
for $\iota<\iota(q(s))$ and $\{t_{\zeta} : \zeta \in W_{q(s),\iota}\}\subseteq L$.

\textbf{Observation 3.3: }Let $\mathbf m \in \mathbf M$, $L\subseteq L_{\mathbf m}$
and $q\in \mathbb{P}_{\mathbf m}$.

a. The projection $p=\pi_L(q)$ exists and $p\in \mathbb{P}_{\mathbf m}(L)$.

b. $\pi_L(q)\leq q$.

\textbf{Definition 3.4: }Let $\mathbf m \in \mathbf M$, denote by $\mathcal{F}_{\mathbf m}$
the collection of functions $f$ having the following properties:

a. There are $L_1,L_2 \subseteq L_{\mathbf m}$ such that $f$ is an
isomorphism from $\mathbf m \restriction L_1$ onto $\mathbf{m} \restriction L_2$.

b. $M_{\mathbf m}\subseteq L_1 \cap L_2$.

c. For every $t\in L_{\mathbf m} \setminus M_{\mathbf m}$, if $t\in L_l$
$(l=1,2)$ then $t/E_{\mathbf m}' \subseteq L_l$.

d. $|\{t/E_{\mathbf m}' : t\in L_l \setminus M_{\mathbf m}\}|\leq \lambda$.

e. $f$ is the identity on $M_{\mathbf m}$.

\textbf{Claim 3.5: }A. Let $\mathbf m \in \mathbf M$ be wide. For every
$f\in \mathcal{F}_{\mathbf m}$ and $X\subseteq L_{\mathbf m}$, if $|X|\leq \lambda$
then there is $g\in \mathcal{F}_{\mathbf m}$ such that:

1. $f\subseteq g$.

2. $Dom(g)=Ran(g)$.

3. $X\subseteq Dom(g)$.

B. If $g\in \mathcal{F}_{\mathbf m}$ satisfies $Dom(g)=Ran(g)$, then
$g^+:=g\cup id_{L_{\mathbf m}\setminus Dom(g)}$ is an automorphim of
$\mathbf m$.

\textbf{Proof: }A. By the proof of claim 1 in $2.19$, $f$ can be
extended to a function $f'\in \mathcal{F}_{\mathbf m}$ such that $X\subseteq Dom(f')$.
It's enough to show that for every $f'\in \mathcal{F}_{\mathbf m}$
there is $f'\subseteq g\in \mathcal{F}_{\mathbf m}$ such that $Dom(g)=Ran(g)$.
The argument is simiar to claim $1$ in $2.19$. Obviously, $Dom(f')$
and $Ran(f')$ are each a union of $M_{\mathbf m}$ with pairwise disjoint
sets of the form $t/E_{\mathbf m}''$, and for each such $t/E_{\mathbf m}''$
exactly one of the following holds:

a. $t/E_{\mathbf m}''\subseteq Dom(f')\cap Ran(f')$.

b. $t/E_{\mathbf m}'' \subseteq Dom(f')$ is disjoint to $Ran(f')$.

c. $t/E_{\mathbf m}''\subseteq Ran(f')$ is disjoint to $Dom(f')$.

As $\mathbf m$ is wide, for every $t/E_{\mathbf m}''$ as in $(b)$ there
are $\lambda^+$ $t_{\alpha}\in L_{\mathbf m}\setminus M_{\mathbf m}$
as in definition $3.1$. Therefore there is $f'\subseteq f_1 \in \mathcal{F}_{\mathbf m}$
such that $Dom(f')\subseteq Ran(f_1)$ and $Ran(f')\subseteq Dom(f_1)$.
Proceed by induction to get a sequence $f'\subseteq f_1 \subseteq...f_n \subseteq...$
of functions in $\mathcal{F}_{\mathbf m}$ such that $Dom(f_n)\subseteq Ran(f_{n+1})$
and $Ran(f_n)\subseteq Dom(f_{n+1})$ for every $n$. Obviously, $g:=\underset{n<\omega}{\cup} f_n \in \mathcal{F}_{\mathbf m}$
is as required.

B. This is easy to check. $\square$

Remark: By the last claim, given $f\in \mathcal{F}_{\mathbf m}$, we
may extend it to $g\in \mathcal{F}_{\mathbf m}$ such that $Dom(g)=Ran(g)$,
and $g$ may be extended to automorphism $h:=g^+$ of $\mathbf m$.
As in claim $3$ of $2.19$, $h$ induces an automorphism $\hat{h}$
of $\mathbb{P}_{\mathbf m}$, and obviously $\hat{f}:= \hat{h} \restriction \mathbb{P}_{\mathbf m}(Dom(f))$
is an isomorphism of $\mathbb{P}_{\mathbf m}(Dom(f))$ to $\mathbb{P}_{\mathbf m}(Ran(f))$.

\textbf{Definition 3.6: }Given $\mathbf m \in \mathbf M$, $\zeta<\lambda^+$,
$t_l \in L_{\mathbf m} \setminus M_{\mathbf m}$ $(l=1,2)$ and sequences
$\bar{s}_l$ of length $\zeta$ of elements of $t_l/E_{\mathbf m}''$,
we shall define by induction on $\gamma$ when $(t_1,\bar{s}_1)$
and $(t_2,\bar{s}_2)$ are $\gamma$-equivalent in $\mathbf m$. We
may write $\bar{s}_l$ instead of $(t_l,\bar{s}_l)$, as the choice
of $t_l$ doesn't matter as long as it's $E_{\mathbf m}''$-equivalent
to the elements of $\bar{s}_l$ (and $\bar{s}_l \neq ()$).

A. $\gamma=0:$ Let $L_l=cl(M_{\mathbf m}\cup Ran(\bar{s}_l))$ (recalling
Definition 1.9 for $l=1,2$. $(t_1,\bar{s}_1)$ is $0-$equivalent
to $(t_2,\bar{s}_2)$ if there is a function $h: L_1 \rightarrow L_2$
such that the following hold:

1. $h$ is an isomorhism from $\mathbf{m} \restriction L_1$ to $\mathbf{m} \restriction L_2$.

2. $h$ maps $\bar{s}_1$ onto $\bar{s}_2$.

3. $h$ is the identity on $M_{\mathbf m}$.

4. $h$ induces an isomorphism from $\mathbb{P}_{\mathbf m}(L_1)$ to
$\mathbb{P}_{\mathbf m}(L_2)$.

B. $\gamma$ is a limit ordinal: $\bar{s}_1$ is $\gamma$-equivalent
to $\bar{s}_2$ iff they're $\beta$-equivalent for every $\beta<\gamma$.

C. $\gamma=\beta+1$: $\bar{s}_1$ is $\gamma$-equivalent to $\bar{s}_2$
if for every $\epsilon<\lambda^+$, $l\in \{1,2\}$ and a sequence
$\bar{s}_l'$ of length $\epsilon$ of elements of $t_l/E_{\mathbf m}''$,
there exists a sequence $\bar{s}_{3-l}'$ of length $\epsilon$ of
elements of $t_{3-l}/E_{\mathbf m}''$ such that $\bar{s}_1 \hat{} \bar{s}_1'$
and $\bar{s}_2 \hat{} \bar{s}_2'$ are $\beta$-equivalent.

\textbf{Definition 3.7: }Let $\beta$ be a limit ordinal, $\mathcal{F}_{\mathbf m,\beta}$
is the collection of functions $f$ such that there is a sequence
$(t_i^l, \bar{s}_i^l : 1\leq l\leq 2, i<i(*))$ satisfying the following
conditions:

A. $i(*)<\lambda^+$.

B. For $l=1,2$, $(t_i^l : i<i(*))$ is a sequence of elements of
$L_{\mathbf m} \setminus M_{\mathbf m}$ such that for every $i<j<i(*)$,
$t_i^l$ and $t_j^l$ are not $E_{\mathbf m}''$-equivalent.

C. $\bar{s}_i^l$ is a sequence of length $\zeta(i)<\lambda^+$ of
elements of $t_i^l/E_{\mathbf m}''$.

D. $\bar{s}_i^1$ and $\bar{s}_i^2$ are $\beta$-equivalent.

E. $f$ is an isomorphism from $\mathbf m \restriction L_1$ to $\mathbf m \restriction L_2$
where $L_l=\underset{i<i(*)}{\cup}Ran(\bar{s}_i^l) \cup M_{\mathbf m}$
$(l=1,2)$.

F. For every $i<i(*)$, $f$ maps $\bar{s}_i^1$ onto $\bar{s}_i^2$.

G. $f$ is the identity on $M_{\mathbf m}$.

\textbf{Claim 3.8: }Let $\mathbf m \in \mathbf M$ be wide and suppose
that:

A. $\mathbf{m}_1 \leq \mathbf m$.

B. For every $t\in L_{\mathbf m} \setminus L_{\mathbf{m}_1}$, $\zeta<\lambda^+$
and a sequence $\bar{s}$ of length $\zeta$ of elements of $t/E_{\mathbf m}''$,
there is a sequence $(t_i, \bar{s}_i : i<\lambda^+)$ such that:

1. $t_i \in L_{\mathbf{m}_1}\setminus M_{\mathbf{m}_1}$.

2. If $i<j<\lambda^+$ then $t_i/E_{\mathbf m}' \neq t_j/E_{\mathbf{m}_1}'$.

3. $\bar{s}_i$ is a sequence of length $\zeta$ of elements of $t_i/E_{\mathbf{m}_1}''$.

4. $(t_i, \bar{s}_i)$ is $1-$equivalent to $(t,\bar{s})$ in $\mathbf m$.

Then $\mathbb{P}_{\mathbf{m}_1}\lessdot \mathbb{P}_{\mathbf m}$.

\textbf{Proof: }We shall freely use the results from Section 4 (of course, it should be noted that none of the relevant results in Section 4 relies on the current claim). Specifically,
we shall use the fact that a function $f\in \mathcal{F}_{\mathbf m, \beta}$
induces an isomorphism $\hat{f}$ from $\mathbb{P}_{\mathbf m}(L_1)$
to $\mathbb{P}_{\mathbf m}(L_2)$ for $L_1$ and $L_2$ as in definition
$3.7$ (see Claim 4.3). Now, note that if $f\in \mathcal{F}_{\mathbf m,\beta}$
for $0<\beta$ and $L\subseteq L_{\mathbf m}$ such that $|L|\leq \lambda$,
then by the definition of $1-$equivalence, $f$ can be extended to
a function $g\in \mathcal{F}_{\mathbf m,0}$ such that $L\subseteq Dom(g)$.
Hence $\hat{g}$ is an isomorphism with domain $\mathbb{P}_{\mathbf m}(L_1 \cup L)$
such that $\hat{f} \subseteq \hat{g}$.

\textbf{Claim 1: }If $0<\beta$ then $\hat{f}$ preserves compatibility
and incompatibility.

\textbf{Proof: }Assume that $p,q \in Dom(\hat{f})$ and $r$ is a
common upper bound in $\mathbb{P}_{\mathbf m}$. If $r\in Dom(\hat{f})$,
then since $\hat{f}$ is order preserving, then $\hat{f}(p)$ and
$\hat{f}(q)$ have a common upper bound. If $r\notin Dom(\hat{f})$,
then use the definition of $\mathcal{F}_{\mathbf m,\beta}$ to extend
$\hat{f}$ to a function $\hat{g}$ such that $\hat{g}(r)$ is defined
(and $g\in \mathcal{F}_{\mathbf m, 0}$), and repeat the previous argument.
The proof in the other direction repeats the same arguments for $f^{-1}$.

\textbf{Claim 2: }Suppose that $i(*)<\lambda^+$, $p_i \in \mathbb{P}_{\mathbf{m}_1}$
$(i<i(*))$ and $p\in \mathbb{P}_{\mathbf m}$, then there is $p^* \in \mathbb{P}_{\mathbf{m}_1}$
such that:

1. $\mathbb{P}_{\mathbf m} \models p_i \leq p$ iff $\mathbb{P}_{\mathbf m} \models p_i \leq p^*$.

2. For every $i<i(*)$, $p$ and $p_i$ are incompatible in $\mathbb{P}_{\mathbf m}$
iff $p^*$ and $p_i$ are incompatible in $\mathbb{P}_{\mathbf m}$. 

\textbf{Proof: }Note that if $p\in \mathbb{P}_{\mathbf m}$ then $p\in \mathbb{P}_{\mathbf{m}_1}$
iff $fsupp(p)\subseteq L_{\mathbf{m}_1}$, therefore we need to find
$p^* \in \mathbb{P}_{\mathbf m}$ satisfying the requirements of the
claim such that $fsupp(p^*)\subseteq L_{\mathbf{m}_1}$. Let $L_1\subseteq L_{\mathbf{m}_1}$
be a set containing $(\underset{i<i(*)}{\cup} fsupp(p_i))\cup M_{\mathbf m}$
and closed under weak memory, such that $|L_1\setminus M_{\mathbf m}|\leq \lambda$
(such $L_1$ exists, recalling that $i(*)<\lambda^+$ and $|w_t^0|\leq \lambda$),
then $\{p_i : i<i(*)\} \subseteq \mathbb{P}_{\mathbf m}(L_1)$. For
every $p_i$ that is compatible with $p$ in $\mathbb{P}_{\mathbf m}$,
let $q_i$ be a common upper bound. As before, there is $L_2 \subseteq L_{\mathbf m}$
containing $L_1 \cup (\cup fsupp(q_i)) \cup fsupp(p)$ and closed
under weak memory such that $|L_2 \setminus M_{\mathbf m}|\leq \lambda$
and $\mathbb{P}_{\mathbf m}(L_2)$ contains $p$ and all of the $q_i$.
We shall prove that it's enough to show that there is $f\in \mathcal{F}_{\mathbf m,1}$
such that $L_2\subseteq Dom(f)$, $Ran(f)\subseteq L_{\mathbf{m}_1}$
and $f$ is the identity on $L_1$. For such $f$ define $p^*:= \hat{f}(p)$.
Now $\hat{f}$ is the identity on $\{p_i : i<i(*)\}$ and $\hat{f}(p)\in \mathbb{P}_{\mathbf{m}_1}$.
By a previous claim, $\hat{f}$ preserves order and incompatibility,
hence $p^*$ is as required. It remains to find $f$ as above. WLOG
$L_2 \cap L_{\mathbf{m}_1}\subseteq L_1$. Let $(t_j : j<j(*))$ be
a sequence of representatives of pairwise $E_{\mathbf m}''$-inequivalent
members of $L_{\mathbf m} \setminus M_{\mathbf m}$ such that every $t\in L_2 \setminus L_1$
is $E_{\mathbf m}''$-equivalent to some $t_j$. For every such $t_j$,
let $\bar{s}_j$ be the sequence of members of $t_j/E_{\mathbf m}''$
in $L_2 \setminus L_1$. By the assumption, for every pair $(\bar{s}_j,t_j)$
as above there exist $\lambda^+$ pairs $((\bar{s}_{j,i}, t_{j,i}) : i<\lambda^+)$
which are $1-$equivalent as in the assumption of the above claim.
By induction on $j<j(*)<\lambda^+$ choose the pair $(\bar{s}_{j,i(j)},t_{j,i(j)})$
such that $t_{j,i(j)}/E_{\mathbf{m}_1}''$ are with no repetitions (this
is possible as $j(*)<\lambda^+)$. Now define $f\in \mathcal{F}_{\mathbf m, 1}$
as the function extending $id \restriction L_1$ witnessing the equivalence
of the pairs we chose. Obviously, $f$ is as required.

\textbf{Claim 3:} $\mathbb{P}_{\mathbf{m}_1}\lessdot \mathbb{P}_{\mathbf m}$.

Remark: We shall use Section 4 in the following proof.

\textbf{Proof: }We shall prove by induction on $\gamma$ that $\mathbb{P}_{\mathbf{m}_1}(L_{\mathbf{m}_1,\gamma}^{dp}) \lessdot \mathbb{P}_{\mathbf m}(L_{\mathbf m,\gamma}^{dp})$.
For $\gamma$ large enough we'll get $\mathbb{P}_{\mathbf{m}_1}\lessdot \mathbb{P}_{\mathbf m}$.

\textbf{First case: $\gamma=0$.}

Denote $E=E_{\mathbf m}''\restriction L_{\mathbf m,\gamma}^{dp}$. $E$
is an equivalence relation and $E\restriction L_{\mathbf{m}_1,\gamma}^{dp}=E_{\mathbf{m}_1}''\restriction L_{\mathbf{m}_1,\gamma}^{dp}$.
Now the claim follows by the fact that $\mathbb{P}_{\mathbf m}(L_{\mathbf{m},\gamma}^{dp})$
(and similarly $\mathbb{P}_{\mathbf{m}_1}(L_{\mathbf{m}_1,\gamma}^{dp})$)
can be represented as a product with $<\lambda$ support of $\{\mathbb{P}_{\mathbf m}(t/E) : t\in L_{\mathbf{m},\gamma}^{dp}\}$.

\textbf{Second case: $\gamma=\beta+1$.}

Denote $M_{\beta}:= \{t\in M_{\mathbf m} : dp_{\mathbf m}^{*}(t)=\beta\}$,
then $M_{\beta}$'s members are pairwise incomparable.

\textbf{Claim: }$\mathbb{P}_{\mathbf{m}_1}(L_{\mathbf{m}_1,\beta}^{dp} \cup M_{\beta}) \lessdot \mathbb{P}_{\mathbf m}(L_{\mathbf m,\beta}^{dp} \cup M_{\beta})$.

\textbf{Proof: }We shall prove the claim by a series of subclaims.

\textbf{Subclaim: }Given $p,q\in \mathbb{P}_{\mathbf{m}_1}(L_{\mathbf{m}_1,\beta}^{dp} \cup M_{\beta})$,
$\mathbb{P}_{\mathbf{m}_1}(L_{\mathbf{m}_1,\beta}^{dp} \cup M_{\beta}) \models p\leq q$
if and only if $\mathbb{P}_{\mathbf m}(L_{\mathbf m,\beta}^{dp} \cup M_{\beta}) \models p\leq q$.

\textbf{Proof: }Note that $L_{\mathbf{m}_1,\beta}^{dp} \cup M_{\beta}$
and $L_{\mathbf m,\beta}^{dp} \cup M_{\beta}$ are initial segments
of $L_{\mathbf{m}_1}$ and $L_{\mathbf m}$, respectively. Note also that
if $\mathbf n \in \mathbf M$ and $L_1\subseteq L_2\subseteq L_{\mathbf n}$,
then $\mathbb{P}_{\mathbf{n} \restriction L_1}\lessdot \mathbb{P}_{\mathbf{n}\restriction L_2}$,
and if $L\subseteq L_{\mathbf n}$ is an initial segment then $\mathbb{P}_{\mathbf n}(L)=\mathbb{P}_{\mathbf n \restriction L}$.
Obviously, $L_{\mathbf{m}_1,\beta}^{dp}$ and $L_{\mathbf m,\beta}^{dp}$
are initial segments of $L_{\mathbf{m}_1}$ and $L_{\mathbf m}$, respectively.
Now the claim follows by the definition of the forcing's partial order
(definition 1.8) and the induction hypothesis.

\textbf{Subclaim: }Given $p_1,p_2 \in \mathbb{P}_{\mathbf{m}_1}(L_{\mathbf{m}_1,\beta}^{dp} \cup M_{\beta})$,
$p_1$ and $p_2$ are compatible in $\mathbb{P}_{\mathbf{m}_1}(L_{\mathbf{m}_1,\beta}^{dp} \cup M_{\beta}$
if and only if theey're compatible in $\mathbb{P}_{\mathbf m}(L_{\mathbf m,\beta}^{dp} \cup M_{\beta})$.

\textbf{Proof: }By the previous subclaim, if $p_1$ and $p_2$ are
compatible in $\mathbb{P}_{\mathbf{m}_1}(L_{\mathbf{m}_1,\beta}^{dp} \cup M_{\beta})$
then they're compatible in $\mathbb{P}_{\mathbf m}(L_{\mathbf m, \beta}^{dp} \cup M_{\beta})$.
Let us now prove the other direction. Suppose that $p\in \mathbb{P}_{\mathbf m}(L_{\mathbf m,\beta}^{dp} \cup M_{\beta})$
is a common upper bound of $p_1$ and $p_2$ in $\mathbb{P}_{\mathbf m}(L_{\mathbf m,\beta}^{dp} \cup M_{\beta})$.
As in the proof of claim 2 above, find $f\in \mathcal{F}_{\mathbf m,1}$
such that $fsupp(p)\cup fsupp(p_1)\cup fsupp(p_2)\subseteq Dom(f)$,
$f\restriction (fsupp(p_1)\cup fsupp(p_2)\cup M_{\beta})$ is the
identity and $Ran(f)\subseteq L_{\mathbf{m}_1}$. Note that if $t\in Dom(f) \cap L_{\mathbf m,\beta}^{dp}$
then $f(t) \in L_{\mathbf{m}_1,\beta}^{dp}$. Since $f((Dom(f)\cap L_{\mathbf m,\beta}^{dp})\cup M_{\beta})\subseteq L_{\mathbf{m}_1,\beta}^{dp} \cup M_{\beta}$,
it follows that $\hat{f}(p) \in \mathbb{P}_{\mathbf{m}_1}(L_{\mathbf{m}_1,\beta}^{dp} \cup M_{\beta})$,
and as before, it's a common upper bound as required.

\textbf{Claim: }$\mathbb{P}_{\mathbf{m}_1}(L_{\mathbf{m}_1,\beta}^{dp} \cup M_{\beta}) \lessdot \mathbb{P}_{\mathbf m}(L_{\mathbf{m},\beta}^{dp} \cup M_{\beta})$.

\textbf{Proof: }Let $I\subseteq \mathbb{P}_{\mathbf{m}_1}(L_{\mathbf{m}_1,\beta}^{dp} \cup M_{\beta})$
be a maximal antichain and suppose towards contradiction that $p\in \mathbb{P}_{\mathbf m}(L_{\mathbf m,\beta}^{dp} \cup M_{\beta})$
contradicts in $\mathbb{P}_{\mathbf m}(L_{\mathbf m,\beta}^{dp} \cup M_{\beta})$
all elements of $I$. As before, choose $f\in \mathcal{F}_{\mathbf m,1}$
which is the identity on $M_{\beta}$ and on $fsupp(q)$ for every
$q\in I$, such that $Ran(f)\subseteq L_{\mathbf{m}_1}$ (hence $f(Dom(f) \cap L_{\mathbf m,\beta}^{dp})\subseteq L_{\mathbf{m}_1,\beta}^{dp}$).
Now $\hat{f}(p) \in \mathbb{P}_{\mathbf{m}_1}(L_{\mathbf{m}_1,\beta}^{dp} \cup M_{\beta})$
and $\hat{f}$ is order preserving, hence $\hat{f}(p)$ contradicts
all members of $I$ in $\mathbb{P}_{\mathbf{m}_1}(L_{\mathbf{m}_1,\beta}^{dp} \cup M_{\beta})$,
contradicting our assumption. Therefore $I$ is a maximal antichain
in $\mathbb{P}_{\mathbf m}(L_{\mathbf m,\beta}^{dp} \cup M_{\beta})$
and $\mathbb{P}_{\mathbf{m}_1}(L_{\mathbf{m}_1,\beta}^{dp}) \lessdot \mathbb{P}_{\mathbf m}(L_{\mathbf m,\beta}^{dp} \cup M_{\beta})$.

We shall now continue with the proof of the induction.

Denote $L_*= L_{\mathbf m,\gamma}^{dp} \setminus (L_{\mathbf m,\beta}^{dp} \cup M_{\beta})$
and denote by $\mathcal{E}$ the collection of pairs $(s_1,s_2)$
such that $s_1,s_2 \in L_{\mathbf m,\gamma}^{dp} \setminus (L_{\mathbf m,\beta}^{dp} \cup M_{\beta})$
and $s_1/E_{\mathbf m}''=s_2/E_{\mathbf m}''$, so $\mathcal{E}$ is an
equivalence relation. Note also that if $s_1$ and $s_2$ are not
$\mathcal{E}$-equivalent, then they're incomparable. Now observe that
the following are true:

1. Suppose that $s\in L_*$, $t\in L_{\mathbf m}$ and $t<s$. If $t\notin L_{\mathbf m,\beta}^{dp}$,
then there is $r\in M_{\beta}$ such that $r\leq t$. Therefore, either
$t\in M_{\beta}$ or $t\in L_*$ and $t\mathcal{E}s$, hence $L_{\mathbf m,<s} \subseteq L_{\mathbf{m},\beta}^{dp} \cup M_{\beta} \cup (s/ \mathcal{E})$.

2. Similarly, if $s\in L_* \cap L_{\mathbf{m}_1}$, then $L_{\mathbf{m}_1,<s} \subseteq L_{\mathbf{m}_1,\beta}^{dp} \cup M_{\beta} \cup (s/\mathcal{E})$.

Let $\{X_{\epsilon} : \epsilon<\epsilon(*)\}$ be the collection of
$\mathcal{E}$-equivalence classes and let $U_1=\{\epsilon : X_{\epsilon} \subseteq L_{\mathbf{m}_1,\gamma}^{dp}\}$,
$Z=L_{\mathbf m,\beta}^{dp} \cup \{X_{\epsilon} : \epsilon \notin U_1\} \cup M_{\beta}$,
$Y=L_{\mathbf m,\beta}^{dp} \cup \{X_{\epsilon : \epsilon \in U_1}\} \cup M_{\beta}$.

It's easy to see that:

1. $L_{\mathbf{m}_1,\gamma}^{dp}=\cup \{X_{\epsilon} : \epsilon \in U_1\} \cup L_{\mathbf{m}_1,\beta}^{dp} \cup M_{\beta}$.

2. $Z\cap L_{\mathbf{m}_1,\gamma}^{dp}=L_{\mathbf{m}_1,\beta}^{dp} \cup M_{\beta}$.

3. $Z\cup L_{\mathbf{m}_1,\gamma}^{dp}=L_{\mathbf m,\gamma}^{dp} \cup M_{\beta}$.

4. $Z\cap Y=L_{\mathbf m,\beta}^{dp} \cup M_{\beta}$.

5. $Z\cup Y=L_{\mathbf m,\gamma}^{dp}$.

By observation $(1)$ (the first one), $Y$ and $Z$ are initial segments
of $L_{\mathbf m}$, and if $s\in Z\setminus Y$ and $t\in Y\setminus Z$,
then $t$ and $s$ are incomparable. Note also that $\mathbb{P}_{\mathbf m}(Y\cup Z)=\mathbb{P}_{\mathbf m}(L_{\mathbf m,\gamma}^{dp})$.
Since $Y$ is an initial segment, $\mathbb{P}_{\mathbf m}(Y) \lessdot \mathbb{P}_{\mathbf m}(Y\cup Z)$.
Let $Y_1=L_{\mathbf{m}_1,\gamma}^{dp} \cup M_{\beta}$, $Y_2=L_{\mathbf m,\beta}^{dp} \cup M_{\beta}$,
obviously $Y_2$ and $Y_1 \cup Y_2$ are initial segments of $L_{\mathbf m}$.
Let $Y_0=Y_1 \cap Y_2$, then $\mathbb{P}_{\mathbf{m}_1}(Y_0)=\mathbb{P}_{\mathbf{m}_1}(L_{\mathbf{m}_1,\beta}^{dp} \cup M_{\beta}) \lessdot \mathbb{P}_{\mathbf m}(L_{\mathbf m,\beta}^{dp} \cup M_{\beta})=\mathbb{P}_{\mathbf m}(Y_2)$.
Since $\mathbb{P}_{\mathbf{m}_1}(Y_0)=\mathbb{P}_{\mathbf m}(Y_0)$, we
get $\mathbb{P}_{\bold m}(Y_0) \lessdot \mathbb{P}_{\mathbf m}(Y_2)$.
Note also that $Y_1 \setminus Y_0$ is disjoint to $M_{\mathbf m}$,
$Y_0$ is an initial segment of $Y_1$ and if $t\in Y_1 \setminus M_{\mathbf m}$
then $(t/E_{\mathbf m}'') \cap L_{\mathbf m,<s} \subseteq Y_1$.

Finally, the desired conclusion will be derived from the following
two claims:

\textbf{Claim 3 (1) }Suppose that $Y_1,Y_2, Y_3 \subseteq L_{\mathbf m}$
and $Y_0=Y_1 \cap Y_2$, then $\mathbb{P}_{\mathbf m}(Y_1) \lessdot \mathbb{P}_{\mathbf m}(Y_3)$
if the following conditions hold:

1. $Y_2 \subseteq Y_3$ are initial segments of $L_{\mathbf m}$.

2. $Y_1 \subseteq Y_2$ and $Y_0$ is an initial segment of $Y_1$.

3. $\mathbb{P}_{\mathbf m}(Y_0) \lessdot \mathbb{P}_{\mathbf m}(Y_2)$.

4. $Y_1 \setminus Y_0 \cap M_{\mathbf m}=\emptyset$.

5. If $t\in Y_1 \setminus M_{\mathbf m}$ then $t/E_{\mathbf m}'' \cap L_{\mathbf m,<t} \subseteq Y_1$.

\textbf{Claim 3 (2): }$\mathbb{P}_{\mathbf{m}_1}(L_1)=\mathbb{P}_{\mathbf{m}_2}(L_1) \lessdot \mathbb{P}_{\mathbf{m}_2}$
if the following conditions hold:

1. $\mathbf{m}_1 \leq \mathbf{m}_2$.

2. $L_0 \subseteq L_1 \subseteq L_{\mathbf{m}_1}$.

3. $L_0$ is an initial segment of $L_1$.

4. $\mathbb{P}_{\mathbf{m}_1}(L_0)=\mathbb{P}_{\mathbf{m}_2}(L_0)$.

5. $\mathbb{P}_{\mathbf{m}_l}(L_0) \lessdot \mathbb{P}_{\mathbf{m}_l}$
for $l=1,2$.

6. if $t\in L_1 \setminus L_0$ then $t\notin M_{\mathbf{m}_2}$ and
$L_{\mathbf{m}_1,<t} \cap (t/E_{\mathbf{m}_1})=L_{\mathbf{m}_2,<t} \cap (t/E_{\mathbf m}) \subseteq L_1$.

By claim $3 (2)$, with $(\mathbf{m}_1, \mathbf m, Y_0, Y_1)$ standing
for $(\mathbf{m}_1,\mathbf{m}_2, L_0,L_1)$ in the claim, we get $\mathbb{P}_{\mathbf{m}_1}(Y_1)=\mathbb{P}_{\mathbf m}(Y_1) \lessdot \mathbb{P}_{\mathbf m}$.
By claim $3 (1)$, it follows that $\mathbb{P}_{\mathbf m}(L_{\mathbf{m}_1,\gamma}^{dp})=\mathbb{P}_{\mathbf m}(Y_1) \lessdot \mathbb{P}_{\mathbf m}(Y_1 \cup Y_2)=\mathbb{P}_{\mathbf m}(Y) \lessdot \mathbb{P}_{\mathbf m}(Y \cup Z)=\mathbb{P}_{\mathbf m}(L_{\mathbf m, \gamma}^{dp})$.
Together we get $\mathbb{P}_{\mathbf{m}_1}(L_{\mathbf{m}_1,\gamma}^{dp})=\mathbb{P}_{\mathbf{m}_1}(Y_1)=\mathbb{P}_{\mathbf m}(Y_1) \lessdot \mathbb{P}_{\mathbf m}(L_{\mathbf m,\gamma}^{dp})$.

\textbf{Proof of claim 3 (1): }We shall prove by induction on $\gamma$
that if $(Y_0, Y_1, Y_2, Y_3)$ are as in the claim's assumptions
and $dp(Y_1)\leq \gamma$ then:

1. $\mathbb{P}_{\mathbf m}(Y_1)\lessdot \mathbb{P}_{\mathbf m}(Y_3)$.

2. If A) then B) where:

A) 1. $p_3 \in \mathbb{P}_{\mathbf m}(Y_3)$.

2. $p_0 \in \mathbb{P}_{\mathbf m}(Y_0)$.

3. If $p_0 \leq q_0 \in \mathbb{P}_{\mathbf m}(Y_0)$ then $p_2=p_3 \restriction Y_2$
and $q_0$ are compatible.

4. $p_1=p_0 \cup (p_3 \restriction (Y_1 \setminus Y_0))$.

B) If $p_1 \leq q_1 \in \mathbb{P}_{\mathbf m}(Y_1)$ then $q_1$ and
$p_3$ are compatible in $\mathbb{P}_{\mathbf m}(Y_3)$.

Suppose we arrived at stage $\gamma$:

For part 2 of the induction claim: By assumption $5$ and the definition
of the conditions in the iteration, $fsupp(p_3 \restriction (Y_1\ \setminus Y_0)) \subseteq Y_1$,
hence $p_1 \in \mathbb{P}_{\mathbf m}(Y_1)$. Suppose towards contradiction
that A) does not hold for some $p_1 \leq q_1 \in \mathbb{P}_{\mathbf m}(Y_1)$,
then there are $s\in Dom(q_1) \cap Dom(p_3)$ and $p_3^+ \in \mathbb{P}_{\mathbf m}(L_{\mathbf m,<s})$
such that $p_3 \restriction L_{\mathbf{m},<s}, q_1 \restriction L_{\mathbf m,<s} \leq p_3^+$
and $p_3^+ \restriction L_{\mathbf{m},<s} \Vdash "q_1(s)$ and $p_3(s)$
are incompatible''. Since $s\in Dom(q_1) \subseteq Y_1$ and $Y_2$
is an initial segment, then necessarily $s\notin Y_0$ (otherwise
we get a contradiction to assumption A)(3)). $\mathbb{P}_{\mathbf m} \models p_1 \leq q_1$,
hence $q_1 \restriction L_{\mathbf m,<s} \Vdash p_1(s) \leq q_1(s)$.
As $q_1 \restriction L_{\mathbf{m},<s} \leq p_3^+$, it follows that
$p_4^+ \restriction L_{\mathbf m,<s} \Vdash p_1(s) \leq q_1(s)$. Now
$s\in Y_1 \setminus Y_0$, hence $p_1(s)=p_3(s)$, hence $p_3^+ \restriction L_{\mathbf m,<s} \Vdash p_3(s) \leq q_1(s)$,
contradicting the choice of $p_3^+$. This proves part 2.

For part 1 of the induction claim: Obviously, $\mathbb{P}_{\mathbf m}(Y_1) \subseteq \mathbb{P}_{\mathbf m}(Y_3)$
and $\mathbb{P}_{\mathbf m}(Y_1) \models p\leq q$ iff $\mathbb{P}_{\mathbf m}(Y_3) \models p\leq q$.
Suppose now that $q_1,q_2 \in \mathbb{P}_{\mathbf m}(Y_1)$ and $p_3 \in \mathbb{P}_{\mathbf m}(Y_3)$
is a common upper bound, we shall prove the existence of a common
upper bound in $\mathbb{P}_{\mathbf m}(Y_1)$. Since $Y_2$ is an initial
segment, it follows that $fsupp(p_3 \restriction Y_2) \subseteq Y_2$,
hence $p_3 \restriction Y_2 \in \mathbb{P}_{\mathbf m}(Y_2)$. Since
$\mathbb{P}_{\mathbf m}(Y_0) \lessdot \mathbb{P}_{\mathbf m}(Y_2)$, it
follows that there exists $p_0 \in \mathbb{P}_{\mathbf m}(Y_0)$ such
that if $p_0 \leq q \in \mathbb{}_{\mathbf m}(Y_0)$, then $q$ and
$p_3 \restriction Y_2$ are compatible. Let $p_1:=p_0 \cup (p_3 \restriction Y_1 \setminus Y_0)$.
As in the proof of part (2), $p_1 \in \mathbb{P}_{\mathbf m}(Y_1)$.
If $p_1 \leq p_1' \in \mathbb{P}_{\mathbf m}(Y_1)$, then by part (2)
of the induction claim, $p_1'$ is compatible with $p_3$. We shall
prove that $p_1$ is a common upper bound of $q_1$ and $q_2$. As
we may replace $p_0$ by $p_0 \leq p_0' \in \mathbb{P}_{\mathbf m}(Y_0)$,
we may assume WLOG that $Dom(q_l) \cap Y_0 \subseteq Dom(p_0) \subseteq Dom(p_1)$
$(l=1,2)$. Also $Dom(q_l) \setminus Y_0 \subseteq Dom(p_3) \setminus Y_0$.
As $Y_2$ is an initial segment, it follows from our assumptions that
$\mathbb{P}_{\mathbf m}(Y_0) \lessdot \mathbb{P}_{\mathbf m}(Y_2) \lessdot \mathbb{P}_{\mathbf m}$.
Since $p_0$ is compatible with $p_3 \restriction Y_0$ in $\mathbb{P}_{\mathbf m}$,
they're compatible in $\mathbb{P}_{\mathbf m}(Y_0)$, hence there is
a common upper bound for $p_0, q_1 \restriction Y_0$ and $q_2 \restriction Y_0$.
Therefore WLOG $q_l \restriction Y_0 \leq p_0$ $(l=1,2)$. Assume
towards contradiction that $q_l \leq p_1$ doesn't hold, then there
is $s\in Dom(q_l)$ such that $q_l \restriction L_{\mathbf m,<s} \leq p_1 \restriction L_{\mathbf m,<s}$
but $p_1 \restriction L_{\mathbf m,<s} \nVdash q_l(s) \leq p_1(s)$.
If $s\in Y_0$, then as $Y_0$ is an initial segment of $Y_1$, it
follows that $p_0 \restriction L_{\mathbf m,<s}=p_1 \restriction L_{\mathbf m,<s}$
and $p_0(s)=p_1(s)$, contradicting the fact that $q_l \leq p_0$.
Therefore $s\in Y_1 \setminus Y_0$. Let $Y_0'=Y_0$, $Y_1'=Y_0 \cup (Y_1 \cap L_{\mathbf m,<s})$,
$Y_2'=Y_2$ and $Y_3'=Y_3$, then $(Y_0',Y_1',Y_2',Y_3')$ satisfy
the assumptions of claim 3 (1) and $dp_{\mathbf m}(Y_1')=dp_{\mathbf m}(s)<\gamma$.
By the induction hypothesis, $\mathbb{P}_{\mathbf m}(Y_1') \lessdot \mathbb{P}_{\mathbf m}(Y_3')$.
As $s\in Y_1 \setminus Y_0$ (and by the assumption, $s\notin M_{\mathbf m}$),
it follows from the assumption that $(s/E_{\mathbf m}) \cap L_{\mathbf m,<s} \subseteq Y_1'$.
Therefore by the definition of the conditions in the iteration, $fsupp(p_1 \restriction \{s\}), fsupp(q_l \restriction \{s\}) \subseteq Y_1'$.
Therefore $p_1(s)$ and $q_l(s)$ are $\mathbb{P}_{\mathbf m}(Y_1')$-names.
Recall that $p_1 \restriction L_{\mathbf m,<s} \nVdash q_1(s) \leq p_1(s) $,
$L_{\mathbf m,<s} \subseteq Y_3=Y_3'$ are initial segments and $\mathbb{P}_{\mathbf m}(Y_1') \lessdot \mathbb{P}_{\mathbf m}(Y_3')$.
Therefore $\mathbb{P}_{\mathbf m}(Y_1' \cap L_{\mathbf m,<s}) \lessdot \mathbb{P}_{\mathbf m}(Y_3' \cap L_{\mathbf m,<s})$
and $fsupp(p_1 \restriction L_{\mathbf m,<s}) \subseteq Y_1 \cap L_{\mathbf m,<s}$.
Therefore $p_1 \restriction (Y_1' \cap L_{\mathbf m,<s}) \nVdash_{\mathbb{P}_{\mathbf m}(Y_1' \cap L_{\mathbf m,<s})} q_l(s) \leq p_1(s)$,
hence there exists $p_1 \restriction (Y_1' \cap L_{\mathbf m,<s}) \leq p_1^+ \in \mathbb{P}_{\mathbf m}(Y_1' \cap L_{\mathbf m,<s})$
such that $p_1^+ \Vdash_{\mathbb{P}_{\mathbf m}(Y_1' \cap L_{\mathbf m,<s})} \neg q_l(s) \leq p_1(s)$,
hence $p_1^+ \Vdash_{\mathbb{P}_{\mathbf m}(Y_3' \cap L_{\mathbf m,<s})} \neg q_l(s) \leq p_1(s)$.
By part (2) of the induction hypothesis with $\gamma_1=dp_{\mathbf m}(s)$
as $\gamma$ and $(p_1 \restriction (Y_1' \cap L_{\mathbf m,<s}), p_1^+, p_3 \restriction L_{\mathbf m,<s})$
standing for $(p_1,q_1,p_3)$ there, $p_1^+$ is compatible with $p_3 \restriction L_{\mathbf m,<s}$
in $\mathbb{P}_{\mathbf m}(L_{\mathbf m,<s})$. Let $p_3^+$ be a common
upper bound. As $q_l \leq p_3$, $p_3^+ \Vdash_{\mathbb{P}_{\mathbf m}(Y_1' \cap L_{\mathbf m,<s})} q_l(s) \leq p_3(s)=p_1(s)$
(recalling that $s\notin Y_0$). As $p_1^+ \Vdash_{\mathbb{P}_{\mathbf m}(Y_1' \cap L_{\mathbf m,<s})} \neg q_l(s) \leq p_1(s)$,
we get $p_3^+ \Vdash_{\mathbb{P}_{\mathbf m}(Y_1' \cap L_{\mathbf m,<s})} \neg q_l(s) \leq p_1(s)$.
Together we got a contradiction, hence $p_1$ is the desired common
upper bound and $\mathbb{P}_{\mathbf m}(Y_1) \subseteq_{ic} \mathbb{P}_{\mathbf m}(Y_3)$.
In order to show that $\mathbb{P}_{\mathbf m}(Y_1) \lessdot \mathbb{P}_{\mathbf m}(Y_3)$,
note that for every $p_3 \in \mathbb{P}_{\mathbf m}(Y_3)$ we can repeat
the argument in the beginning of the proof and get $p_0 \in \mathbb{P}_{\mathbf m}(Y_0)$
and $p_1 \in \mathbb{P}_{\mathbf m}(Y_1)$ that satisfy the requirements
in part (2) of the induction. Hence, part (2) holds for $(p_0,p_1,p_3)$
hence $\mathbb{P}_{\mathbf m}(Y_1) \lessdot \mathbb{P}_{\mathbf m}(Y_3)$.

\textbf{Proof of claim 3 (2): }For $l=1,2$ define the sequence $\bar{L}_l=(L_{l,i} : i<4)$
as follows: $L_{l,0}=L_0$, $L_{l,1}=L_1$, $L_{l,3}=L_{\mathbf{m}_l}$
and $L_{l,2}$ will be defined as the set of $s\in L_{\mathbf{m}_l}$
such that $s\leq t$ for some $t\in L_0$. It's easy to see that $(\mathbf{m}_l, \bar{L}_l)$
satisfies the assumptions of claim 3 (1), therefore $\mathbb{P}_{\mathbf{m}_l}(L_1)=\mathbb{P}_{\mathbf{m}_l}(L_{l,1}) \lessdot \mathbb{P}_{\mathbf{m}_l}(L_{l,3})=\mathbb{P}_{\mathbf{m}_l}$,
so $\mathbb{P}_{\mathbf{m}_2}(L_1) \lessdot \mathbb{P}_{\mathbf{m}_2}$,
as required. We shall now prove the remaining part of the claim. Let
$(s_{\alpha} : \alpha<\alpha(*))$ be an enumeration of the elements
of $L_1 \setminus L_0$ such that if $s_{\alpha}<s_{\beta}$ then
$\alpha \leq \beta$. For every $\alpha \leq \alpha(*)$ define $L_{0,\alpha}=L_0 \cup \{s_{\beta} : \beta<\alpha\}$.
We shall prove by induction on $\alpha \leq \alpha(*)$ that $\mathbb{P}_{\mathbf{m}_1}(L_{0,\alpha})=\mathbb{P}_{\mathbf{m}_2}(L_{0,\alpha})$.
For $\alpha=\alpha(*)$ we'll have $\mathbb{P}_{\mathbf{m}_1}(L_1)=\mathbb{P}_{\mathbf{m}_2}(L_1)$
as required. 

\textbf{First case $(\alpha=0)$: }In this case $L_0=L_{0,\alpha}$
and the claim follows from assumption $(4)$.

\textbf{Second case ($\alpha$ is a limit ordinal): }Obviously $\mathbb{P}_{\mathbf{m}_1}(L_{0,\alpha})=\mathbb{P}_{\mathbf{m}_2}(L_{0,\alpha})$
as sets. By the definition of the partial order and the induction
hypothesis, it follows that $\mathbb{P}_{\mathbf{m}_1}(L_{0,\alpha})=\mathbb{P}_{\mathbf{m}_2}(L_{0,\alpha})$
as partial orders. 

\textbf{Third case ($\alpha=\beta+1$): }Obiously $\mathbb{P}_{\mathbf{m}_1}(L_{0,\alpha})=\mathbb{P}_{\mathbf{m}_2}(L_{0,\alpha})$
as sets. Suppose that $\mathbb{P}_{\mathbf{m}_1}(L_{0,\alpha}) \models p\leq q$.
If $s_{\beta} \notin Dom(q)$, then $p,q\in \mathbb{P}_{\mathbf{m}_1}(L_{0,\beta})$
and the claim follows from the induction hypothesis. If $s_{\beta} \in Dom(p) \cap Dom(q)$,
then by the definition of the iteration, $\mathbb{P}_{\mathbf{m}_1}(L_{0,\beta}) \models p\restriction L_{0,\beta} \leq q\restriction L_{0,\beta}$
and $q\restriction L_{0,\beta} \Vdash_{\mathbb{P}_{\mathbf{m}_1}(L_{0,\beta})} p(s_{\beta}) \leq q(s_{\beta})$.
Now note that $fsupp(p\restriction \{s_{\beta}\}), fsupp(q\restriction \{s_{\beta}\}) \subseteq L_{0,\beta}$,
hence $p(s_{\beta})$ and $q(s_{\beta})$ are $\mathbb{P}_{\mathbf{m}_2}(L_{0,\beta})$-names.
In addition, $p\restriction L_{0,\beta}, q\restriction L_{0,\beta} \in \mathbb{P}_{\mathbf{m}_1}(L_{0,\beta})=\mathbb{P}_{\mathbf{m}_2}(L_{0,\beta})$,
therefore by the induction hypothesis $\mathbb{P}_{\mathbf{m}_2}(L_{0,\beta}) \models p\restriction L_{0,\beta \leq q\restriction L_{0,\beta}}$
and $q\restriction L_{0,\beta} \Vdash_{\mathbb{P}_{\mathbf{m}_2}(L_{0,\beta})} p(s_{\beta}) \leq q(s_{\beta})$.
Therefore $\mathbb{P}_{\mathbf{m}_2}(L_{0,\alpha}) \models p\leq q$.
The other direction is proved similarly. This concludes the proof
of the induction and claim 3 (2).

We shall now return to the original induction proof.

\textbf{Third case: }$\gamma$ is a limit ordinal.

By claim 2, $\mathbb{P}_{\mathbf m}(L_{\mathbf{m}_1}) \lessdot \mathbb{P}_{\mathbf m}$.
Apply that claim to $(\mathbf{m}_1 \restriction L_{\mathbf{m}_1,\gamma}^{dp}, \mathbf{m} \restriction L_{\mathbf m,\gamma}^{dp})$
instead of $(\mathbf{m}_1,\mathbf m)$ and get $\mathbb{P}_{\mathbf m}(L_{\mathbf{m}_1,\gamma}^{dp}) \lessdot \mathbb{P}_{\mathbf m}(L_{\mathbf m,\gamma}^{dp})$.
Note that $\mathbb{P}_{\mathbf{m}_1}(L_{\mathbf{m}_1,\gamma}^{dp})=\mathbb{P}_{\mathbf m}(L_{\mathbf{m}_1,\gamma}^{dp})$
as sets, and the definition of the order depends only on $\mathbb{P}_{\mathbf{m}_1}(L_{\mathbf{m}_1,\beta}^{dp})$
for $\beta<\gamma$, therefore by the induction hypothesis $\mathbb{P}_{\mathbf{m}_1}(L_{\mathbf{m}_1,\gamma}^{dp})=\mathbb{P}_{\mathbf m}(L_{\mathbf{m}_1,\gamma}^{dp})$.
Therefore $\mathbb{P}_{\mathbf{m}_1}(L_{\mathbf{m}_1,\gamma}^{dp}) \lessdot \mathbb{P}_{\mathbf m}(L_{\mathbf m, \gamma}^{dp})$.
$\square$

\textbf{Definition 3.9: }Let $\mathbf m \in \mathbf{M}_{\leq \lambda_2}$
and $M\subseteq M_{\mathbf m}$ such that, as always, $w_t^0 \subseteq M$
for every $t\in M$. Define $\mathbf n=\mathbf{m}(M) \in \mathbf{M}_{\leq \lambda_2}$
as follows:

1. $\mathbf{q_n}=\mathbf{q_m}$.

2. $M_{\mathbf n}=M$.

3. $E_{\mathbf n}'=\{(s,t) : s\neq t \wedge \{s,t\} \nsubseteq M\}$.

4. $\bar{v}_{\mathbf n}=\bar{v}_{\mathbf m}$.

It's easy to check that $\mathbf n$ satisfies all of the requirements
in Definition 2.2 and is equivalent to $\mathbf m$, therefore $\mathbb{P}_{\mathbf m}=\mathbb{P}_{\mathbf n}$.

\textbf{Claim 3.10: }Let $\mathbf m \in \mathbf{M}_{\leq \lambda_2}$
and $M\subseteq M_{\mathbf m}$ such that, as always, $w_t^0 \subseteq M$
for every $t\in M$.

A. If $\mathbf n:=\mathbf{m}(M) \leq \mathbf{n}_1$ then there exists $\mathbf{m}_1 \in \mathbf M$
such that $\mathbf m \leq \mathbf{m}_1$ and $\mathbf{m}_1$ is equivalent
to $\mathbf{n}_1$.

B. If $\mathbf m \in \mathbf{M}_{ec}$ then $\mathbf{m}(M)=\mathbf n \in \mathbf{M}_{ec}$.

\textbf{Proof: }A) Define $\mathbf{m}_1 \in \mathbf{M}_{ec}$ as follows: 

1. $\mathbf{q_{m_1}}:=\mathbf{q_{n_1}}$.

2. $M_{\mathbf{m}_1}:=M_{\mathbf m}$.

3. $E_{\mathbf{m}_1}':=E_{\mathbf m}' \cup \{(s,t) : sE_{\mathbf{n}_1}'t \wedge \{s,t\} \subseteq (L_{\mathbf{n}_1} \setminus L_{\mathbf n}) \cup M\}$.

We shall show that $\mathbf{m}_1 \in \mathbf M$. $E_{\mathbf{m}_1}'$ is
an equivalence relation on $L_{\mathbf{m}_1} \setminus M_{\mathbf{m}_1}$:
Suppose that $s,t,r \in L_{\mathbf{m}_1} \setminus M_{\mathbf{m}_1}$
such that $sE_{\mathbf{m}_1}'t \wedge tE_{\mathbf{m}_1}'r$. If $sE_{\mathbf m}'t \wedge tE_{\mathbf m}'r$
or $sE_{\mathbf{n}_1}'t \wedge tE_{\mathbf{n}_1}'r \wedge \{s,t,r\} \subseteq (L_{\mathbf{n}_1}\setminus L_{\mathbf n})$,
then $sE_{\mathbf{m}_1}'r$, therefore we may assume WLOG that $sE_{\mathbf m}'t \wedge tE_{\mathbf{n}_1}'r \wedge \{t,r\} \subseteq L_{\mathbf{n}_1} \setminus L_{\mathbf n}$,
but this is impossible as $sE_{\mathbf m}'t$ hence $t\in L_{\mathbf m}=L_{\mathbf n}$.
Therefore $E_{\mathbf{m}_1}'$ is a transitive relation on $L_{\mathbf{m}_1} \setminus M_{\mathbf{m}_1}$
and obviously it's an equivalence relation. Suppose now that $s,t\in L_{\mathbf{m}_1} \setminus M_{\mathbf{m}_1}$
are not $E_{\mathbf{m}_1}'$-equivalent. If $s,t \in L_{\mathbf{m}_1} \setminus L_{\mathbf n}$
then $s,t$ are not $E_{\mathbf{n}_1}'$-equivalent, therefore $s<_{\mathbf{n}_1}t$
iff there exists $r\in M_{\mathbf{n}_1}$ such that $s<_{\mathbf{n}_1}r<_{\mathbf{n}_1}t$.
Therefore $s<_{\mathbf{m}_1}t$ iff there exists $r\in M_{\mathbf{m}_1}$
such that $s<_{\mathbf{m}_1}r<_{\mathbf{m}_1}t$. Suppose that $s,t \in L_{\mathbf n} \setminus M_{\mathbf{m}_1}$,
then they're not $E_{\mathbf m}'$-equivalent, therefore $s_{\mathbf m}t$
iff there is $r\in M_{\mathbf m}$ such that $s<_{\mathbf m}r<_{\mathbf m}t$.
Therefore $s_{\mathbf{m}_1}t$ iff there exists $r\in M_{\mathbf{m}_1}$
between them. Finally, suppose WLOG that $s\in L_{\mathbf{m}_1} \setminus L_{\mathbf n} \wedge t\in L_{\mathbf n} \setminus M_{\mathbf{m}_1}$
and $s<t$. If $s$ and $t$ are not $E_{\mathbf{n}_1}$-equivalent,
then as before, $s<_{\mathbf{m}_1}t$ iff there is $r\in M_{\mathbf m}$
between them. If $sE_{\mathbf{n}_1}'t$, then $s\in t/E_{\mathbf{n}_1}'=t/E_{\mathbf n}'$,
hence $s\in L_{\mathbf n}$, contradicting the choice of $s$. This
proves that $\mathbf{m}_1$ satisfies the requirement in definition $2.2(A)(D)(2)$.
It is easy to verify that $\mathbf{m}_1$ satisfies the rest of the
requirements in definition $2.2$. For example, $2.2(A)(6):$ Let
$t\in L_{\mathbf{m}_1} \setminus M_{\mathbf{m}_1}$, if $t\in L_{\mathbf n}=L_{\mathbf m}$
then $u_{\mathbf{q_{m_1}},t}^{0}=u_{\mathbf{q_{n_1}},t}^{0}=u_{\mathbf{q_n},t}^{0}=u_{\mathbf{q_m},t}^{0} \subseteq t/E_{\mathbf m}' \subseteq t/E_{\mathbf{m}_1}'$.
Suppose that $t\in L_{\mathbf{m}_1} \setminus L_{\mathbf m}$, then $u_{\mathbf{q_{m_1}},t}^{0}=u_{\mathbf{q_{n_1}},t}^{0} \subseteq t/E_{\mathbf{n}_1}'$
hence similarly $u_{\mathbf{q_{m_1}},t}^{0} \subseteq t/E_{\mathbf{m}_1}'$.

Suppose that $t\in L_{\mathbf{m}_1}$, $u\in v_{\mathbf{m}_1,t}$ and
$u\nsubseteq M_{\mathbf{m}_1}$, then $u\in v_{\mathbf{n}_1,t}$ and $u\nsubseteq M_{\mathbf{n}_1}$,
hence there is $s\in L_{\mathbf{n}_1} \setminus M$ such that $u\subseteq s/E_{\mathbf{n}_1}'$.
There are now two possibilities:

1. $t\notin M_{\mathbf{m}_1}$. In this case, for every $t\in L_{\mathbf{m}_1} \setminus M_{\mathbf{m}_1}$,
$u\subseteq u_{\mathbf{m}_1,t}^{0} \subseteq t/E_{\mathbf{m}_1}'$.

2. $t\in M_{\mathbf{m}_1}$. Suppose that $s\notin L_{\mathbf n}$. If
there is $r\in u$ such that $r\in L_{\mathbf m} \setminus M_{\mathbf n}$,
then $s\in r/E_{\mathbf{n}_1}'=r/E_{\mathbf{n}}'$, hence $s\in L_{\mathbf n}$,
which is a contradiction. Therefore $u\cup \{s\} \subseteq (L_{\mathbf{n}_1} \setminus L_{\mathbf n}) \cup M$
hence $u\subseteq s/E_{\mathbf{m}_1}'$. Suppose that $s\in L_{\mathbf n}$,
then $u\subseteq s/E_{\mathbf{n}_1}'=s/E_{\mathbf n}' \subseteq L_{\mathbf n}$,
therefore $u\in v_{\mathbf n,t}=v_{\mathbf m,t}$, hence there is $r\in L_{\mathbf m} \setminus M_{\mathbf m}$
such that $u\subseteq r/E_{\mathbf m}'$. Therefore $u\subseteq r/E_{\mathbf{m}_1}'$.
The other requirements of definition 2.2 are easy to verify, therefore
$\mathbf{m}_1 \in \mathbf M$ and obviously $\mathbf m \leq \mathbf{m}_1$
and $\mathbf{m}_1$ is equivalent to $\mathbf{n}_1$.

B) Suppose that $\mathbf n \leq \mathbf{n}_1 \leq \mathbf{n}_2$ and let
$\mathbf m \leq \mathbf{m}_1,\mathbf{m}_2$ be as in part A) for $\mathbf{n}_1$
and $\mathbf{m}_2$. We shall prove that $\mathbf{m} \leq \mathbf{m}_1 \leq \mathbf{m}_2$.
First note that $\mathbf{q_{m_1}}=\mathbf{q_{n_1}} \leq \mathbf{q_{n_2}}=\mathbf{q_{m_2}}$
and $M_{\mathbf{m}_2}=M_{\mathbf m}=M_{\mathbf{m}_1}$. Let $t\in L_{\mathbf{m}_1} \setminus M_{\mathbf{m}_1}$
and suppose that $s\in t/E_{\mathbf{m}_1}'$. By the definition of $\mathbf{m}_1$,
if $t\in L_{\mathbf m}$ then $s\in t/E_{\mathbf m}'\subseteq t/E_{\mathbf{m}_2}'$.
If $t\in L_{\mathbf{m}_1} \setminus L_{\mathbf m}$ then $sE_{\mathbf{n}_1}'t$,
hence $sE_{\mathbf{n}_2}'t$ and it follows that $sE_{\mathbf{m}_2}'t$.
Therefore $t/E_{\mathbf{m}_1}'\subseteq t/E_{\mathbf{m}_2}'$. Suppose
now that $s\in t/E_{\mathbf{m}_2}'$. If $t\in L_{\mathbf m}$ then $s\in t/E_{\mathbf{m}_2}'=t/E_{\mathbf m}'\subseteq t/E_{\mathbf{m}_1}'$.
If $t\in L_{\mathbf{m}_1} \setminus L_{\mathbf m}$ then $sE_{\mathbf{n}_2}'t$,
hence $sE_{\mathbf{n}_1}'t$ and $sE_{\mathbf{m}_1}'t$. Therefore $t/E_{\mathbf{m}_2}'\subseteq t/E_{\mathbf{m}_1}'$.
Similarly, it's easy to verify the rest of the requirements for $"\mathbf{m}_1 \leq \mathbf{m}_2"$,
therefore $\mathbf m \leq \mathbf{m}_1 \leq \mathbf{m}_2$. Now $\mathbf m \in \mathbf{M}_{ec}$,
therefore $\mathbb{P}_{\mathbf{m}_1} \lessdot \mathbb{P}_{\mathbf{m}_2}$.
Since $\mathbf{m}_l$ is equivalent to $\mathbf{n}_l$ $(l=1,2)$, we
get $\mathbb{P}_{\mathbf{n}_1} \lessdot \mathbb{P}_{\mathbf{n}_2}$, hence
$\mathbf n \in \mathbf{M}_{ec}$ as required. $\square$

\textbf{Claim 3.11: }Let $\mathbf m \in \mathbf{M}_{\leq \lambda_2}$,
then there exists $\mathbf n \in \mathbf{M}_{ec}$ such that $\mathbf m \leq \mathbf n$
and $|L_{\mathbf n}| \leq \lambda_2$.

\textbf{Proof: }Use claim 2.19 to pick $\mathbf n \in \mathbf{M}_{\chi}$
for $\chi$ large enough, such that $\mathbf n \in \mathbf{M}_{ec}$ is
very wide and full and $\mathbf m \leq \mathbf n$. We shall try to choose
$\mathbf{m}_{\alpha} \in \mathbf M$ by induction on $\alpha<\lambda_2^+$
such that the following conditions hold:

1. $\mathbf{m}_0=\mathbf m$.

2. $(\mathbf{m}_{\beta} : \beta<\alpha)\hat (\mathbf n)$ is $\leq_{\mathbf M}$-increasing
and continuous.

3. $|L_{\mathbf{m}_{\alpha}}| \leq \lambda_2$.

4. If $\alpha=\beta+1$ then one of the following conditions holds:

A) $\mathbf{m}_{\beta}$ is not wide and $\mathbf{m}_{\alpha}$ is wide.

B) There is $t_1 \in L_{\mathbf n} \setminus M_{\mathbf n}$ and a sequence
$\bar{s}_1$ of elements of $t_1/E_{\mathbf n}''$ such that for every

$t_2 \in L_{\mathbf{m}_{\beta}} \setminus M_{\mathbf m}$ and a sequence
$\bar{s}_2$ of elements of $t_2/E_{\mathbf{m}_{\beta}}''$, $(t_2,\bar{s}_2)$
is not 1-equivalent to $(t_1,\bar{s}_1)$ in $\mathbf n$, but there
is a 1-equivalent pair $(t_2,\bar{s}_2)$ in $L_{\mathbf{m}_{\alpha}}$.

We shall later prove that since $\beth_2(\lambda_1) \leq \lambda_2$,
there exists $\alpha<\lambda_2^+$ for which we won't be able to choose
an appropriate $\mathbf{m}_{\alpha}$. If $\delta$ is a limit ordinal,
then we can we can define $\mathbf{m}_{\delta}=\underset{\gamma<\delta}{\cup} \mathbf{m}_{\gamma}$,
hence necessarily $\alpha$ has the form $\alpha=\beta+1$. We shall
prove that $\mathbf{m}_{\beta}$ is as required. First we shall prove
that the pair $(\mathbf{m}_{\beta},\mathbf n)$ satisfies the assumptions
of claim 3.8 where $(\mathbf{m}_{\beta},\mathbf n)$ here stands for $(\mathbf{m}_1,\mathbf m)$
in 3.8. Obviously, $\mathbf{m}_{
\beta} \leq \mathbf n$. Suppose that $t\in L_{\mathbf n} \setminus L_{\mathbf{m}_{\beta}}$
and $\bar{s}$ is a sequence of $<\lambda^+$ members of $t/E_{\mathbf n}''$.
Let $\mathbf{m}_{\alpha} \in \mathbf M$ be wide such that $\mathbf{m}_{\beta} \leq \mathbf{m}_{\alpha} \leq \mathbf n$,
$|L_{\mathbf{m}_{\alpha}}| \leq \lambda_2$ and $\bar{s}, t$ are from
$L_{\mathbf{m}_{\alpha}}$. As $\mathbf{m}_{\alpha}$ does not satisfy
the induction's requirements, necessarily there are $t_2 \in L_{\mathbf{m}_{\beta}} \setminus M_{\mathbf m}$
and a sequence $\bar{s}_2$ of elements of $t_2/E_{\mathbf{m}_{\beta}}''$
that are 1-equivalent to $(t_1,\bar{s}_1)$ in $\mathbf n$. If $\mathbf{m}_{\beta}$
is wide, then there exists sequence $(r_{\alpha} : \alpha<\lambda^+)$
of elements of $L_{\mathbf{m}_{\beta}} \setminus M_{\mathbf m}$ such
that $r_{\alpha}/E_{\mathbf{m}_{\beta}}''\neq r_{\gamma}/E_{\mathbf{m}_{\beta}}''$
for every $\alpha<\gamma$, and $\mathbf{m}_{\beta} \restriction (r_{\alpha}/E_{\mathbf{m}_{\beta}})$
is isomorphic to $\mathbf{m}_{\beta} \restriction (t_2/E_{\mathbf{m}_{\beta}})$
for every $\alpha<\lambda^+$. For every $\alpha<\lambda^+$, denote
that isomorphism by $f_{\alpha}$ and denote by $\bar{s}_{\alpha}'$
the image of $\bar{s}_2$ under $f_{\alpha}$. Now obviously the sequence
$((r_{\alpha},\bar{s}_{\alpha}') : \alpha<\lambda^+)$ is as required.
If $\mathbf{m}_{\beta}$ is not wide, then since $\mathbf{m}_{\alpha}$
is wide, we get a contradiction to the fact the induction terminated
at $\mathbf{m}_{\beta}$. Therefore $(\mathbf{m}_{\beta},\mathbf n)$ satisfies
the assumptions of claim 3.8. 

Now suppose that $\mathbf{m}_{\beta} \leq \mathbf{n}_1 \leq \mathbf{n}_2$.
First assume that $\mathbf{n}_2 \leq \mathbf n$ and $|L_{\mathbf{n}_2}|\leq \lambda_2$.
Suppose that $t\in L_{\mathbf n} \setminus L_{\mathbf{n}_2}$ and $\bar{s}$
is a sequence of length $\zeta<\lambda^+$ of elements of $t/E_{\mathbf n}''$.
Since $(\mathbf{m}_{\beta},\mathbf n)$ satisfies the assumptions of claim
3.8, there are $t_i \in L_{\mathbf{m}_{\beta}} \setminus M_{\mathbf{m}_{\beta}} \subseteq L_{\mathbf{n}_2} \setminus M_{\mathbf{n}_2}$
and sequences $\bar{s}_i$ from $t_i/E_{\mathbf{m}_{\beta}}''=t_i/E_{\mathbf{n}_2}''$ (for $i<\lambda^+$)
as in the assumptions of claim 3.8. By claim 3.8, $\mathbb{P}_{\mathbf{n}_2} \lessdot \mathbb{P}_{\mathbf n}$.
Similarly, $\mathbb{P}_{\mathbf{n}_1} \lessdot \mathbb{P}_{\mathbf n}$,
therefore $\mathbb{P}_{\mathbf{n}_1} \lessdot \mathbb{P}_{\mathbf{n}_2}$.

Why can we assume WLOG that $|L_{\mathbf{n}_2}| \leq \lambda_2$?

Let $\chi$ be a cardinal large enough such that $\mathbf{m}_{\beta},\mathbf{n}_1,\mathbf{n}_2,\mathbf{n} \in H(\chi)$,
and let $N$ be an elementary submodel of $(H(\chi),\in)$ such that:

1. $\mathbf{m}_{\beta},\mathbf{n}_1,\mathbf{n}_2,\mathbf{n},\mathbf m \in N$.

2. $[N]^{\leq \lambda} \subseteq N$.

3. $||N|| \leq \lambda_2$.

4. $\lambda_2+1 \subseteq N$.

Let $L'=L_{\mathbf{n}_2} \cap N$, $\mathbf{n}_2'=\mathbf{n}_2 \restriction L'$
and $\mathbf{n}_1'=\mathbf{n}_1 \restriction (L'\cap L_{\mathbf{n}_1})$.
Now we may work in $N$ and replace $(\mathbf{n}_1,\mathbf{n}_2)$ by
$(\mathbf{n}_1,\mathbf{n}_2')$, as $|L_{\mathbf{n}'_2}| \leq \lambda_2$,
we get the desired result.

Why can we assume WLOG that $\mathbf{n}_2 \leq \mathbf n$?

As $\mathbf n$ is very wide and full, for every $t\in L_{\mathbf{n}_2} \setminus M_{\mathbf{n}_2}$
there exist $|L_{\mathbf n}|$ members $t_i \in L_{\mathbf n} \setminus M_{\mathbf n}$
such that $\mathbf n \restriction (t_i/E_{\mathbf n})$ is isomorphic
to $\mathbf{n}_2 \restriction (t/E_{\mathbf{n}_2})$ over $M_{\mathbf n}$
(and remember that $|L_{\mathbf{n}_2}|\leq |L_{\mathbf n}|$). Therefore
$\mathbf{n}_2$ is isomorphic to an $\mathbf{n}_3$ that satisfies $\mathbf{n}_3 \leq \mathbf n$,
so WLOG $\mathbf{n}_2 \leq \mathbf n$.

It remains to show that there exists $\alpha<\lambda_2^+$ such that
we can't choose $\mathbf{m}_{\alpha}$ as required by the induction.
Suppose towards contradiction that for every $\alpha<\lambda_2^+$
there is $\mathbf{m}_{\alpha}$ as required, then necessarily there
exist $\lambda_2^+$ ordinals $\alpha<\lambda_2^+$ such that $\mathbf{m}_{\alpha}$
satisfies $4(B)$. Therefore, there exist $\lambda_2^+$ distinct
1-equivalence classes in $\mathbf n$. We shall prove that the number
of 1-equivalence classes in $\mathbf n$ is at most $\beth_3(\lambda_1)$,
and since $\beth_3(\lambda_1) \leq \lambda_2<\lambda_2^+$, we'll
get a contradiction.

Let $\mathbf m \in \mathbf M$. First note that the number of distinct
0-equivalence classes in $\mathbf m$ is at most $\beth_2(\lambda_1)$,
as there exist at most $\beth_1(\lambda_1)$ isomorphism types of
$\mathbf m \restriction L$ for $L$ as in the definition of 0-equivalence,
so by adding the number of possible orderings of $\mathbb{P}_{\mathbf m}(L)$,
we get the desired bound. Now given $\bar{s}_2,\bar{s}_2$ as in the
definition of 1-equivalence, denote by $C_1,C_2$ the 0-equivalence
classes of sequences of the form $\bar{s}_1 \hat \bar{s}_1', \bar{s}_2 \hat \bar{s}_2'$,
respectively, for $\bar{s}_1', \bar{s}_2'$ as in the definition of
1-equivalence. $\bar{s}_1$ is 1-equivalent to $\bar{s}_2$ iff they're
0-equivalent and $C_1=C_2$. Given $\bar{s}$ as in the definition
of 1-equivalence, if $C$ is the collection of 0-equivalence classes
of sequences of the form $\bar{s} \hat \bar{s}'$ as in the definition
of 1-equivalence, then $C$ is contained in the set of 0-equivalence
classes over $\mathbf m$, which has at most $\beth_2(\lambda_1)$ members.
Therefore, there are at most $\beth_3(\lambda_1)$ different choices
for $C$, hence there are at most $\beth_3(\lambda_1)$ distinct 1-equivalence
classes over $\mathbf m$. $\square$

\textbf{\large{}Concluding the proof of the main claim}{\large\par}

\textbf{Conclusion 3.12: }A) Suppose that

0. $\mathbf{m}_l \in \mathbf{M}_{ec}$ $(l=1,2)$ and

1. $M_l \subseteq M_{\mathbf{m}_l}$ $(l=1,2)$ (and as always we assume
that $M_l$ is closed under weak memory).

2. $\mathbf{m}_1 \restriction M_1$ is isomorphic to $\mathbf{m}_2 \restriction M_2$.

3. $|L_{\mathbf{m}_1}|,|L_{\mathbf{m}_2}|\leq \lambda_2$.

Then there exists an isomorphism from $\mathbb{P}_{\mathbf{m}_1}[M_1]$
onto $\mathbb{P}_{\mathbf{m}_2}[M_2]$.

B) Suppose that $\mathbf m \in \mathbf{M}_{\leq \lambda_2}$, $M\subseteq M_{\mathbf m}=L_{\mathbf m}$
and $\mathbf n=\mathbf m \restriction M$, then $\mathbb{P}_{\mathbf n}^{cr} \lessdot \mathbb{P}_{\mathbf m}^{cr}$.

\textbf{Proof: }A) Define $\mathbf{n}_l:=\mathbf{m}_l(M_l)$ for $l=1,2$.
By claim 3.10, $\mathbf{n}_1,\mathbf{n}_2 \in \mathbf{M}_{ec}$. $\mathbf{n}_2 \restriction M_{\mathbf{n}_1}=\mathbf{m}_1 \restriction M_1$
is isomorphic to $\mathbf{n}_2 \restriction M_{\mathbf{n}_2}=\mathbf{m}_2 \restriction M_2$,
hence by claim 2.20, $\mathbb{P}_{\mathbf{n}_1}[M_{\mathbf{n}_1}]$ is
isomorphic to $\mathbb{P}_{\mathbf{n}_2}[M_{\mathbf{n}_2}]$. Therefore,
$\mathbb{P}_{\mathbf{m}_1}[M_1]$ is isomorphic to $\mathbb{P}_{\mathbf{m}_2}[M_2]$.

B) Let $\mathbf{m}_1 \in \mathbf{M}_{ec}$ such that $\mathbf m \leq \mathbf{m}_1$
and $|L_{\mathbf{m}_1}|\leq \lambda_2$. Let $\mathbf{n}_1:=\mathbf{m}_1(M)$,
then by our previous claims, $\mathbf{n}_1 \in \mathbf{M}_{ec}$. Obviously,
$\mathbf n \leq \mathbf{n}_1$, therefore $\mathbb{P}_{\mathbf n}^{cr}=\mathbb{P}_{\mathbf{n}_1}[M]=\mathbb{P}_{\mathbf{m}_1}[M] \lessdot \mathbb{P}_{\mathbf{m}_1}[L_{\mathbf m}]=\mathbb{P}_{\mathbf m}^{cr}$.
$\square$

\textbf{Conclusion 3.13: }In conclusion 2.25 we can add: Suppose that
$U_1,U_2 \subseteq \delta_*$ are closed under weak memory, $(\alpha_i : i<otp(U_1))$
and $(\beta_j : j<otp(U_2))$ are increasing enumerations of $U_1$
and $U_2$, respectively, and $h: U_1 \rightarrow U_2$ is an isomorphism
of $\mathbf m \restriction U_1$ onto $\mathbf m \restriction U_2$, then
there exists a unique generic set $G''\subseteq \mathbb{P}_{\mathbf m}^{cr}[U_2]$
such that $\eta_{\alpha_i}=\underset{\sim}{\eta_{\beta_i}}[G'']$
for every $i<otp(U_1)$.

\textbf{Proof: }In the construction that appears in 2.24 we can take
$\mathbf m \leq \mathbf n \in \mathbf{M}_{ec}$ such that $|L_{\mathbf n}| \leq \lambda_2$.
By $2.25(G+H)$ and $3.12(B)$, it follows that there exists a generic
set $G''\subseteq \mathbb{P}_{\mathbf m}^{cr}[U_2]$ such that $\eta_{\alpha_i}=\underset{\sim}{\eta_{\beta_i}}[G'']$
for every $i<otp(U_1)$. $\square$

\textbf{\large{}4. The properties of the projection and an
addition to the proof of Claim 3.8}{\large\par}

In this section we shall rely on the results of sections 0-2, with the exception of Conclusion 2.26. The results of this section will be used in the proof of Claim 3.8.
\\
\\
\textbf{Claim 4.1: }Let $p\in \mathbb{P}_{\mathbf m}$ and denote $S_p=\{\pi_L(p) : $there
exists $t\in fsupp(p)$ such that $L=t/E_{\mathbf m} \}$, then $\Vdash_{\mathbb{P}_{\mathbf m}} "p\in \underset{\sim}{G}$
iff $S_p \subseteq \underset{\sim}{G}"$.

\textbf{Proof: }If $fsupp(p) \subseteq M_{\mathbf m}$, then for every
$t\in fsupp(p)$, $\pi_{t/E_{\mathbf m}}(p)=p$, hence $S_p=\{p\}$
and there is nothing to prove. Therefore assume that $fsupp(p)\nsubseteq M_{\mathbf m}$.
By the properties of the projection, for every $t\in fsupp(p)$, $\pi_{t/E_{\mathbf m}}(p)\leq p$,
therefore $\Vdash_{\mathbb{P}_{\mathbf m}} "p\in \underset{\sim}{G} \rightarrow S_p \subseteq \underset{\sim}{G}"$.
In the other direction, suppose that $q\Vdash_{\mathbb{P}_{\mathbf m}} "S_p \subseteq \underset{\sim}{G}"$,
it's enough to show that $q$ is compatible with $p$. Assume towards
contradiction that $p$ and $q$ are incompatible. WLOG $Dom(p)\subseteq Dom(q)$.
By the assumption, $q\Vdash_{\mathbb{P}_{\mathbf m}} "\pi_{t/E_{\mathbf m}}(p) \in \underset{\sim}{G}"$
for every $t\in fsupp(p)$ and we may assume that $tr(p(s))\subseteq tr(q(s))$
for every $s\in Dom(p)$. Since $p$ contradicts $q$, there are $s\in Dom(p) \cap Dom(q)$
and $q\restriction L_{\mathbf m,<s} \leq q_1 \in \mathbb{P}_{\mathbf m}(L_{\mathbf m,<s})$
such that $q_1 \Vdash "p(s)$ contradicts $q(s)"$. By the definition
of forcing templates, $q_1 \Vdash "tr(q(s))$ contradicts $p(s)"$.
Therefore, by the definition of forcing templates and by the definition
of the iteration, there is $\iota<\iota(p(s))$ such that $q_1 \Vdash "tr(q(s))$
contradicts $\mathbf{B}_{p(s),\iota}(...,\underset{\sim}{\eta_{t_{\zeta}}}(a_{\zeta}),...)_{\zeta \in W_{p(s),\iota}}"$.
By the definition of the iteration (definition $2.2$), there is $u\in v_s$
such that $\{t_{\zeta} : \zeta \in W_{p(s),\iota}\} \subseteq u$.
By the same definition, there is $t\in fsupp(p)$ such that $\{t_{\zeta} : \zeta \in W_{p(s),\iota}\} \subseteq t/E_{\mathbf m}$.
Therefore $q_1 \Vdash "\pi_{t/E_{\mathbf m}}(p) \notin \underset{\sim}{G}$
or $tr(q(s))\nsubseteq \underset{\sim}{\eta_s}"$. Now define $q_2=q_1 \cup (q\restriction (L_{\mathbf m} \setminus L_{\mathbf m,<s}))$.
$q\leq q_2$, hence $q_2 \Vdash "\pi_{t/E_{\mathbf m}}(p) \in \underset{\sim}{G}"$.
On the other hand, $q(s)=q_2(s)$, hence $q_2 \Vdash tr(q(s)) \subseteq \underset{\sim}{\eta_s}$.
$q_1 \leq q_2$, therefore, every generic set $G$ that contains $q_2$
contains $q_1$ and also $tr(q(s))\subseteq \underset{\sim}{\eta_s}[G]$
and $\pi_{t/E_{\mathbf m}}(p)\in G$, contradicting our observation
about $q_1$. Therefore, $p$ and $q$ are compatible. $\square$

\textbf{Claim 4.2: }Let $\mathbf m \in \mathbf M$ be wide and suppose
that

1. $i(*)<\lambda$.

2. $t_i \in L_{\mathbf m} \setminus M_{\mathbf m}$ for every $i<i(*)$.

3. $t_i$ is not $E_{\mathbf m}''$-equivalent to $t_j$ for every $i<j<i(*)$.

4. $X_i=t_i/E_{\mathbf m}$.

5. $\psi_* \in \mathbb{P}_{\mathbf m}[M_{\mathbf m}]$.

6. $\psi_i \in \mathbb{P}_{\mathbf m}[X_i]$ for $i<i(*)$.

7. If $\mathbb{P}_{\mathbf m}[M_{\mathbf m}] \models \psi_* \leq \phi$,
then $\phi$ is compatible with $\psi_i$ in $\mathbb{P}_{\mathbf m}[L_{\mathbf m}]$
for every $i<i(*)$.

\textbf{then} there exists a common upper bound for $\{\psi_i : i<i(*)\} \cup \{\psi_*\}$
in $\mathbb{P}_{\mathbf m}[L_{\mathbf m}]$.

\textbf{Proof: }In this proof we shall use the notion of $*$-projection
that appears in the next section, as well as the results established independently there (it should be emphasized that this is not the same notion as the previously mentioned projection). Let $p \in \mathbb{P}_{\mathbf m}$
such that $p\Vdash_{\mathbb{P}_{\mathbf m}} "\psi_*[\underset{\sim}{G}]=true"$.
Since $\mathbf m$ is wide, there is an automorphism $f$ of $\mathbf m$
(over $M_{\mathbf m}$) that maps the members of $fsupp(p) \setminus M_{\mathbf m}$
to a set that is disjoint to $\underset{i<i(*)}{\cup}X_i$ (recall
that $|fsupp(p)|<\lambda^+$). Therefore, we may assume WLOG that
$fsupp(p) \cap X_i \subseteq M_{\mathbf m}$ for every $i<i(*)$. By
induction on $i\leq i(*)$ we'll choose conditions $p_i$ such that:

1. $p_i \in \mathbb{P}_{\mathbf m}$.

2. $(p_j : j\leq i)$ is increasing.

3. $p_0=p$.

4. If $i=j+1$ then $p_i \Vdash_{\mathbb{P}_{\mathbf m}} "\psi_j[\underset{\sim}{G}]=true"$.

5. $fsupp(p_i)$ is disjoint to $\cup \{X_j \setminus M_{\mathbf m} : i\leq j<i(*)\}$.

6. $p_i$ is chosen by the winning strategy $\mathbf{st}$ that is guaranteed
by the $(<\lambda)$-strategic completeness of $\mathbb{P}_{\mathbf m}$.

If we succeed to construct the above sequence, then for every $i<i(*)$,
$p_{i(*)} \Vdash_{\mathbb{P}_{\mathbf m}} "\psi_i[\underset{\sim}{G}]=true"$.
In addition, $p_{i(*)} \Vdash_{\mathbb{P}_{\mathbf m}} "\psi_*[\underset{\sim}{G}]=true"$
(recalling that $p\leq p_{i(*)}$), therefore, $p_{i(*)} \Vdash_{\mathbb{P}_{\mathbf m}} "\psi_*[\underset{\sim}{G}]=true \wedge (\underset{i<i(*)}{\wedge} \psi_i[\underset{\sim}{G}]=true)"$.
Therefore, $\psi_* \wedge (\underset{i<i(*)}{\wedge} \psi_i) \in \mathbb{P}_{\mathbf m}[L_{\mathbf m}]$
is the desired common upper bound.

We shall now carry the induction:

\textbf{First stage ($i=0$): }Choose $p_0=p$ (note that (5) holds
by the assumption on $fsupp(p)$).

\textbf{Second stage ($i$ is a limit ordinal): }Let $p_i'$ be an
upper bound to $(p_j : j<i)$ that is chosen according to $\mathbf{st}$.
Since $\mathbf m$ is wide, as before we can find an automorphism $f$
of $\mathbf m$ such that $f(fsupp(p_i') \setminus M_{\mathbf m})$ is
disjoint to $\cup \{X_j \setminus M_{\mathbf m} : i\leq j<i(*)\}$ and
$f$ is the identity on $\underset{j<i}{\cup}fsupp(p_j)$ (this is
possible by (5) in the induction hypothesis). Let $p_i:=\hat{f}(p_i')$.
By the definition of $\hat{f}$, $p_i$ satisfies requirements 1-5,
and as $\mathbf{st}$ is preserved by $\hat{f}$, $p_i$ satsifies (6)
as well.

\textbf{Third stage ($i=j+1$): }Let $\phi_j \in \mathbb{P}_{\mathbf m}[M_{\mathbf m}]$
be the $*$-projection of $p_j$ to $\mathbb{P}_{\mathbf m}[M_{\mathbf m}]$.
We shall first prove that $\psi_* \leq \phi_j$. If it's not true,
then there exists $\phi_j \leq \theta \in \mathbb{P}_{\mathbf m}[L_{\mathbf m}]$
contradicting $\psi_*$. Let $r\in \mathbb{P}_{\mathbf m}$ such that
$r\Vdash_{\mathbb{P}_{\mathbf m}} "\theta[\underset{\sim}{G}]=true"$,
then $r\Vdash_{\mathbb{P}_{\mathbf m}} "\psi_*[\underset{\sim}{G}]=false"$.
Since $r\Vdash_{\mathbb{P}_{\mathbf m}} "\theta[\underset{\sim}{G}]=true"$,
it follows that $\phi_j \leq \theta \leq r$, hence by the definition
of $\phi_j$, $r$ is compatible with $p_j$. By the density of $\mathbb{P}_{\mathbf m}$
in $\mathbb{P}_{\mathbf m}[L_{\mathbf m}]$, $r$ and $p_j$ have a common
upper bound $p\in \mathbb{P}_{\mathbf m}$. $p_0 \leq p_j \leq p$,
hence $p\Vdash_{\mathbb{P}_{\mathbf m}} "\psi_*[\underset{\sim}{G}]=true"$,
which is a contradiction. Therefore, $\psi_* \leq \phi_j$, hence
$\phi_j$ is compatible with $\psi_j$. By the density of $\mathbb{P}_{\mathbf m}$,
they have a common upper bound $q_j^1 \in \mathbb{P}_{\mathbf m}$.
As before, since $\mathbf m$ is wide, we may assume WLOG that $fsupp(q_j^1) \setminus M_{\mathbf m}$
is disjoint to $fsupp(p_j)\setminus M_{\mathbf m}$ and $\cup \{X_{j'} : j+1\leq j'<i(*)\}$.
By claim 4.4 (with $(p_j,q_j^1,\phi_j)$ here standing for $(p,q,\psi)$
there), $p_j$ and $q_j^1$ are compatible in $\mathbb{P}_{\mathbf m}$.
Let $p_i$ be a common upper bound chosen by the strategy. By our
choice, $\psi_j \leq p_i$, hence $p_i \Vdash_{\mathbb{P}_{\mathbf m}} "\psi_j[\underset{\sim}{G}]=true"$.
As before, use thee fact that $\mathbf m$ is wide to assume WLOG that
$fsupp(p_i) \setminus M_{\mathbf m} \cap X_{j'}=\emptyset$ for every
$i\leq j'<i(*)$. As in the previous case, we conclude that $p_i$
is as required. $\square$

\textbf{Claim 4.3: }Suppose that $\mathbf m \in \mathbf M$ is wide. Let
$f\in \mathcal{F}_{\mathbf m,\beta}$ (see definition 3.7) and denote
its domain and range by $L_1$ and $L_2$, respectively, then $f$
induces an isomorphism from $\mathbb{P}_{\mathbf m}(L_1)$ onto $\mathbb{P}_{\mathbf m}(L_2)$.

\textbf{Proof: }Obvivously, $\hat{f}$ is bijective. Now let $p_1,q_1 \in \mathbb{P}_{\mathbf m}(L_1)$
and let $p_2=\hat{f}(p_1),q_2=\hat{f}(q_1) \in \mathbb{P}_{\mathbf m}(L_2)$.
We shall prove that $\mathbb{P}_{\mathbf m} \models p_1 \leq q_1$ iff
$\mathbb{P}_{\mathbf m} \models p_2 \leq q_2$. Let $(t_i^1 : i<i(*))$
be a sequence such that:

1.$t_i^1 \in fsupp(q_1) \setminus M_{\mathbf m}$ for every $i$.

2. $t_i^1$ and $t_j^1$ are not $E_{\mathbf m}''$-equivalent for every
$i<j<i(*)$.

3. Every $t\in fsuppp(q_1) \setminus M_{\mathbf m}$ is $E_{\mathbf m}''$-equivalent
to some $t_i^1$.

For every $i<i(*)$, define $t_i^2=f(t_i^1)$ and let $\bar{t}_l=(t_i^l : i<i(*))$
$(l=1,2)$. Assume WLOG that $fsupp(p_1) \subseteq \cup \{t_i^1/E_{\mathbf m}'' : i<j(*)\} \cup M_{\mathbf m}$
for some $j(*) \leq i(*)$. For every $i<i(*)$, let $q_{1,i}=\pi_{t_i^1/E_{\mathbf m}}(q_1)$
and let $\psi_{1,i}^* \in \mathbb{P}_{\mathbf m}[M_{\mathbf m}]$ be the
$*$-projection of $q_{1,i}$ to $\mathbb{P}_{\mathbf m}[M_{\mathbf m}]$ (in
the sense of section 5). Let $\psi_1^*=\underset{i<i(*)}{\wedge} \psi_{1,i}^*$.
By the properties of the ($*$-)projection, $\psi_{1,i}^* \leq q_{1,i} \leq q_1$
for every $i<i(*)$, therefore $q_1 \Vdash_{\mathbb{P}_{\mathbf m}} "\psi_1^*[\underset{\sim}{G}]=true"$
and $\psi_1^* \in \mathbb{P}_{\mathbf m}[L_{\mathbf m}]$. For every $i<i(*)$
define $\psi_{1,i}^{**}=\psi_{1,i}^* \wedge q_{1,i} \in \mathbb{P}_{\mathbf m}[t_i^1/E_{\mathbf m}]$.
When the above conditions hold, we say that $\psi_1^*$ and $\bar{\psi}_1^*=(\psi_{1,i}^*, \psi_{1,i}^{**},q_{1,i} : i<i(*))$
analyze $q_1$ (or $(q_1,\bar{t}_1)$). Now similarly choose $\phi_1^*$
and $\bar{\phi}_1^*=(\phi_{1,i}^*, \phi_{1,i}^{**},p_{1,i} : i<j(*))$
that analyze $(p_1,(t_i^1 : i<j(*)))$. The function $f$ naturally
induces a function on $\mathbb{P}_{\mathbf m}[L_1]$, which we shall
also denote by $\hat{f}$. Now define: $\psi_2^*=\hat{f}(\psi_1^*)$,
$\psi_{2,i}^*=\hat{f}(\psi_{1,i}^*)$, $\psi_{2,i}^{**}=\hat{f}(\psi_{1,i}^{**})$,
$\phi_2^*=\hat{f}(\phi_1^*)$, $\phi_{2,i}^*=\hat{f}(\phi_{1,i}^*)$,
$\phi_{2,i}^{**}=\hat{f}(\phi_{1,i}^{**})$, $p_{2,i}=\hat{f}(p_{1,i})$,
$q_{2,i}=\hat{f}(q_{1,i})$.

It's easy to see that $(\psi_2, \bar{\psi}_2^*)$ analyze $q_2$ and
$(\phi_2^*, \bar{\phi}_2^*)$ analyze $p_2$.

Claim: Let $A_l$ $(l=1,2)$ be the claim $\mathbb{P}_{\mathbf m} \models p_l \leq q_l$
and let $B_l$ $(l=1,2)$ be the claim ``$\mathbb{P}_{\mathbf m}[t_i^l/E_{\mathbf m}] \models \phi_l^* \wedge p_{l,i} \leq \psi_l^* \wedge q_{l,i}$
for every $i<i(*)$'', then for $l\in \{1,2\}$, $A_l$ is equivalent
to $B_l$.

Proof: Suppose that $B_l$ doesn't hold for some $i$, then there
exists $\theta \in \mathbb{P}_{\mathbf m}[t_i^l/E_{\mathbf m}]$ such
that $\mathbb{P}_{\mathbf m}[t_i^l/E_{\mathbf m}] \models \psi_l^* \wedge q_{l,i} \leq \theta$
and $\theta$ is incompatible with $\phi_l^* \wedge p_{l,i}$ in $\mathbb{P}_{\mathbf m}[t_i^l/E_{\mathbf m}]$,
hence $\theta \wedge \phi_l^* \wedge p_{l,i} \notin \mathbb{P}_{\mathbf m}[t_i^l/E_{\mathbf m}]$.
For every $j$ define $\psi_j'$ as follows: If $j=i$ define $\psi_j':= \theta$.
Otherwise, define $\psi_j'=\psi_l^* \wedge q_{l,j}$. Now let $\phi' \in \mathbb{P}_{\mathbf m}[M_{\mathbf m}]$
be the $*$-projection of $\theta$ to $\mathbb{P}_{\mathbf m}[M_{\mathbf m}]$,
so if $\phi' \leq \phi \in \mathbb{P}_{\mathbf m}[M_{\mathbf m}]$ then
$\phi$ is compatible with $\theta$. Note also that $\psi_l^* \leq \phi'$:
If it wasn't true, then for some $\phi' \leq \chi \in \mathbb{P}_{\mathbf m}[M_{\mathbf m}]$,
$\chi$ contradicts $\psi_l^*$. By the choice of $\phi'$, $\chi$
is compatible with $\theta$ in $\mathbb{P}_{\mathbf m}[L_{\mathbf m}]$.
Let $\chi'$ be a common upper bound, then $\psi_l^* \leq \theta \leq \chi'$,
hence $\chi$ is compatible with $\psi_l^*$, which is a contradiction.
Therefore, $\psi_l^* \leq \phi'$.

For every $j\neq i$, if $\phi' \leq \phi \in \mathbb{P}_{\mathbf m}[M_{\mathbf m}]$,
then $\psi_{l,j}^* \leq \psi_l^* \leq \phi' \leq \phi$, hence $\phi$
is compatible with $q_{l,j}$. Since $\psi_l^* \leq \phi$, $\phi$
is also compatible with $\psi_l^* \wedge q_{l,j}$. By claim 4.2,
there is a common upper bound $q_l^+$ for $\phi'$ and all of the
$\psi_j'$. By the density of $\mathbb{P}_{\mathbf m}$, we may assume
that $q_l^+ \in \mathbb{P}_{\mathbf m}$. As $q_{l,j} \leq q_l^+$ for
every $j$, it follows from from claim 4.1 that $q_l \leq q_l^+$.
Since $\theta \leq q_l^+$ and $\theta$ contradicts $\phi_l^* \wedge p_{l,i}$,
necessarilly $q_l^+ \Vdash_{\mathbb{P}_{\mathbf m}} "(\phi_l^* \wedge p_{l,i})[\underset{\sim}{G}]=false"$.
By the properties of the projection, $p_{l,i} \leq p_l$, and as we
saw before, $\phi_l^* \leq p_l$, hence $p_l \Vdash_{\mathbb{P}_{\mathbf m}}(\phi_l^* \wedge p_{l,i})[\underset{\sim}{G}]=true$.
Now if $G\subseteq \mathbb{P}_{\mathbf m}$ is generic such that $q_l^+ \in G$,
then $q_l \in G$ and $p_l \notin G$, therefore ``$p_l \leq q_l$''
doesn't hold.

In the other direction, suppose that $B_l$ is true. Suppose towards
contradiction that $A_l$ doesn't hold. By the assumption, there is
$q_l \leq q_l^+ \in \mathbb{P}_{\mathbf m}$ contradicting $p_l$. For
$\psi_l^*$ and $\bar{\psi}_l^*$ that analyze $q_l$ we have $\mathbb{P}_{\mathbf m}[L_{\mathbf m}] \models \psi_l^* \wedge q_{l,i} \leq q_l \leq q_l^+$
for every $i$. By $B_l$, $\mathbb{P}_{\mathbf m}[L_{\mathbf m}] \models \phi_l^* \wedge p_{l,i} \leq q_l^+$
for every $i$. By claim 4.1, $p_l \leq q_l^+$, contradicting the
choice of $q_l^+$.

Therefore, $A_l$ $(l=1,2)$ is equivalent to $B_l$ $(l=1,2)$. Obviously,
$B_1$ is equivalent to $B_2$, therefore, $A_1$ is equivalent to
$A_2$. $\square$

\textbf{Claim 4.4: }Let $p,q \in \mathbb{P}_{\mathbf m}$, then $p$
and $q$ are compatible in $\mathbb{P}_{\mathbf m}$ if there exists
$\psi$ such that the following conditions hold (we shall denote this
collection of statements by $\square_{p,q,\psi}$):

1. $\psi \in \mathbb{P}_{\mathbf m}[M_{\mathbf m}]$.

2. $fsupp(p) \cap fsupp(q) \subseteq M_{\mathbf m}$, and for every
$t\in fsupp(q) \setminus M_{\mathbf m}$ and $s\in fsupp(p) \setminus M_{\mathbf m}$,
$s/E_{\mathbf m}''\neq t/E_{\mathbf m}''$.

3. If $\psi \leq \phi \in \mathbb{P}_{\mathbf m}[M_{\mathbf m}]$, then
$\phi$ is compatible with $p$ in $\mathbb{P}_{\mathbf m}[L_{\mathbf m}]$.

4. $q$ and $\psi$ are compatible in $\mathbb{P}_{\mathbf m}[L_{\mathbf m}]$.

\textbf{Proof: }We choose $(p_n, q_,n, \psi_n)$ by induction on $n<\omega$
such that the following conditions hold:

1. If $n$ is even then $\square_{p_n,q_n,\psi_n}$ holds.

2. If $n$ is odd then $\square_{q_n,p_n,\psi_n}$ holds.

3. $(p_0,q_0,\psi_0)=(p,q,\psi)$.

4. If $n=2m+1$ and $s\in Dom(p_{2m}) \cap M_{\mathbf m}$ then $s\in Dom(q_{2m+1})$
and $tr(p_{2m}(s))\subseteq tr(q_{sm+1}(s))$.

5. If $n=2m+2$ and $s\in Dom(q_{2m+1}) \cap M_{\mathbf m}$ then $s\in Dom(p_{2m+2})$
and $tr(q_{2m+1}(s))\subseteq tr(p_{2m+2}(s))$.

6. If $m<n$ then $p_m \leq p_n$ and $q_m \leq q_n$.

For $n=0$ there is no probem. Suppose that $n=2m+1$ and $(p_{2m},q_{2m}, \psi_{2m})$
has been chosen. Let $u_{2m}=Dom(p_{2m}) \cap M_{\mathbf m}$ and for
every $s\in u_{2m}$, let $\nu_s=tr(p_{2m}(s))$ and denote by $p_{s,\nu_s} \in \mathbb{P}_{\mathbf m}$
the condition $\underset{a\in Dom(\nu_s)}{\wedge}p_{s,a,\nu_s(a)}$.
Obviously, $\mathbb{P}_{\mathbf m}[L_{\mathbf m}] \models p_{s,\nu_s} \leq p_{2m}$.
Let $s\in u_{2m}$ and suppose towards contradiction that $p_{s,\nu_s}\leq \psi_{2m}$
doesn't hold, then $\psi_{2m}$ is compatible with $\neg p_{s,\nu_s}$.
Let $\phi$ be a common upper bound in $\mathbb{P}_{\mathbf m}[M_{\mathbf m}]$.
By the induction hypothesis and $\square_{p_{2m},q_{2m},\psi_{2m}}$,
$\phi$ is compatible with $p_{2m}$. Therefore, $p_{2m}$ is compatible
with $\neg p_{s,\nu_s}$, contradicting the fact that $\mathbb{P}_{\mathbf m}[L_{\mathbf m}] \models p_{s,\nu_s} \leq p_{2m}$.
Therefore, $p_{s,\nu_s} \leq \psi_{2m}$.

By the induction hypothesis and condition (4) of $\square_{p_{2m},q_{2m},\psi_{2m}}$,
there is a common upper bound $q_{2m}'$ for $q_{2m}$ and $\psi_{2m}$,
and by the density of $\mathbb{P}_{\mathbf m}$, we may suppose that
$q_{2m}'\in \mathbb{P}_{\mathbf m}$. For every $s\in u_{2m}$, since
$p_{s,\nu_s} \leq \psi_{2m}$, it follows that $\nu_s \subseteq tr(q_{2m}')$
and $s\in Dom(q_{2m}')$. Let $\psi_{2m}' \in \mathbb{P}_{\mathbf m}[M_{\mathbf m}]$
be the $*$-projection of $q_{2m}'$ to $\mathbb{P}_{\mathbf m}[M_{\mathbf m}]$.
So if $\psi_{2m}'\leq \phi \in \mathbb{P}_{\mathbf m}[M_{\mathbf m}]$,
then $\phi$ and $q_{2m}'$ are compatible in $\mathbb{P}_{\mathbf m}[L_{\mathbf m}]$.
Note also that $\psi_{2m} \leq \psi_{2m}'$: Otherwise, there is $\psi_{2m}'\leq \phi \in \mathbb{P}_{\mathbf m}[M_{\mathbf m}]$
contradicting $\psi_{2m}$. Let $\chi \in \mathbb{P}_{\mathbf m}[L_{\mathbf m}]$
be a common upper bound for $q_{2m}'$ and $\phi$, so $\psi_{2m} \leq \chi$,
therefore $\phi$ is compatible with $\psi_{2m}$, which is a contradiction.
Therefore, $\psi_{2m} \leq \psi_{2m}'$, so $p_{s,\nu_s} \leq \psi_{2m} \leq \psi_{2m}'$
for every $s\in u_{2m}$.

Since $\mathbf m$ is wide, we may assume WLOG that $fsupp(q_{2m}') \cap fsupp(p_{2m}) \subseteq M_{\mathbf m}$
and similarly for the second part of condition (2). By the induction
hypothesis and $\square_{p_{2m},q_{2m},\psi_{2m}}$, since $\psi_{2m} \leq \psi_{2m}'$,
there is a common upper bound $p_{2m}'\in \mathbb{P}_{\mathbf m}$ for
$p_{2m}$ and $\psi_{2m}'$. Since $fsupp(q_{2m}') \cap fsuppp(p_{2m}) \subseteq M_{\mathbf m}$
and $\mathbf m$ is wide, WLOG $fsupp(p_{2m}') \cap fsupp(q_{2m}') \subseteq M_{\mathbf m}$
and similarly with the second part of condition (2). Now define $p_n=p_{2m}'$,
$q_n=q_{2m}'$, $\psi_n=\psi_{2m}'$. Obviously $\square_{q_n,p_n,\psi_n}$
holds, $p_{2m} \leq p_{2m+1}$ and $q_{2m} \leq q_{2m+1}$. If $s\in Dom(p_{2m}) \cap M_{\mathbf m}$,
then $s\in Dom(q_{2m}')=Dom(q_n)$ and $tr(p_{2m}(s))=\nu_s \subseteq tr(q_{2m}'(s))=tr(q_n(s))$.
This completes the induction step for odd stages. If $n=2m+2$, the
proof is the same, alternating the roles of the $p$'s and the $q$'s.
Now choose $p_*$ and $q_*$ as the upper bounds of $(p_n : n<\omega)$
and $(q_n : n<\omega)$, repsectively, such that:

1. $Dom(p_*)=\underset{n<\omega}{\cup}Dom(p_n)$.

2. $Dom(q_*)=\underset{n<\omega}{\cup}Dom(q_n)$.

3. If $s\in Dom(p_n)$ then $tr(p_*(s))=\underset{n\leq k}{\cup}tr(p_k(s))$.

4. If $s\in Dom(q_n)$ then $tr(q_*(s))=\underset{n\leq k}{\cup}tr(q_k(s))$.

Claim: $p_*,q_* \in \mathbb{P}_{\mathbf m}$ satisfy the following conditions:

1. $Dom(p_*) \cap Dom(q_*) \subseteq M_{\mathbf m}$.

2. $Dom(p_*) \cap M_{\mathbf m}=Dom(q_*) \cap M_{\mathbf m}$.

3. If $s\in Dom(p) \cap M_{\mathbf m}$ then $tr(p_*(s))=tr(q_*(s))$
(so $p_*$ and $q_*$ are strongly compatible).

Proof: 1. Since $(p_n : n<\omega)$ and $(q_n : n\omega)$ are increasing,
then so are $(Dom(p_n) : n<\omega)$ and $(Dom(q_n) : n<\omega)$.
Since $fsupp(p_n) \cap fsupp(q_n) \subseteq M_{\mathbf m}$, it follows
that $Dom(p_*) \cap Dom(q_*) \subseteq M_{\mathbf m}$.

2. If $t\in Dom(p_*) \subseteq M_{\mathbf m}$, then $t\in Dom(p_n)$
for some even $n$. By the inductive construction, $t\in Dom(q_{n+1}) \subseteq Dom(q_*)$,
therefore $Dom(p_*) \cap M_{\mathbf m} \subseteq Dom(q_*) \cap M_{\mathbf m}$,
and the other direction is proved similarly.

3. Suppose that $s\in Dom(p_*) \cap M_{\mathbf m}$, then by the previous
claim, $s\in Dom(p_*) \cap Dom(q_*)$. Let $n<\omega$ such that $s\in Dom(p_n) \cap Dom(q_n)$,
then $tr(p_*(s))=\underset{n\leq k}{\cup}tr(p_k(s))$ and $tr(q_*(s))=\underset{n\leq k}{\cup}tr(q_k(s))$.
By conditions 4+5 of the induction, it follows that $tr(p_*(s))=tr(q_*(s))$.

By the above claim, $p_*$ and $q_*$ are compatible in $\mathbb{P}_{\mathbf m}$.
As $p=p_0 \leq p_*$ and $q=q_0 \leq q_*$, it follows that $p$ and
$q$ are compatible in $\mathbb{P}_{\mathbf m}$ as well. $\square$

\textbf{\large{}5. The existence of $*$-projections for $\mathbb{P}_{\mathbf m}[L]$}{\large\par}

\textbf{Remark: }1. The results of this section are used in the proofs of 4.2-4.4.
\\
\\
2. Note again that the notion of projection to be introduced
in the next definition is not the same as the one previously used (hence the distinction between "$*$-projection" and "projection").

\textbf{Definition 5.1: }Let $\phi \in \mathbb{P}_{\mathbf m}[L_{\mathbf m}]$.
$\psi \in \mathbb{P}_{\mathbf m}[L]$ will be called the $*$-projection
of $\phi$ to $\mathbb{P}_{\mathbf m}[L]$ if the following conditions
hold:

1. If $\mathbb{P}_{\mathbf m}[L] \models \psi \leq \theta$, then $\theta$
and $\phi$ are compatible in $\mathbb{P}_{\mathbf m}[L_{\mathbf m}]$.

2. If $\psi^* \in \mathbb{P}_{\mathbf m}[L]$ satisfies (1), then $\mathbb{P}_{\mathbf m}[L] \models \psi \leq \psi^*$.

\textbf{Claim 5.2: }Let $L\subseteq L_{\mathbf m}$. For every $\phi \in \mathbb{P}_{\mathbf m}[L]$
there exists $\psi \in \mathbb{P}_{\mathbf m}[L]$ which is the $*$-projection
of $\phi$.

\textbf{Proof: }Given $\psi_1,\psi_2 \in\mathbb{P}_{\mathbf m}[L]$,
obviously they're compatible in $\mathbb{P}_{\mathbf m}[L]$ iff they're
compatible in $\mathbb{P}_{\mathbf m}[L_{\mathbf m}]$. Let $\Lambda_1$
be the set of $\psi \in \mathbb{P}_{\mathbf m}[L]$ that contradict
$\phi$ and let $\Lambda_2$ be the set of $\psi \in \mathbb{P}_{\mathbf m}[L]$
such that $\psi$ contradicts all members of $\Lambda_1$. Let $\psi \in \mathbb{P}_{\mathbf m}[L]$.
If $\psi$ is compatible with some $\psi_1 \in \Lambda_1$, let $\psi_2$
be a common upper bound, so $\psi_2 \in \Lambda_1$. If $\psi$ contradicts
all members of $\Lambda_1$, then $\psi \in \Lambda_2$, so $\Lambda_1 \cup \Lambda_2$
is dense in $\mathbb{P}_{\mathbf m}[L]$. Note that if $\psi_1 \in \Lambda_1$
and $\psi_2 \in \Lambda_2$, then $\psi_1$ contradicts $\psi_2$.
Let $\{\psi_i : i<i(*)\}$ be a maximal antichain of elements of $\Lambda_2$.
By $\lambda^+-c.c.$, $i(*)<\lambda^+$. Define $\psi_*=\neg (\underset{i<i(*)} \wedge \neg \psi_i) \in \mathbb{P}_{\mathbf m}[L]$.
We shall prove that $\psi_*$ is a $*$-projection as desired. Suppose that
$\psi_* \leq \theta \in \mathbb{P}_{\mathbf m}[L]$ and suppose towards
contradiction that $\theta$ is incompatible with $\phi$, then $\theta \in \Lambda_1$.
Let $G\subseteq \mathbb{P}_{\mathbf m}$ be a generic set such that
$\theta[G]=true$, then for some $i$, $\psi_i[G]=true$, hence $\psi_i$
and $\theta$ are compatible. Now recall that $\psi_i \in \Lambda_2$
and $\theta \in \Lambda_1$, so we got a contradiction. Therefore
$\psi_*$ satisfies the requirement in $(1)$.

Suppose now that $\chi \in \mathbb{P}_{\mathbf m}[L]$ satisfies part
$(1)$ in Definition 5.1. Suppose towards contradiction that $\psi_* \leq \chi$
does not hold, then for some $\chi \leq \chi_*$, $\chi_*$ contradicts
$\psi_*$. Since $\Lambda_1 \cup \Lambda_2$ is dense in $\mathbb{P}_{\mathbf m}[L]$,
there is $\theta \in \Lambda_1 \cup \Lambda_2$ such that $\chi_* \leq \theta$.
Since $\chi \leq \theta$, necessarily $\theta \in \Lambda_2$. Therefore,
for some $i<i(*)$, $\theta$ is compatible with $\psi_i$, hence
this $\psi_i$ is compatible with $\chi_*$. Recall that $\psi_* \leq \psi_i$,
hence $\chi_*$ and $\psi_*$ are compatible, contradicting the choice
of $\chi_*$. Therefore, $\psi_* \leq \chi$.

\textbf{Observation 5.3: }If $\psi_1,\psi_2 \in \mathbb{P}_{\mathbf m}[L]$
are $*$-projections of $\phi \in \mathbb{P}_{\mathbf m}[L_{\mathbf m}]$,
then $\mathbb{P}_{\mathbf m}[L] \models \psi_1 \leq \psi_2 \wedge \psi_2 \leq \psi_1$.
$\square$

\textbf{Observation 5.4: }If $\psi \in \mathbb{P}_{\mathbf m}[L]$ is
the $*$-projection of $\phi \in \mathbb{P}_{\mathbf m}[L_{\mathbf m}]$, then
$\psi \leq \phi$. $\square$

\textbf{\large{}References}\textbf{ }

{[}HwSh:1067{]} Haim Horowitz and Saharon Shelah, Saccharinity with ccc, arXiv:1610.02706

{[}JuSh:292{]} Haim Judah and Saharon Shelah, Souslin forcing. J.
Symbolic Logic, 53(4), 1188--1207

{[}KeSh:872{]} Jakob Kellner and Saharon Shelah, Decisive creatures and large continuum, J. Symbolic Logic, 74(1), 73-104

{[}RoSh:628{]} Andrzej Roslanowski and Saharon Shelah, Norms on possibilities II: More ccc ideals on $2^{\omega}$, J. Appl. Anal., 3(1), 103-127

{[}Sh:587{]} Saharon Shelah, Not collapsing cardinals $\leq \kappa$
in $(<\kappa)$-support iterations, Israel Journal of Mathematics
\textbf{136} (2003), 29-115

{[}Sh:630{]} Saharon Shelah, Properness without elementaricity. J.
Appl. Anal., 10(2), 169--289

{[}Sh:945{]} Saharon Shelah, On $Con(\mathfrak{d}_{\lambda}>cov_{\lambda}(meagre))$,
Trans. Amer. Math. Soc., 373(8), 5351--5369

{[}Sh:1036{]} Saharon Shelah, Forcing axioms for $\lambda$-complete $\mu^+$-cc, Math. Log. Q., 68(1), 6-26 

{[}Sh:1126{]} Saharon Shelah, Corrected iterations, arXiv:2108.03672

$\\$

(Haim Horowitz)

E-mail address: haim.horowitz@mail.huji.ac.il

$\\$

(Saharon Shelah) Einstein Institute of Mathematics 

Edmond J. Safra Campus, 

The Hebrew University of Jerusalem.

Givat Ram, Jerusalem 91904, Israel.

Department of Mathematics 

Hill Center - Busch Campus, 

Rutgers, The State University of New Jersey. 

110 Frelinghuysen Road, Piscataway, NJ 08854-8019 USA 

E-mail address: shelah@math.huji.ac.il

\end{document}